# Gaussian Approximation of General Nonparametric Posterior Distributions

Zuofeng Shang*

*Department of Mathematical Sciences, IUPUI, Indianapolis, IN 46202, USA*

*Corresponding author: shangzf@iu.edu

and

Guang Cheng

*Department of Statistics, Purdue University, West Lafayette, IN 47907, USA*

*chengg@purdue.edu*

**In memory of Professor Jayanta K. Ghosh**

**Abstract:** In a general class of Bayesian nonparametric models, we prove that the posterior distribution can be asymptotically approximated by a Gaussian process. Our results apply to nonparametric exponential family that contains both Gaussian and non-Gaussian regression, and also hold for both efficient (root-$n$) and inefficient (non root-$n$) estimation. Our general approximation theorem does not rely on posterior conjugacy, and can be verified in a class of Gaussian process priors that has a smoothing spline interpretation [59, 44]. In particular, the limiting posterior measure becomes prior-free under a Bayesian version of "under-smoothing" condition. Finally, we apply our approximation theorem to examine the asymptotic frequentist properties of Bayesian procedures such as credible regions and credible intervals.

## 1. Introduction

A common practice in quantifying Bayesian uncertainty is to construct credible regions that cover a large fraction of posterior mass. In some cases, it is of interest to investigate the probability that the true parameter (that generates observations) is covered by the credible regions, the so-called frequentist property. Such frequentist studies on Bayesian procedures often rely on the asymptotic shapes of posterior distributions, which may be characterized by the Bernstein-von Mises (BvM) theorem.

In nonparametric settings, Freedman [14, 15] found that "almost all" Bayesian prior distributions yield inconsistent posteriors. After three decades, Cox [10] and Freedman [16] found that credible regions for nonparametric function cover the truth with probability approaching to zero. In the decades since their seminal work, achievements were made mostly in Gaussian settings. For instance, BvM theorem has been established for mean sequences (or signals) in Gaussian sequence models (equivalently, Gaussian white noise models); see [26, 22, 23, 50, 51, 7, 8, 31, 30]. In Gaussian regression with fixed design, [47, 48] proposed adaptive credible regions for regression functions; with random design, [61] proposed credible sets under sieved priors. In models where efficient estimation (at $\sqrt{n}$-rate) is possible, [41, 7, 8, 9] proposed credible intervals for functionals of infinite-dimensional parameters. As far as we are aware, posterior approximation in a general nonparametric context remains an open problem, despite the importance of non-Gaussian data.

The major goal of this paper is to prove a Gaussian approximation result in more general Bayesian nonparametric settings without relying on posterior conjugacy. Specifically, we consider a nonparametric exponential family that covers both Gaussian and non-Gaussian regression. As far as we know, even for the special Gaussian setup, the random design case was not investigated in the literature. Also, our framework is applicable even when the efficient estimation is unavailable. As explained later, our result is established based on substantially different techniques from those in the aforementioned literature.

Under total variation distance, we prove Gaussian process approximation of general posterior distributions, which significantly generalizes the (total variation) BvM result obtained by Leahu [26] in the special Gaussian white noise model. This posterior approximation result is useful in studying the frequentist properties of finite sample (or asymptotic) valid credible regions for regression functions. For instance, the frequentist coverage of the credible regions is proven to approach one given any credibility level, and can be further corrected to the credibility level by invoking a weaker topology (inspired by [7, 8]). We point out that, different from the bounded Lipschitz metric [7, 8], our approximation results hold under (stronger) total variation metric that also applies to $L^2$ credible balls and point-wise credible intervals. Our result can be viewed as complementary to [33, 36, 34, 35] who showed that although Bayesian methods are robust with finite information, they could be brittle when handling continuous systems. Rather, our positive results rely on the facts that the statistical models in consideration are correctly specified and the assigned priors charge the function space (with proper topological and geometrical details) with full mass.

Our general approximation theorem can be verified in a class of Gaussian process priors that implicitly controls the magnitude of higher-order derivatives of regression functions through a (non-random) hyper-parameter. Also, in the special Gaussian regression, these Gaussian processes match with the sequence priors considered in [59]. This leads to an interesting smoothing spline interpretation [59], which can be rigorously justified by an application of Hájek's Lemma ([19]). More importantly, this allows us to develop new technical tools based on recent progress in smoothing spline inferences (e.g.,[44]). For the above reasons, this class of Gaussian process priors can be viewed as "tuning prior." We will suggest a practical means for selecting priors via generalized cross validation (GCV). Simulation results in Section 6 strongly support this proposal. As mentioned by one referee, empirical Bayes approach for determining priors has been considered by [42] which is different from our GCV approach.

A somewhat surprising discovery in our paper is that the hyper-parameter affects the limiting posterior measure in a very subtle manner. Explicitly, we find that prior information persists in the Gaussian approximation measure under the (nearly) optimal choice of hyper-parameter. By optimal choice of hyper-parameter, we mean the one that leads to optimal contraction rate. Nevertheless, when the hyper-parameter is sub-optimal, the Gaussian approximation measure becomes prior free. This is consistent with the folklore that "data wash out prior effect" in the parametric models; see BvM theorem in [54]. The sub-optimal choice of hyper-parameter can be viewed as a Bayesian analog of the "under-smoothing" condition in the nonparametric literature.

The rest of this article is organized as follows. In Section 2, we present a general nonparametric exponential framework covering Gaussian regression and non-Gaussian regression. Section 3 includes the main results of the article. Specifically, Section 3.1 presents a general formulation of nonparametric posterior distribution. Under this formulation, Section 3.2 derives a Gaussian approximation theorem, and Section 3.3 constructs a class of Gaussian process priors for this theorem. Section 4 presents a series of applications of our main theorem that include credible region of the regression function and credible interval of a general class of linear functionals. Frequentist validity is also investigated in this section. Section 5 develops a prior-free approximation theorem and relevant inferential methods. Section 6 includes a simulation study. All proof details are postponed to a Supplementary Document [45]. A set of contraction rate results is included in the latter, and may be of independent interest.

## 2. Nonparametric Exponential Family

In this section, we present a general class of nonparametric regression models beyond Gaussian regression. Let $Y \in \mathcal{Y} \subseteq \mathbb{R}$ be response variable and $X \in \mathbb{I} := [0, 1]$ be covariate variable. Our model lies in an (natural) exponential family where given a functional parameter $f$, the random pair $(Y, X)$ follows:

$$
\begin{aligned}
p_f(y, x) &= p_f(y|x)\pi(x) \\
&= \exp\{yf(x) - A(f(x)) + c(y, x)\}\pi(x),
\end{aligned}
\tag{2.1}
$$

where $A(\cdot)$ is a known function defined upon $\mathbb{R}$, $c(y, x)$ is a quantity depending on $y, x$ to make (2.1) a valid density, and $\pi(x)$ represents marginal density of $X$. For technical convenience, we assume $\underline{\pi} \le \inf_{x \in \mathbb{I}} \pi(x) \le \sup_{x \in \mathbb{I}} \pi(x) \le \bar{\pi}$, for constants $\underline{\pi}, \bar{\pi} > 0$. The above framework (2.1) covers many non-Gaussian models; see Examples 2.1–2.4.

Assume that there exists a "true" parameter $f_0$ under which the sample is drawn from (2.1), and that $f_0$ belongs to an $m$-th order Sobolev space:

$$
S^m(\mathbb{I}) = \{f \in L^2(\mathbb{I})|f^{(j)} \text{are abs. cont. for } j = 0, 1, \ldots, m-1, \text{ and } f^{(m)} \in L^2(\mathbb{I})\}.
$$

Throughout the paper, we let $m \ge 1$ such that $S^m(\mathbb{I})$ is a RKHS; see [37].

The primary model assumption is given below. Let $\dot{A}, \ddot{A}, \dddot{A}$ be the first-, second- and third-order derivatives of $A$. Denote $\|f\|_\infty$ as the sup-norm of $f$. For any fixed $C > 0$, define $\mathcal{F}(C) = \{f \in S^m(\mathbb{I}) : \|f\|_\infty \le C\}$. Let $P_f^n$ denote the probability of the data under $f$, and $E_f$ is the expectation under $f$.

**Assumption A1.** *$A$ is three-times continuously differentiable on $\mathbb{R}$. For any $z \in \mathbb{R}$, $\ddot{A}(z) > 0$. Moreover, for any constant $C > \|f_0\|_\infty$, there exist positive constants $C_0, C_1, C_2$ (possibly depending on $C$) such that*

$$\sup_{f \in \mathcal{F}(C)} E_f \left\{ \exp(|Y - \dot{A}(f(X))|/C_0) \Big| X \right\} \le C_1, \ a.s., \tag{2.2}$$

*and for any $z \in [-2C, 2C]$,*

$$1/C_2 \le \ddot{A}(z) \le C_2, \ \text{and} \ |\dddot{A}(z)| \le C_2. \tag{2.3}$$

Assumption A1 can be easily verified in various settings including the following examples.

**Example 2.1** (*Normal regression*)**.** Suppose that under $f$, $(Y, X)$ follows normal regression:

$$Y = f(X) + \epsilon,$$

where $\epsilon \sim N(0, 1)$. Then $A(z) = z^2/2$. For any $f \in S^m(\mathbb{I})$,

$$E_f \left\{ \exp(|Y - \dot{A}(f(X))|) \Big| X \right\} = E\{\exp(|\epsilon|)\} = \frac{2}{\sqrt{e}}(1 - \Phi(1)),$$

where $\Phi(\cdot)$ is the cumulative distribution function of $\epsilon$. Therefore, (2.2) holds for $C_0 = 1$ and $C_1 = \frac{2}{\sqrt{e}}(1 - \Phi(1))$. It is easy to see that (2.3) holds for $C_2 = 1$.

**Example 2.2** (*Logistic regression*)**.** Suppose that under $f$, $(Y, X)$ follows logistic regression:

$$p_f(y|x) = \frac{\exp(yf(x))}{1 + \exp(f(x))}, \ \text{ for } y = 0, 1.$$

Here, $A(z) = \log(1 + \exp(z))$. For any $C > \|f\|_\infty$ and $f \in \mathcal{F}(C)$, $|\dot{A}(f(X))| \le (1 + \exp(-C))^{-1}$, leading to that

$$\sup_{f \in \mathcal{F}(C)} E_f \left\{ \exp(|Y - \dot{A}(f(X))|) \Big| X \right\} \le \exp\left( \frac{2 + \exp(-C)}{1 + \exp(-C)} \right),$$

and for any $z \in [-2C, 2C]$,

$$\frac{\exp(2C)}{(1 + \exp(2C))^2} \le \ddot{A}(z) \le \frac{1}{4}, \ \text{ and } \ |\dddot{A}(z)| \le \frac{1}{4},$$

which means that (2.2) holds for $C_0 = 1$ and $C_1 = \exp\left( \frac{2 + \exp(-C)}{1 + \exp(-C)} \right)$, (2.3) holds for $C_2 = \max\{\frac{1}{4}, (1 + \exp(2C))^2 \exp(-2C)\}$.

**Example 2.3** (*Binomial regression*). Suppose that under $f$, $(Y, X)$ follows binomial regression:

$$p_f(y|x) = \binom{\mathfrak{a}}{y} \frac{\exp(yf(x))}{(1 + \exp(f(x)))^{\mathfrak{a}}}, \quad \text{for } y = 0, 1, \ldots, \mathfrak{a},$$

where $\mathfrak{a}$ is a known positive integer. In particular, $\mathfrak{a} = 1$ reduces to logistic regression in Example 2.2. It is easy to see that $A(z) = \mathfrak{a} \log(1 + \exp(z))$. Similar to Example 2.2, it can be shown that (2.2) holds for $C_0 = 1$ and $C_1 = \exp\left(\frac{\mathfrak{a}+1+\exp(-C)}{1+\exp(-C)}\right)$, (2.3) holds for $C_2 = \max\{\frac{\mathfrak{a}}{4}, (1 + \exp(2C))^2 \mathfrak{a}^{-1} \exp(-2C)\}$.

**Example 2.4** (*Poisson regression*). Suppose that under $f$, $(Y, X)$ follows Poisson regression:

$$p_f(y|x) = \frac{\exp(yf(x))}{y!} \exp(-\exp(f(x))), \quad \text{for } y = 0, 1, 2, \ldots$$

Therefore, $A(z) = \exp(z)$. For any $C > \|f_0\|_\infty$ and $f \in \mathcal{F}(C)$,

$$
\begin{aligned}
E_f\left\{\exp(|Y - \dot{A}(f(X))|) \Big| X\right\} &\leq \exp(\exp(C)) E_f\left\{\exp(Y) \Big| X\right\} \\
&= \exp(\exp(C)) \times \exp(\exp(C)(e-1)) = \exp(\exp(C)e),
\end{aligned}
$$

and for any $z \in [-2C, 2C]$, $\exp(-2C) \leq \ddot{A}(z) \leq \exp(2C)$ and $|\dddot{A}(z)| \leq \exp(2C)$, implying that (2.2) holds for $C_0 = 1$ and $C_1 = \exp(\exp(C)e)$, (2.3) holds for $C_2 = \exp(2C)$.

**Remark 2.1.** *With stronger assumptions (e.g., stronger smoothness condition on $f$) and more tedious technical arguments, the results in this paper can be generalized to the following model:*

$$p_f(y|x) \sim \exp(yA_1(f(x)) - A_2(f(x)) + c(y, x)),$$

*where $A_1, A_2$ are known functions.*

Under the model Assumption A1, there exists an underlying eigen-system $(\varphi_\nu(\cdot), \rho_\nu)$ that simultaneously diagonalizes two bilinear forms $V$ and $U$, where

$$V(g, \tilde{g}) := E\{\ddot{A}(f_0(X))g(X)\tilde{g}(X)\} \quad \text{and} \quad U(g, \tilde{g}) := \int_0^1 g^{(m)}(x)\tilde{g}^{(m)}(x)dx \qquad (2.4)$$

for any $g, \tilde{g} \in S^m(\mathbb{I})$, where the expectation in the definition of $V$ is taken with respect to $\pi$, the design density. For simplicity, denote $V(g) = V(g, g)$ and $U(g) = U(g, g)$ from now on. It follows by Proposition 2.2 of [44] that $(\varphi_\nu, \rho_\nu)$ is a solution of the following ordinary differential system (whose existence and uniequeness is guaranteed by [3]):

$$
\begin{aligned}
(-1)^m \varphi_\nu^{(2m)}(\cdot) &= \rho_\nu \ddot{A}(f_0(\cdot))\pi(\cdot)\varphi_\nu(\cdot), \\
\varphi_\nu^{(j)}(0) &= \varphi_\nu^{(j)}(1) = 0, \quad j = m, m+1, \ldots, 2m-1.
\end{aligned}
\qquad (2.5)
$$

This eigen-system is building blocks of Gaussian process priors considered in this paper.

The following proposition summarizes some useful properties of $(\varphi_\nu(\cdot), \rho_\nu)$. Its proof can be found in [44, Proposition 2.2]. Two positive sequences $a_\nu, b_\nu$ are asymptotically equivalent, denoted $a_\nu \asymp b_\nu$, if $a_\nu/b_\nu$ is bounded below from zero and above from infinity. Define an inner product on $S^m(\mathbb{I})$: $\langle g, \tilde{g} \rangle_{U,V} = U(g, \tilde{g}) + V(g, \tilde{g})$. Let $\| \cdot \|_{U,V}$ be the corresponding norm, i.e., $\|g\|_{U,V} = \sqrt{\langle g, g \rangle_{U,V}}$.

**Proposition 2.1.** *Let Assumption A1 be satisfied. Then there exist a nondecreasing sequence $\rho_\nu$ and a sequence of functions $\varphi_\nu \in S^m(\mathbb{I})$ such that $\rho_1 = \cdots = \rho_m = 0$, $\rho_\nu > 0$ for $\nu > m$, $\rho_\nu \asymp \nu^{2m}$ and*

$$V(\varphi_\mu, \varphi_\nu) = \delta_{\mu\nu}, \quad U(\varphi_\mu, \varphi_\nu) = \rho_\mu \delta_{\mu\nu}, \; \mu, \nu \in \mathbb{N}, \tag{2.6}$$

*where $\delta_{\mu\nu}$ is the Kronecker's delta. In particular, any $f \in S^m(\mathbb{I})$ admits a Fourier expansion $f = \sum_\nu V(f, \varphi_\nu)\varphi_\nu$ with convergence held in the $\| \cdot \|_{U,V}$-norm.*

For any $f, g \in S^m(\mathbb{I})$, define $J(f, g) = \sum_\nu \gamma_\nu V(f, \varphi_\nu) V(g, \varphi_\nu)$, where

$$\gamma_\nu = \begin{cases} 1, & \nu = 1, 2, \ldots, m, \\ \rho_\nu, & \nu > m. \end{cases} \tag{2.7}$$

Obviously, the null space of $J$ is trivial in the sense that $J(g) := J(g, g) = 0$ if and only if $g = 0$. Furthermore, $J(g) = \int_0^1 g^{(m)}(x)^2 dx$ if $V(g, \varphi_1) = \cdots = V(g, \varphi_m) = 0$.

# 3. Main Results

## 3.1. Nonparametric Posterior Distribution

In this section, we introduce a general nonparametric Bayesian framework. Generically, $f$ is assumed to follow a probability measure $\Pi_\lambda$ (possibly involving a hyper-parameter $\lambda$). The specification of $\Pi_\lambda$ can be naturally carried out through its Radon-Nikodym (RN) derivative with respect to a base measure $\Pi$. Here, we assume that $\Pi$ is any (not necessarily Gaussian) probability measure on $(S^m(\mathbb{I}), \mathcal{B})$, where $\mathcal{B}$ is the smallest $\sigma$-algebra that contains all open subsets in $(S^m(\mathbb{I}), \| \cdot \|_{U,V})$.

The posterior distribution of $f$ can be written as

$$P(B|\mathbf{D}_n) = \frac{\int_B \exp(n\ell_n(f)) d\Pi_\lambda(f)}{\int_{S^m(\mathbb{I})} \exp(n\ell_n(f)) d\Pi_\lambda(f)}, \quad B \in \mathcal{B}. \tag{3.1}$$

for any $\Pi$-measurable subset $B \subseteq S^m(\mathbb{I})$, where the log-likelihood $\ell_n(f) = \frac{1}{n} \sum_{i=1}^n [Y_i f(X_i) - A(f(X_i))]$. Here, $\mathbf{D}_n \equiv \{Z_1, \ldots, Z_n\}$ and $Z_i = (Y_i, X_i)$, $i = 1, \ldots, n$ are *iid* copies of $Z = (Y, X)$. As for nonparametric priors, we choose the RN derivative as a function of roughness penalty $J(f)$:

$$\frac{d\Pi_\lambda}{d\Pi}(f) \propto \exp\left(-\frac{n\lambda}{2} J(f)\right), \tag{3.2}$$

where $\lambda > 0$ is a hyper-parameter. We remark that the nonparametric prior (3.2) implicitly depends on $U$ defined in (2.4), which controls the growth of $f^{(m)}$.

We will discuss in Section 3.3 a class of Gaussian process priors satisfying (3.2); see Lemma 3.3. Under (3.2), we can re-write

$$P(B|\mathbf{D}_n) = \frac{\int_B \exp(n\ell_{n,\lambda}(f))d\Pi(f)}{\int_{S^m(\mathbb{I})} \exp(n\ell_{n,\lambda}(f))d\Pi(f)}, \quad B \in \mathcal{B}. \tag{3.3}$$

Here, $\ell_{n,\lambda}(f)$ represents the penalized likelihood

$$\ell_n(f) - \frac{\lambda}{2}J(f,f),$$

which is often used in the smoothing spline literature [59]. We remark that (3.1) and (3.3) hold universally irrespective of the model assumption A1.

The hyper-parameter $\lambda$ induces a new inner product for $S^m(\mathbb{I})$ defined by

$$\langle f,g \rangle = V(f,g) + \lambda J(f,g), \quad f,g \in S^m(\mathbb{I}). \tag{3.4}$$

Let $\|f\| = \sqrt{\langle f,f \rangle}$ denote the corresponding norm. Both $\langle \cdot, \cdot \rangle$ and $\| \cdot \|$ will be very useful for subsequent theoretical analysis. For any $g = \sum_{\nu=1}^\infty g_\nu \varphi_\nu \in S^m(\mathbb{I})$, it can be seen that

$$\|g\|_{U,V}^2 = \sum_{\nu=1}^\infty g_\nu^2(1+\rho_\nu), \quad \|g\|^2 = \sum_{\nu=1}^\infty g_\nu^2(1+\lambda\gamma_\nu).$$

Therefore,

$$\min\{1,\lambda\}\|g\|_{U,V}^2 \le \|g\|^2 \le (1+\lambda)\|g\|_{U,V}^2, \quad g \in S^m(\mathbb{I}).$$

Thus, we have proved the following lemma.

**Lemma 3.1.** *For any $\lambda > 0$, $\| \cdot \|_{U,V}$ and $\| \cdot \|$ are equivalent norms (in the sense of [40]) for $S^m(\mathbb{I})$.*

By Lemma 3.1, $S^m(\mathbb{I})$ has the same topology under $\| \cdot \|_{U,V}$ and $\| \cdot \|$. Therefore, $\mathcal{B}$ can also be viewed as the Borel $\sigma$-algebra in $(S^m(\mathbb{I}), \| \cdot \|)$. Moreover, it follows by [44] that the space $(S^m(\mathbb{I}), \langle \cdot, \cdot \rangle)$ is a reproducing kernel Hilbert space (RKHS), with $K(\cdot, \cdot)$ being the reproducing kernel function.

**Proposition 1.** *Under Assumption A1, for any $f \in S^m(\mathbb{I})$ and $z \in \mathbb{I}$, we have $\|f\|^2 = \sum_\nu |V(f,\varphi_\nu)|^2(1+\lambda\gamma_\nu)$ and $K_z(\cdot) \equiv K(z,\cdot) = \sum_\nu \frac{\varphi_\nu(z)}{1+\lambda\gamma_\nu}\varphi_\nu(\cdot)$.*

### 3.2. A General Approximation Result

In this section, we show that $P(\cdot|\mathbf{D}_n)$ expressed in (3.3) can be asymptotically approximated by a posterior measure, denoted as $P_0$. Furthermore, if the imposed prior $\Pi_\lambda$ is Gaussian, $P_0$ is also Gaussian as shown in Section 3.3.

We start from a prior concentration condition (Assumption A2) on $\Pi$, which will be verified by Lemma 3.4 for a class of Gaussian process priors specified in Section 3.3. Assumption A2 is typical in Bayesian nonparametric literature; see [18]. It requires suitably large prior mass on the $\varepsilon$-ball centering at the true function.

**Assumption A2.** *There exist positive constants $c_0, c_1, \psi$ such that, for any $\varepsilon \geq \lambda^{\frac{2m+\psi}{4m}}$,*

$$\Pi(\|f - f_0\| \leq \varepsilon) \geq c_1 \exp(-c_0 \varepsilon^{-\frac{2}{2m+\psi}}).$$

Our next assumption is on the smoothness of $f_0$, expressed through its Fourier coefficients, i.e., $f_0(\cdot) = \sum_{\nu=1}^{\infty} f_\nu^0 \varphi_\nu(\cdot)$:

Condition (**S**): $\sum_{\nu=1}^{\infty} |f_\nu^0|^2 \rho_\nu^{1+\frac{\psi}{2m}} < \infty$.

Heuristically, Condition (**S**) means that $f_0 \in S^{m+\frac{\psi}{2}}(\mathbb{I})$, which requires the regularity of $f_0$ to be higher than that of $S^m(\mathbb{I})$. Such a requirement is usually needed for deriving the optimal rate of contraction; see [56]. This condition is also used to quantify the remainder term of the quadratic approximation to likelihood ratio.

The following theorem says that $P(B|\mathbf{D}_n)$ can be well approximated by

$$P_0(B) = \frac{\int_B \exp(-\frac{n}{2}\|f - \widehat{f}_{n,\lambda}\|^2) d\Pi(f)}{\int_{S^m(\mathbb{I})} \exp(-\frac{n}{2}\|f - \widehat{f}_{n,\lambda}\|^2) d\Pi(f)}, \text{ for any } B \in \mathcal{B}, \tag{3.5}$$

where $\widehat{f}_{n,\lambda}$ is a smoothing spline estimate defined as

$$\widehat{f}_{n,\lambda} = \arg \max_{f \in S^m(\mathbb{I})} \ell_{n,\lambda}(f). \tag{3.6}$$

We can view $P_0$ as a posterior measure obtained by replacing the penalized likelihood $\ell_{n,\lambda}(f)$ in (3.3) by its quadratic approximation $-\|f - \widehat{f}_{n,\lambda}\|^2/2$. The validity of this quadratic approximation is guaranteed by Assumption A2.

Let $h = \lambda^{1/(2m)}$ and $h_* \equiv n^{-\frac{1}{2m+\psi+1}}$.

**Theorem 3.2.** *(Nonparametric Posterior Approximation) Suppose prior (3.2) is imposed on $f$, Assumptions A1 and A2 hold, and $f_0 = \sum_{\nu=1}^{\infty} f_\nu^0 \varphi_\nu$ satisfies Condition (**S**). Furthermore, suppose $m > 1 + \frac{\sqrt{3}}{2} \approx 1.866$, $0 < \psi < m - \frac{1}{2}$, and $h \asymp n^{-a}$ with $a$ being a constant satisfying*

$$\max\left\{\frac{2}{6m+3\psi-1}, \frac{2m}{2m(4m+2\psi-3)+1}, \frac{1}{4m}\right\} < a \leq \frac{1}{2m+1}. \tag{3.7}$$

*Then we have, as $n \to \infty$,*

$$\sup_{B \in \mathcal{B}} |P(B|\boldsymbol{D}_n) - P_0(B)| = o_{P_{f_0}^n}(1). \tag{3.8}$$

We remark that the asymptotic posterior $P_0$ implicitly depends on the prior $\Pi$ and is not necessarily a Gaussian measure. A prior-free Gaussian approximation can be obtained under suitable choices of $h$; see Section 5.

We sketch the proof of Theorem 3.2. According to a contraction rate result (see Proposition A.1), the posterior mass is mostly concentrated on an $M\widetilde{r}_n$-ball of $f_0$, denoted as $\mathbb{B}_{M\widetilde{r}_n}(f_0)$, where $\widetilde{r}_n = (nh)^{-1/2} + h^{m+\frac{\psi}{2}}$ and $M > 0$ is a suitably large constant. Hence, for any $B \in \mathcal{B}$, we decompose $P(B|D_n) = P(B \cap \mathbb{B}_{M\widetilde{r}_n}(f_0)|D_n) + P(B \cap \mathbb{B}_{M\widetilde{r}_n}^c(f_0)|D_n)$. The second term is uniformly negligible for all $B \in \mathcal{B}$. By applying Taylor expansion to the penalized likelihood $\ell_{n,\lambda}$ (in terms of Fréchet derivatives) and empirical processes techniques, we can show that the first term is asymptotically close to $P_0(B)$ uniformly for $B \in \mathcal{B}$.

**Remark 3.1.** *It holds trivially that $h \asymp h_* := n^{-1/(2m+\psi+1)}$, among others, satisfies Condition (3.7). The choice $h \asymp h_*$ can simultaneously yield the optimal contraction rate of the credible balls as will be seen in Section 4. We also remark that Theorem 3.2 still holds by replacing Condition (3.7) with the following more general Rate Condition **(R)**:*

$$r_n = o(h^{3/2}), \ h^{1/2}\log n = o(1), \ nh^{2m+1} \geq 1, \ D_n = O(\widetilde{r}_n),$$

$$\widetilde{r}_n b_{n1} \leq 1, \ b_{n2} \leq 1, \ r_n^3 b_{n1} \leq \widetilde{r}_n^2, \ r_n^2 b_{n2} \leq \widetilde{r}_n^2, n\widetilde{r}_n^2(\widetilde{r}_n b_{n1} + b_{n2}) = o(1),$$

*where $r_n = (nh)^{-1/2} + h^m$, $\widetilde{r}_n = (nh)^{-1/2} + h^{m+\frac{\psi}{2}}$, $D_n = n^{-1/2}h^{-\frac{6m-1}{4m}}r_n\log n + h^{-1/2}r_n^2\log n$, $b_{n1} = n^{-1/2}h^{-\frac{8m-1}{4m}}(\log n)^2 + h^{-1/2}(\log n)^{3/2}$ and $b_{n2} = n^{-1/2}h^{-\frac{6m-1}{4m}}(\log n)^{3/2}$.*

**Remark 3.2.** *The TV-distance used in Theorem 3.2 is stronger than the bounded Lipschitz metric used by [7]. Hence, Theorem 3.2 can treat inferential problems with stronger topological structures, typically leading to non-root-n rate, such as the construction of $L^2$ credible region and pointwise credible interval (see Sections 4.1 and 4.3). Of course, Theorem 3.2 can also treat problems with weaker topological structures such as those in Section 4.2.*

### *3.3. Gaussian Process Prior*

In this section, we demonstrate that the general approximation Theorem 3.2 holds for the probability measures $\Pi$ and $\Pi_\lambda$ induced by a class of Gaussian process (GP) priors. In other words, we will show these $\Pi$ and $\Pi_\lambda$ satisfy (3.2). Under this class of GP priors, the limiting posterior $P_0$ is shown to be Gaussian, whose explicit characterization is given.

Let $\Pi$ be a probability measure induced by the following GP:

$$G(t) = \sum_{\nu=1}^{\infty} v_\nu \varphi_\nu(t), \tag{3.9}$$

in the sense that $\Pi(B) = P(G \in B)$ for any $B \in \mathcal{B}$. Here, $\{v_\nu\}_{\nu=1}^{\infty}$ is a sequence of independent random variables (independent of $\mathbf{D}_n$) satisfying

$$v_\nu \sim N(0, \tau_\nu^{-2}), \ \text{with} \ \tau_\nu^2 = \begin{cases} \sigma_\nu^{-2}, & \nu = 1, 2, \ldots, m, \\ \theta_\nu, & \nu > m, \end{cases} \tag{3.10}$$

$\sigma_1^2, \ldots, \sigma_m^2$ are fixed constants and $\theta_\nu \asymp \rho_\nu^{1+\beta/(2m)}$ for a constant $\beta > 1$.

We next define another GP inducing $\Pi_\lambda$:

$$G_\lambda(t) = \sum_{\nu=1}^{\infty} w_\nu \varphi_\nu(t), \tag{3.11}$$

where $w_\nu$'s are independent of the observations with

$$w_\nu \sim \begin{cases} N(0, \sigma_\nu^2/(1 + n\lambda\sigma_\nu^2)), & \nu = 1, 2, \ldots, m, \\ N(0, 1/(\theta_\nu + n\lambda\rho_\nu)), & \nu > m. \end{cases}$$

Note that $G = G_\lambda$ when $\lambda = 0$. Hence, we form a prior family $\mathcal{G} \equiv \{G_\lambda(\cdot) : \lambda \geq 0\}$. In fact, $G(\cdot)$ can be viewed as an envelope of $\mathcal{G}$ since their prior variances are the largest. The requirement $\beta > 1$ is necessary for $G_\lambda$ being a valid prior on $S^m(\mathbb{I})$. In fact, if $\beta = 1$, then the path of $G_\lambda$ does not belong to $S^m(\mathbb{I})$ almost surely (see [62, pp. 541]). However, if $\beta > 1$, we have $E\{U(G_\lambda)\} = \sum_{\nu > m} \rho_\nu/(\theta_\nu + n\lambda\rho_\nu) < \infty$, indicating that the path of $G_\lambda$ belongs to $S^m(\mathbb{I})$.

**Remark 3.3.** *Recall that the regularity of the parameter space $S^m(\mathbb{I})$ is characterized through $\rho_\nu \asymp \nu^{2m}$, in comparison with $\theta_\nu \asymp \nu^{2m+\beta}$ for that of the GP prior $G_\lambda$ for any $\lambda \geq 0$. Then, it follows from [56] that the RKHS of $G_\lambda$ is $S^{m+\frac{\beta}{2}}(\mathbb{I})$, while the parameter space $S^m(\mathbb{I})$ can be viewed as its completion in $\|\cdot\|_{U,V}$-norm. Therefore, the parameter $\beta$ represents the "relative smoothness" of the prior to the parameter space. Similar correspondence between the parameter space and prior can be found in [56, 57]. It will be seen in Section 4 that the optimal contraction rate of the posterior distribution is determined by both smoothness, i.e., the values of $m$ and $\beta$.*

**Remark 3.4.** *The GPs $G$ and $G_\lambda$ can also be written in terms of mean functions and covariance kernels. Specifically, they both have zero mean functions with covariance kernels $R(s,t) = E\{G(s)G(t)\} = \sum_{\nu \geq 1} \varphi_\nu(s)\varphi_\nu(t)/\tau_\nu^2$ and $R_\lambda(s,t) = E\{G_\lambda(s)G_\lambda(t)\} = \sum_{\nu=1}^m \varphi_\nu(s)\varphi_\nu(t)/(\sigma_\nu^{-2} + n\lambda) + \sum_{\nu > m} \varphi_\nu(s)\varphi_\nu(t)/(\theta_\nu + n\lambda\rho_\nu)$. Similar to [2], our GP prior can be also viewed as a Gaussian measure with covariance being a positive, self-adjoint and trace-class operator on the functional space. Note that the covariance might not be a Matérn kernel.*

In the following lemma, we show that (3.2) holds under the above class of GP priors.

**Lemma 3.3.** *[19, Hájek's Lemma] With $f \in S^m(\mathbb{I})$, we have the following Radon-Nikodym derivative of $\Pi_\lambda$ with respect to $\Pi$:*

$$\begin{aligned} \frac{d\Pi_\lambda}{d\Pi}(f) &= \prod_{\nu=1}^m (1 + n\lambda\sigma_\nu^2)^{1/2} \prod_{\nu=m+1}^\infty (1 + n\lambda\rho_\nu/\theta_\nu)^{1/2} \exp\left(-\frac{n\lambda}{2} J(f)\right) \\ &\propto \exp\left(-\frac{n\lambda}{2} J(f)\right). \end{aligned}$$

It should be mentioned that Wahba [59] designed a set of Gaussian sequence priors to estimate the smoothing parameter. Her prior yields a similar RN derivative as the one in Lemma 3.3.

Our next lemma shows that Assumption A2 holds under the above $\Pi$ induced by (3.9).

**Lemma 3.4.** *Assumption A2 holds for $\psi = \beta - 1$ for $\Pi$ induced by (3.9).*

Lemma 3.4 will be needed in constructing Bayesian inference procedures; see Section 4.

Given the above $\Pi$ and $\Pi_\lambda$, we next show that $P_0$ is essentially a Gaussian measure $\Pi_W$, induced by the following Gaussian process $W$, in the sense that $P_0(B) = P(W \in B|\mathbf{D}_n) := \Pi_W(B)$. Given the Fourier expansion $\widehat{f}_{n,\lambda}(\cdot) = \sum_{\nu=1}^{\infty} \widehat{f}_\nu \varphi_\nu(\cdot)$, we define a GP

$$W(\cdot) = \sum_{\nu=1}^{\infty} (a_{n,\nu} \widehat{f}_\nu + b_{n,\nu} \xi_\nu) \varphi_\nu(\cdot), \tag{3.12}$$

where $a_{n,\nu} = n(1+\lambda\gamma_\nu)\{\tau_\nu^2 + n(1+\lambda\gamma_\nu)\}^{-1}$, $b_{n,\nu} = \{\tau_\nu^2 + n(1+\lambda\gamma_\nu)\}^{-1/2}$, $\xi_\nu = \tau_\nu v_\nu$ are iid standard normal variables with $v_\nu$ and $\tau_\nu^2$ satisfying (3.10), and the sequence $\gamma_\nu$ is defined in (2.7). For better illustration, we decompose $W$ as

$$W = \widetilde{f}_{n,\lambda} + W_n, \tag{3.13}$$

where $\widetilde{f}_{n,\lambda}(\cdot) := \sum_{\nu=1}^{\infty} n(1+\lambda\gamma_\nu)\{\tau_\nu^2 + n(1+\lambda\gamma_\nu)\}^{-1} \widehat{f}_\nu \varphi_\nu(\cdot)$ and $W_n(\cdot) := \sum_{\nu=1}^{\infty} \{\tau_\nu^2 + n(1+\lambda\gamma_\nu)\}^{-1/2} \xi_\nu \varphi_\nu(\cdot)$ is a zero-mean GP. Note that the posterior mode $\widetilde{f}_{n,\lambda}$ is asymptotically equivalent to the efficient linear estimate $\widehat{f}_{n,\lambda}$ since $\|\widetilde{f}_{n,\lambda} - \widehat{f}_{n,\lambda}\| = o_{P_{f_0}^n}(1)$. This is consistent with the traditional BvM theorem in the parametric setup ([54]).

Theorem 3.5 presents the Gaussian characterization of $P_0$.

**Theorem 3.5.** *With $f \in S^m(\mathbb{I})$, the Radon-Nikodym derivative of $\Pi_W$ with respect to $\Pi$ is*

$$\frac{d\Pi_W}{d\Pi}(f) = \frac{\exp(-\frac{n}{2}\|f - \widehat{f}_{n,\lambda}\|^2)}{\int_{S^m(\mathbb{I})} \exp(-\frac{n}{2}\|f - \widehat{f}_{n,\lambda}\|^2) d\Pi(f)}.$$

*Hence, we have*

$$\frac{dP_0}{d\Pi}(f) = \frac{d\Pi_W}{d\Pi}(f).$$

Together with Theorem 3.2, Theorem 3.5 implies that the posterior distribution $P(\cdot|\mathbf{D}_n)$ and $\Pi_W(\cdot)$ are asymptotically close under the total variation distance. This approximation result greatly facilitates the construction and theoretical analysis of nonparametric Bayesian inference procedures in Section 4. For example, from (3.13), we can tell that $\widetilde{f}_{n,\lambda}$ is approximately the posterior mode of $P(\cdot|\mathbf{D}_n)$, and can be used as the center of credible region, e.g., (4.1).

## 4. Bayesian Inference Procedures

In this section, we consider Bayesian inference procedures such as credible balls and point-wise credible intervals. These inference procedures are fully driven by posterior samples, so-called finite sample construction. For example, the radius of the credible ball is directly drawn from MCMC samples so that the posterior coverage is exact. We also comment on the asymptotic construction where the radius is obtained by asymptotic theory in Remark 4.1. Under a proper choice of $\lambda$, these Bayesian inference procedures are shown to possess frequentist validity.

Throughout this section, we choose $\Pi, \Pi_\lambda$ as GPs designed in Section 3.3 for technical convenience. We also suppose that $f_0$ satisfies Condition (**S**), and let $h \asymp h_*$.

### 4.1. Credible Region in Strong Topology

We consider the construction of credible region for $f$ in terms of $L^2$-norm, and also study its frequentist property. The existing literature mostly focuses on Gaussian setup: [26, 22, 7, 23, 50, 8, 31, 51] for Gaussian white noise, [47, 48] for Gaussian regression with fixed design, and [61] for Gaussian regression with sieved priors. In contrast, our results are established in the more general nonparametric exponential family.

For any $f, g \in S^m(\mathbb{I})$, define $\langle f, g \rangle_2 = V(f, g)$, an inner product on $S^m(\mathbb{I})$, and $\|f\|_2 = V(f)^{1/2}$ its corresponding norm, which is a type of $L^2$-norm. For any $\alpha \in (0, 1)$, let $r_n(\alpha) > 0$ satisfy $P(f \in S^m(\mathbb{I}) : \|f - \widetilde{f}_{n,\lambda}\|_2 \le r_n(\alpha) | \mathbf{D}_n) = 1 - \alpha$. In practice, $r_n(\alpha)$ can be computed as the $(1 - \alpha)$ posterior quantile of $\|f - \widetilde{f}_{n,\lambda}\|_2$ through MCMC samples of $f$; see [29] for more introduction. A credible region with an exact credibility level $(1 - \alpha)$ is constructed as

$$R_n(\alpha) = \left\{ f \in S^m(\mathbb{I}) : \|f - \widetilde{f}_{n,\lambda}\|_2 \le r_n(\alpha) \right\}. \tag{4.1}$$

We next examine the frequentist property of $R_n(\alpha)$.

**Theorem 4.1.** *Suppose that Assumption A1 holds, $f_0$ satisfies Condition (**S**), $m > 1 + \frac{\sqrt{3}}{2}$, $1 < \beta < m + 1/2$, and $h \asymp h_*$. Then for any $\alpha \in (0, 1)$, $\lim_{n \to \infty} P^n_{f_0}(f_0 \in R_n(\alpha)) = 1$.*

It is easy to see that the $L^2$-diameter of $R_n(\alpha)$ is $2r_n(\alpha)$. Remark 4.1 reveals that $r_n(\alpha)$ achieves the optimal rate $n^{-\frac{2m+\beta-1}{2(2m+\beta)}}$ when $h \asymp h_*$, therefore, the $L^2$-diameter of $R_n(\alpha)$ attains optimality.

A relevant result in [22] says that the credible region of the mean sequence in Gaussian sequence models has coverage probability approaching one when the hyper-parameter is properly selected as an order of $n$. Theorem 4.1 generalizes their result to nonparametric exponential family.

### 4.2. Credible Region in Weak Topology

The frequentist coverage of the credible region (4.1) asymptotically approaches one regardless of the credibility level. This motivates us to construct a modified credible region using a *weaker* topology such that the truth can be covered with probability approaching *exactly* the credibility level. Besides Theorem 3.2, our proof also relies on a strong approximation result ([52]).

We first define a weaker metric by following [7, 8]. For any $f \in S^m(\mathbb{I})$ with $f(\cdot) = \sum_{\nu=1}^{\infty} f_\nu \varphi_\nu(\cdot)$, define $\|f\|_\omega^2 = \sum_{\nu=1}^{\infty} \omega_\nu f_\nu^2$, where $\omega_\nu$ is a given positive sequence satisfying $\omega_\nu = \nu^{-1} (\log 2\nu)^{-\tau}$ for a constant $\tau > 1$. Since $\omega_\nu \le 1$ for all $\nu \ge 1$, it is easy to see that $\|f\|_\omega \le \|f\|_2$. Therefore, $\|\cdot\|_\omega$ is weaker than $\|\cdot\|_2$. We will show that under this weaker norm, any $(1 - \alpha)$ credible region can recover exactly $(1 - \alpha)$ frequentist coverage.

For any $\alpha \in (0, 1)$, let $r_{\omega,n}(\alpha) > 0$ satisfy $P(f \in S^m(\mathbb{I}) : \|f - \widetilde{f}_{n,\lambda}\|_\omega \le r_{\omega,n}(\alpha) | \mathbf{D}_n) = 1 - \alpha$. We construct a credible region with credibility level $(1 - \alpha)$:

$$R_n^\omega(\alpha) = \left\{ f \in S^m(\mathbb{I}) : \|f - \widetilde{f}_{n,\lambda}\|_\omega \le r_{\omega,n}(\alpha) \right\}. \tag{4.2}$$

Theorem 4.2 proves that $R_n^\omega(\alpha)$ asymptotically possesses the frequentist coverage $(1 - \alpha)$.

**Theorem 4.2.** *Suppose that Assumption [A1] holds, $f_0$ satisfies Condition (**S**), $m > 1 + \frac{\sqrt{3}}{2}$, $1 < \beta < \min\{m + \frac{1}{2}, \frac{(2m-1)^2}{2m}\}$, and $h \asymp h_*$. Then for any $\alpha \in (0,1)$, $\lim_{n \to \infty} P^n_{f_0}(f_0 \in R^\omega_n(\alpha)) = 1 - \alpha$.*

[7, 8] consider credible regions with similar frequentist property in Gaussian sequence models an density estimation. Theorem 4.2 generalize their results to nonparametric exponential family.

We note that the $L^2$-diameter of $R^\omega_n(\alpha)$ is infinity (see Section A.4 of appendix). But we can impose a restriction to make its $L^2$-diameter being finite, by using a strategy of [7]. Specifically, define

$$R^{\star\omega}_n(\alpha) = R^\omega_n(\alpha) \cap \{f \in S^m(\mathbb{I}) : J(f) \le M\},$$

for a constant $M > 0$. It can be shown that the $L^2$-diameter of $R^{\star\omega}_n(\alpha)$ is $O_{P^n_{f_0}}(n^{-\frac{2m+\beta-1}{2(2m+\beta)}}\sqrt{\log n})$ (see Section A.4 of appendix). The leading factor $n^{-\frac{2m+\beta-1}{2(2m+\beta)}}$ is the optimal contraction rate under Sobolev norm (see Section A.6). So the $L^2$-diameter is now rate optimal upto a logarithmic factor.

### *4.3. Linear Functionals on the Regression Function*

We construct credible intervals for a general class of linear functionals in nonparametric exponential family. Frequentist coverage of the proposed credible interval is also investigated. In particular, we consider two important special cases: (i) evaluation functional: $F_z(f) = f(z)$, where $z \in \mathbb{I}$ is a fixed number; (ii) integral functional: $F_\omega(f) = \int_0^1 f(z)\omega(z)dz$, where $\omega(\cdot)$ is a known deterministic integrable function such as an indicator function. We find that the former leads to an interval contracting at slower than root-$n$ rate, while the latter leads to root-$n$ rate.

The existing literature mostly focus on functionals where efficient estimation with $\sqrt{n}$-rate is available ([41, 7, 8, 9]). The more general inefficient estimation with slower than root-$n$ rate (e.g., evaluation functional) is only treated recently by [49] in Gaussian white noise model. As will be seen, our theory treat efficient and inefficient estimation in a unified framework.

Let $F : S^m(\mathbb{I}) \mapsto \mathbb{R}$ be a linear $\Pi$-measurable functional, i.e., $F(af + bg) = aF(f) + bF(g)$ for any $a, b \in \mathbb{R}$ and $f, g \in S^m(\mathbb{I})$. We say that $F$ satisfies Condition (**F**) if there exist constants $\kappa > 0$ and $r \in [0,1]$ such that for any $f \in S^m(\mathbb{I})$,

$$|F(f)| \le \kappa h^{-r/2}\|f\|. \tag{4.3}$$

Lemma 4.3 below (given in [44]) implies that both $F_z$ and $F_\omega$ satisfy (4.3).

**Lemma 4.3.** *There exists a universal constant $c > 0$ s.t. for any $f \in S^m(\mathbb{I})$, $\|f\|_\infty \le ch^{-1/2}\|f\|$.*

Let $r_{F,n}(\alpha) > 0$ satisfy $P(f \in S^m(\mathbb{I}) : |F(f) - F(\widetilde{f}_{n,\lambda})| \le r_{F,n}(\alpha)|\mathbf{D}_n) = 1 - \alpha$. Define $(1 - \alpha)$ credible interval for $F(f)$ as

$$CI^F_n(\alpha) : F(\widetilde{f}_{n,\lambda}) \pm r_{F,n}(\alpha). \tag{4.4}$$

Theorem 4.4 below shows that $CI^F_n(\alpha)$ covers the true value $F(f_0)$ with probability asymptotically at least $(1 - \alpha)$ for any $F$ satisfying Condition (**F**). Our result holds for both efficient estimation such as $F = F_\omega$ and inefficient case such as $F = F_z$ in contrast with the existing literature.

For $k \geq 1$, define

$$\theta_{k,n}^2 = \sum_{\nu=1}^{\infty} \frac{F(\varphi_\nu)^2}{(\tau_\nu^2 + n(1 + \lambda\gamma_\nu))^k}.$$

**Theorem 4.4.** *Suppose that Assumption A1 holds, $f_0 = \sum_{\nu=1}^{\infty} f_\nu^0 \varphi_\nu$ satisfies Condition* (**S′**): $\sum_{\nu=1}^{\infty} |f_\nu^0|^2 \nu^{2m+\beta} < \infty$, $m > 1 + \frac{\sqrt{3}}{2}$, $1 < \beta < \min\{m + \frac{1}{2}, \frac{(2m-1)^2}{2m}\}$, *and $h \asymp h_*$. Meanwhile,*

$$n^k \theta_{k,n}^2 \asymp h^{-r} \quad \text{for } k = 1, 2. \tag{4.5}$$

*Then for any $\alpha \in (0,1)$, we have*

$$\liminf_{n \to \infty} P_{f_0}^n(F(f_0) \in CI_n^F(\alpha)) \geq 1 - \alpha, \tag{4.6}$$

*given that Condition* (**F**) *holds. Moreover, if $0 < \sum_{\nu=1}^{\infty} F(\varphi_\nu)^2 < \infty$, then $\lim_{n \to \infty} P_{f_0}^n(F(f_0) \in CI_n^F(\alpha)) = 1 - \alpha$.*

By carefully examining the proof of Theorem 4.4, we find that when $F = F_z$, the inequality (4.6) is actually strict, and the length of $CI_n^F(\alpha)$ satisfies $r_{F,n}(\alpha) \asymp n^{-\frac{2m+\beta-1}{2(2m+\beta)}}$. When $F = F_\omega$, $CI_n^F(\alpha)$ covers the truth with probability approaching $1 - \alpha$, and its length satisfies $r_{F,n}(\alpha) \asymp n^{-1/2}$. Please see Remark 4.1 for more details. Therefore, there exists a subtle difference between the two types of functionals. Simulation results in Section 6 empirically confirm this subtle distinction.

Remark that Condition (**S′**) is slightly stronger than Condition (**S**), which is used to correct certain bias arising from the prior. Such a condition can be understood as the "under-smoothing" condition in smoothing spline; see [44]. Condition (4.5) is not restrictive and can be verified in concrete settings; see Proposition 4.5 below. The proof of Proposition 4.5 relies on a nice closed form of $\varphi_\nu$ and a careful analysis of the trigonometric functions.

**Proposition 4.5.** *Suppose $m = 2$, $X \sim Unif[0,1]$, and $Y|f, X \sim N(f(X), 1)$.*

   (i) *If $F = F_z$ for any $z \in (0,1)$, then (4.5) holds for $r = 1$;*

   (ii) *If $F = F_\omega$ for any $\omega \in L^2(\mathbb{I}) \backslash \{0\}$, then $0 < \sum_{\nu=1}^{\infty} F(\varphi_\nu)^2 < \infty$ and (4.5) holds for $r = 0$.*

**Remark 4.1.** *The radii $r_n(\alpha)$, $r_{\omega,n}(\alpha)$ and $r_{F,n}(\alpha)$ are determined by posterior samples of $f$ This might be time-consuming in practice. In fact, the proofs of Theorems 4.1, 4.2 and 4.4 reveal that the radii satisfy the following large-sample (data-free) limits:*

$$
\begin{aligned}
r_n(\alpha) &= \left( \sqrt{\frac{\zeta_{1,n} + \sqrt{2\zeta_{2,n}} z_\alpha}{n}} \right) \left( 1 + o_{P_{f_0}^n}(1) \right), \quad \text{where } \zeta_{k,n} = \sum_{\nu=1}^{\infty} \frac{1}{(1 + \lambda\gamma_\nu + n^{-1}\tau_\nu^2)^k}, \\
r_{\omega,n}(\alpha) &= \sqrt{\frac{c_\alpha}{n}} \left( 1 + o_{P_{f_0}^n}(1) \right), \\
r_{F,n}(\alpha) &= \theta_{1,n} z_{\alpha/2} \left( 1 + o_{P_{f_0}^n}(1) \right),
\end{aligned}
\tag{4.7}
$$

*where $c_\alpha > 0$ satisfies $P(\sum_{\nu=1}^{\infty} \omega_\nu \eta_\nu^2 \leq c_\alpha) = 1 - \alpha$ with $\eta_\nu$ being independent standard Gaussian, and $z_\alpha = \Phi^{-1}(1 - \alpha)$ with $\Phi$ being the standard Gaussian c.d.f. Replacing the radii by the above limits (4.7), one can establish the asymptotic proxies of (4.1), (4.2) and (4.4), which can reduce computational burden. The frequentist coverage of these asymptotic regions/intervals remains the same as the original ones. Proof of Theorem 4.1 indicates that $r_n(\alpha)$ attains the optimal rate of contraction $n^{-\frac{2m+\beta-1}{2(2m+\beta)}}$ when $h \asymp h_*$ and $\zeta_{k,n} \asymp n^{1/(2m+\beta)}$ for $k = 1, 2$. The optimal contraction rate can be viewed as a Bayesian counterpart of the minimax estimation rate in classic frequentist literature, e.g., [53].*

## 5. Nonparametric BvM Theorem

The traditional BvM theorem ([54]) in parametric models requires the limit posterior measure to be prior free. However, the posterior approximation in Section 3 still contains some prior information, i.e., the sequences $\sigma_\nu^2$ and $\tau_\nu^2$ in $W$. Similar phenomenon has also been observed in Bayesian sparse linear models; see Theorem 6 in [4].

In this section, we derive a limit Gaussian posterior that is nonetheless prior free, and thus establish nonparametric BvM theorem in the traditional sense. This can be achieved by simply choosing a sub-optimal $h$, in contrast with the optimal choice of $h$ in Section 3. This finding can be viewed as a Bayesian analog of the well known "under-smoothing" idea in the nonparametric literature. Additionally, we prove that some other choices of $h$ lead to the failure of BvM. Construction of posterior balls together with their asymptotic validity are also investigated based on the new nonparametric BvM theorem.

The intuition behind our prior-free limit Gaussian measure is quite simple: we set the prior information $\tau_\nu = 0$ in the expression of $W$ given in (3.13). The resulting GP becomes

$$W^\star = \widehat{f}_{n,\lambda} + W_n^\star,$$

where $W_n^\star(\cdot) = \sum_{\nu > m} \{n(1 + \lambda \gamma_\nu)\}^{-1/2} \xi_\nu \varphi_\nu(\cdot)$ and $\xi_\nu \overset{iid}{\sim} N(0,1)$. Notably, $W^\star$ depends only on the smoothing spline estimate $\widehat{f}_{n,\lambda}$ and the sequence $\gamma_\nu$. The latter depends on differential equations (2.5) which involve only the function $A(\cdot)$ and true $f_0$. Hence, $W^\star$ contains no prior information. Define $P_\star$ as the probability measure of $W^\star$ (conditional on $\mathbf{D}_n$).

We next show that $P_\star$ is indeed an asymptotic posterior measure. Based on Theorem 3.2, it suffices to show that the deviation between $P_0$ and $P_\star$ is sufficiently small. Unfortunately, this cannot be achieved if we choose $h \asymp h_*$ or its small neighborhood, i.e., (5.3). In this case, the mean of $P_0$, i.e., $\widetilde{f}_{n,\lambda}$, and that of $P_\star$, i.e., $\widehat{f}_{n,\lambda}$, are found not to converge to each other fast enough. This leads to the failure of BvM theorem. However, if we choose $h$ converging to zero significantly slower than $h_*$ in the sense of (5.1), a prior-free nonparametric BvM theorem holds.

**Theorem 5.1.** *Suppose Assumption A1 holds and $f_0 = \sum_{\nu=1}^{\infty} f_\nu^0 \varphi_\nu$ satisfies Condition (S) with*

$\psi = \beta - 1$. Let $m > 1 + \frac{\sqrt{3}}{2} \approx 1.866$, $1 < \beta < m - \frac{1}{2}$ and $h \asymp n^{-a}$ with

$$\max\left\{\frac{2}{6m + 3\beta - 4}, \frac{2m}{2m(4m + 2\beta - 5) + 1}, \frac{1}{4m}\right\} < a < \frac{2}{4m + 2\beta + 1}. \tag{5.1}$$

*Then we have, as $n \to \infty$,*

$$\sup_{t \geq 0} |P(\|f - \widehat{f}_{n,\lambda}\|_2 \leq t|\boldsymbol{D}_n) - P_\star(\|f - \widehat{f}_{n,\lambda}\|_2 \leq t)| = o_{P^n_{f_0}}(1). \tag{5.2}$$

*Moreover, there exists a $f_0$ satisfying Condition (**S**) such that, for any $h \asymp n^{-a}$ with a satisfying*

$$\frac{2}{4m + 2\beta + 1} \leq a < \frac{8m + 4\beta + 2}{(8m + 4\beta + 1)(2m + \beta)}, \tag{5.3}$$

*(5.2) does not hold.*

Based on Theorem 5.1, we are ready to construct (prior-free) credible balls with the center $\widehat{f}_{n,\lambda}$ that asymptotically attain desirable credibility levels. Unfortunately, the corresponding radii in this case converge to zero even faster than the optimal rate of contraction such that the truth will be excluded from the credible balls. A simple remedy is to "blow up" the radius (see the similar idea of [51]). To be more specific, we construct a ball centering at $\widehat{f}_{n,\lambda}$ with radius $(1 + \varepsilon)r_n^\star$:

$$R_n^\star(\varepsilon) = \{f \in S^m(\mathbb{I}) : \|f - \widehat{f}_{n,\lambda}\|_2 \leq (1 + \varepsilon)r_n^\star\}, \text{ for any } \varepsilon > 0, \tag{5.4}$$

where

$$r_n^\star = \sqrt{\frac{1}{nh}\int_0^\infty \frac{1}{(1 + x^{2m})^2}dx + h^{2m}}.$$

The above choice of $r_n^\star$ is of the order $r_n = (nh)^{-1/2} + h^m$ which can achieve the rate $n^{-m/(2m+1)}$ by using $h = n^{-1/(2m+1)}$; see Remark 3.1. Note that we may use generalized cross validation to select such a $h$; see [58]. In practice, one may replace $r_n^\star$ by a finite-sample counterpart, e.g., $r_n(\alpha)^{\frac{2m(2m+\beta)}{(2m+\beta-1)(2m+1)}}$, which can be shown to also achieve the rate $n^{-m/(2m+1)}$ (recalling that $r_n(\alpha)$ is the radius of the ball $R_n(\alpha)$ determined in Section 4.1).

A direct consequence of Theorem 5.1 implies that $R_n^\star(\varepsilon)$ asymptotically possesses large credibility level and frequentist coverage.

**Corollary 5.2.** *Suppose Assumption A1 holds and $f_0$ satisfies Condition (**S**) with $\psi = \beta - 1$. Let $m > 1 + \frac{\sqrt{3}}{2} \approx 1.866$, $1 < \beta < m - \frac{1}{2}$ and $h \asymp n^{-a}$ with a satisfying (5.1). Then for any $\varepsilon > 0$, as $n \to \infty$, $P(R_n^\star(\varepsilon)|\boldsymbol{D}_n) \geq 1 - \alpha$ with $P^n_{f_0}$-probability approaching one, and $P^n_{f_0}(f_0 \in R_n^\star(\varepsilon)) = 1 + o(1)$.*

## 6. Simulations

In this section, we empirically investigate the frequentist coverage probabilities of the credible region (4.1) and modified credible region (4.2), and credible intervals for evaluation functional and integral functional. As for the choice of $\lambda$ (equivalently, $h$) in the GP prior, we suggest employing the generalized cross validation (GCV) method. For example, let $h_{GCV}$ be the GCV-selection of $h$, which is known to achieve an rate $n^{-1/(2m+1)}$; see [58]. Then, we can set $h_* = h_{GCV}^{(2m+1)/(2m+\beta)}$. This method works very well as demonstrated in the simulations.

We generated data from the following model

$$Y_i = f_0(X_i) + \epsilon_i, \ i = 1, 2, \ldots, n, \tag{6.1}$$

where $X_i$ are *iid* uniform over $[0, 1]$, and $\epsilon_i$ are *iid* standard normal random variables independent of $X_i$. The true regression function $f_0$ was chosen as $f_0(x) = 3\beta_{30,17}(x) + 2\beta_{3,11}(x)$, where $\beta_{a,b}$ is the probability density function for $Beta(a, b)$. Figure 1 displays the true function $f_0$, from which it can be seen that $f_0$ has both peaks and troughs. GP prior defined in Section 3.3 was used with $m = \beta = 2$. The $h$ was selected based on GCV proposed by [59]. MCMC algorithms were employed to draw posterior samples of $f$.

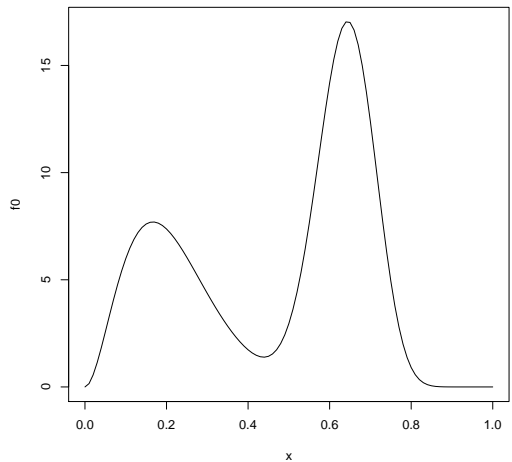

FIG 1. *Plot of the true function $f_0$ used in model (6.1).*

To examine the coverage property of the credible regions, we chose $n$ ranging from 20 to 2000. For each $n$, 1,000 independent trials were conducted. From each trial, a credible region (CR) based on (4.1) and a modified credible region (MCR) based on (4.2) were constructed. Proportions of the CR and MCR covering $f_0$ were calculated, and were displayed against the sample sizes. Results are summarized in Figure 2. It can be seen that for different $1 - \alpha$, i.e., the credibility levels, the coverage proportions (CP) of CR are greater than $1 - \alpha$ when $n$ is large enough. They even tend to one for large sample sizes. However, the CP of the MCR tends to exactly $1 - \alpha$ when $n$ increases. Thus, the numerical results confirm our theory developed in Sections 4.1 and 4.2.

To examine the coverage property of credible intervals, we chose $n = 2^5, 2^7, 2^8, 2^9$ to demonstrate the trend of coverage along with increasing sample sizes. For evaluation functional, we considered $F = F_z$ for 15 evenly-spaced $z$ points in $[0, 1]$. For each $z$, a credible interval based on (4.4) was constructed. We then calculated the coverage probability of this interval based on 1,000 independent experiments, that is, the empirical proportion of the intervals (among the 1,000 intervals) that cover the true value $f_0(z)$. Figure 3 summarizes the results for different credibility levels $\alpha$, where coverage probabilities are plotted against the corresponding points $z$. It can be seen that the coverage probability of the pointwise intervals is a bit larger than $1 - \alpha$ for all $\alpha$ and $n$ being considered. This is consistent with Proposition 4.5 (i), except for the points near the right peak of $f_0$. Indeed, at those points near the right peak, under-coverage has been observed. This is a common phenomenon in the frequentist literature: the peak and trouts may affect the coverage property of the pointwise interval; see [32, 44]. This is also possible due to the mismatch of smoothness level between the prior and true parameter.

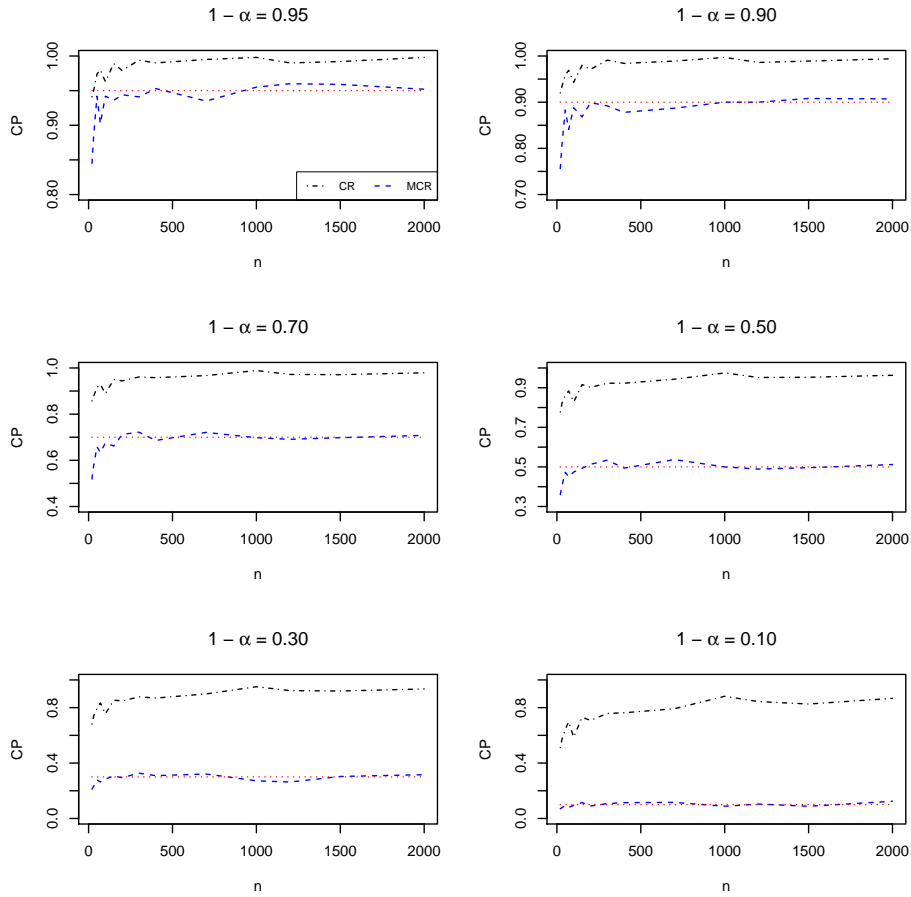

FIG 2. *Coverage proportion (CP) of CR and MCR, constructed by (4.1) and (4.2) respectively, based on various sample sizes and credibility levels. The dotted red line indicates the position of the $1 - \alpha$ credibility level.*

For integral functional, we considered $F = F_{\omega_{z_0}}$ for $\omega_{z_0}(z) = I(0 \leq z \leq z_0)$ with 15 evenly-spaced $z_0$ points in $[0, 1]$. We evaluated the coverage probability at each $z_0$ based on 1,000 exper-iments. Figure 4 summarizes the results for different credibility levels $\alpha$, where coverage probabilities are plotted against the corresponding points $z_0$. It can be seen that, as $n$ increases, the coverage probability of the integral intervals tends to $1 - \alpha$ for all $\alpha$. This phenomenon is consistent with our theory, i.e., Proposition 4.5 (ii).

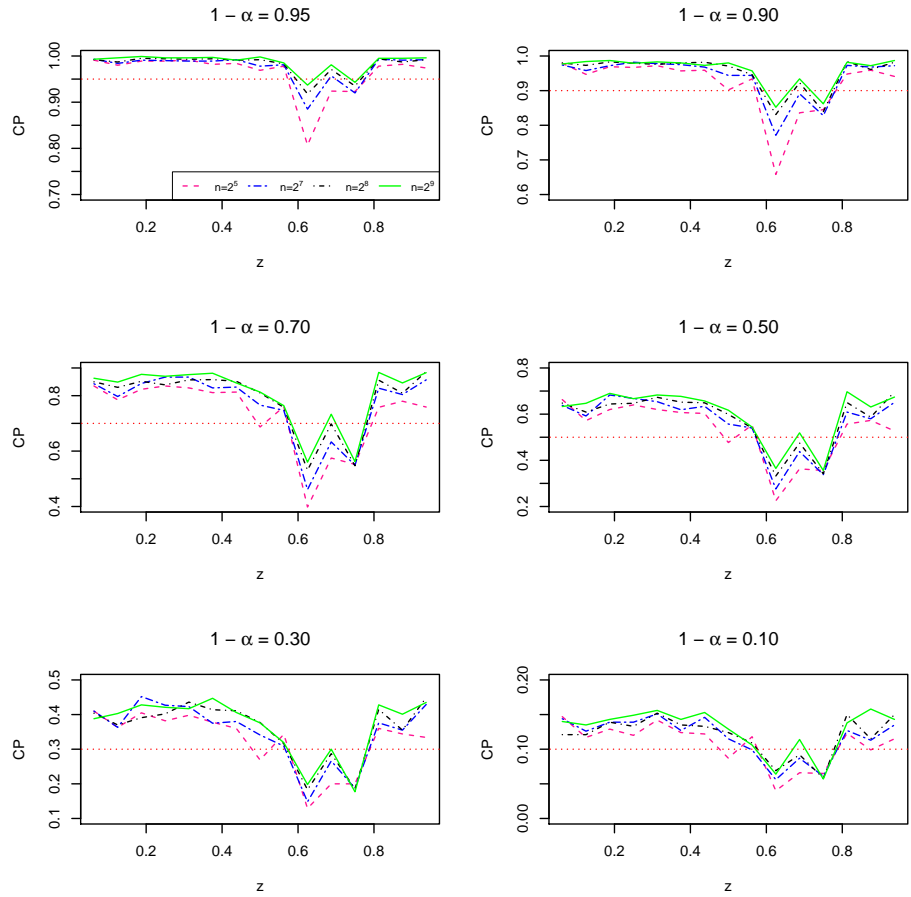

FIG 3. *Coverage proportion (CP) of the credible interval for $F_z(f_0)$ versus $z$. The dotted red line indicates the position of the $1 - \alpha$ credibility level.*

## Funding

Shang's research was sponsored by NSF DMS-1764280. Cheng's research was sponsored by NSF CAREER Award DMS-1151692, DMS-1418042, DMS-1712919, Simons Fellowship in Mathematics and Office of Naval Research (ONR N00014-15-1-2331).

**Acknowledgements**. We thank Prof. Jayanta Ghosh for careful reading and comments, and also thank PhD student Meimei Liu at Purdue for help with the simulation study.

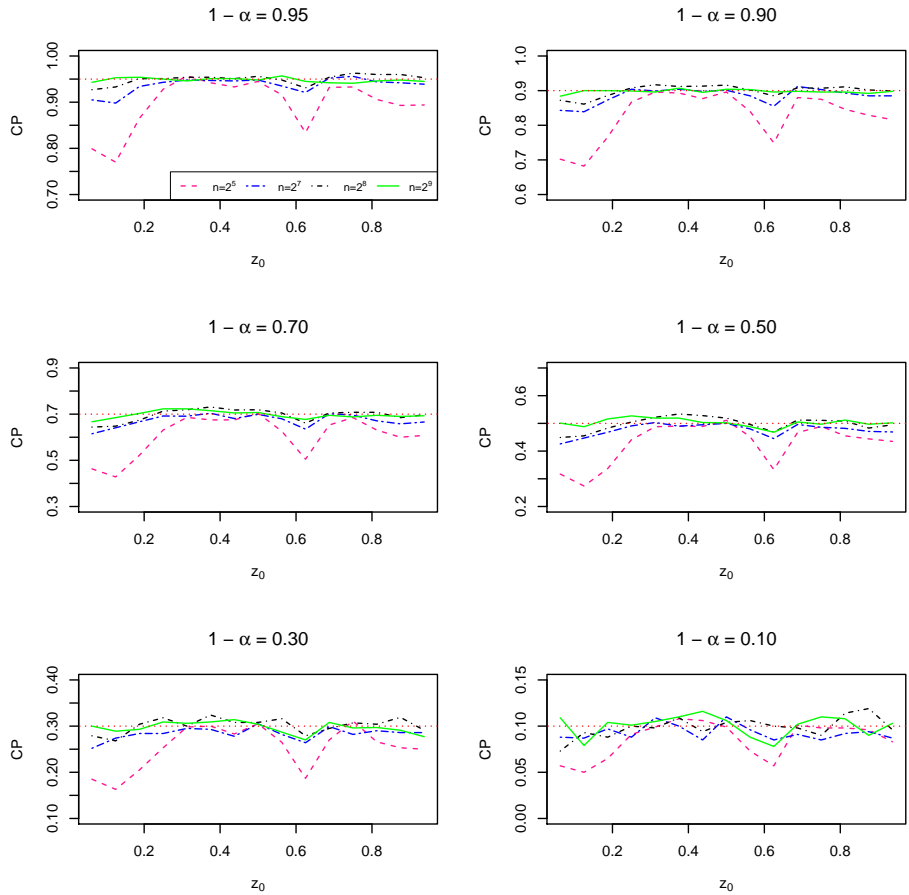

FIG 4. *Coverage proportion (CP) of the credible interval for* $F_{\omega_{z_0}}(f_0)$ *versus* $z_0$. *The dotted red line indicates the position of the* $1-\alpha$ *credibility level.*

# Gaussian Approximation of General Nonparametric Posterior Distributions

This supplementary document consists of two parts. Part I contains the proofs of the main results in this paper. Part II contains the proofs of auxiliary results.

## Supplementary Document: Part I

This section contains proofs of main results in Sections 3, 4 and 5.

### *A.1. Proofs in Section 3*

Before proving Theorem 3.2, let us state several preliminary results.

**Lemma A.1.** *Under Condition (**S**), we have $\|\widehat{f}_{n,\lambda} - f_0\| = O_{P_{f_0}}(\widetilde{r}_n)$.*

**Proposition A.1.** *(Contraction Rate) Suppose Assumption A1 holds, and $f_0 = \sum_{\nu=1}^{\infty} f_\nu^0 \varphi_\nu$ satisfies Condition (**S**). Furthermore, the following Rate Condition holds:*

$$r_n = o(h^{3/2}), \ h^{1/2} \log n = o(1), \ nh^{2m+1} \geq 1, \ D_n = O(\widetilde{r}_n),$$
$$\widetilde{r}_n b_{n1} \leq 1, \ b_{n2} \leq 1, \ r_n^3 b_{n1} \leq \widetilde{r}_n^2, \ r_n^2 b_{n2} \leq \widetilde{r}_n^2.$$

*Then, for any $\varepsilon_1, \varepsilon_2 \in (0,1)$, there exist positive constants $M', N'$ s.t. for any $n \geq N'$,*

$$P_{f_0}^n \left( P(\|f - f_0\| \geq M' \widetilde{r}_n | \boldsymbol{D}_n) \geq \varepsilon_2 \right) \leq \varepsilon_1, \tag{A.1}$$

*where $P_f^n$ denotes the probability measure induced by $\boldsymbol{D}_n$ under $f$.*

Lemmas A.1 is a direct consequence of Lemma A.12 in online supplementary and $\|S_{n,\lambda}(f_0)\| = O_{P_{f_0}^n}(h^{m+\psi/2})$, where $S_{n,\lambda}(f_0)$ is defined according to [43], i.e.,

$$S_{n,\lambda}(f_0) = \frac{1}{n} \sum_{i=1}^{n} (Y_i - \dot{A}(f_0(X_i))) K_{X_i} - \mathcal{P}_\lambda f_0.$$

**Lemma A.2.** *It holds that*

$$\ell_{n,\lambda}(f) - \ell_{n,\lambda}(\widehat{f}_{n,\lambda}) + \frac{1}{2} \|f - \widehat{f}_{n,\lambda}\|^2 = T_1(f) + T_2(f), \tag{A.2}$$

*where*

$$
\begin{aligned}
T_1(f) &= -\frac{1}{n} \int_0^1 \int_0^1 s \sum_{i=1}^n [\ddot{A}(\widehat{f}_{n,\lambda}(X_i) + ss'(f - \widehat{f}_{n,\lambda})(X_i))(f - \widehat{f}_{n,\lambda})(X_i)^2 \\
&\qquad -\ddot{A}(f_0(X_i))(f - \widehat{f}_{n,\lambda})(X_i)^2]dsds', \\
T_2(f) &= -\frac{1}{2n} \sum_{i=1}^n [\ddot{A}(f_0(X_i))(\Delta f)(X_i)^2 - E_{f_0}^X \{\ddot{A}(f_0(X))(f - \widehat{f}_{n,\lambda})(X)^2\}].
\end{aligned}
\tag{A.3}
$$

*Proof of Lemma A.2.* By Taylor expansion in terms of Fréchet derivative, the result holds. See Section A.6. $\square$

*Proof of Theorem 3.2.* It follows by Remark 3.1 that Rate Condition (**R**) holds. We will prove the theorem under Rate Condition (**R**).

Let $\varepsilon_1, \varepsilon_2$ be arbitrarily small positive constants. Let $\varepsilon_3$ be small fixed with $0 < \varepsilon_3 < \log 2$ and $4\varepsilon_3 \exp(\varepsilon_3) + 2\varepsilon_3 \le \varepsilon_2/3$. Consider three events:

$$
\begin{aligned}
\mathcal{E}'_n &= \{\|\widehat{f}_{n,\lambda} - f_0\| \le M_1 \widetilde{r}_n\} \\
\mathcal{E}''_n &= \{P(\|f - f_0\| \ge M_2 \widetilde{r}_n | \mathbf{D}_n) \le \varepsilon_3\} \\
\mathcal{E}'''_n &= \{P_0(\|f - f_0\| \ge M_2 \widetilde{r}_n) \le \varepsilon_3\},
\end{aligned}
$$

where $M_1, M_2$ are large enough constants. It follows from Lemma A.1 and Proposition A.1 that we can choose $M_2 > M_1$ (both large enough) s.t. $P_{f_0}^n(\mathcal{E}'_n \cap \mathcal{E}''_n) \ge 1 - \varepsilon_1/2$. On $\mathcal{E}'_n$, we have

$$
\|f - f_0\| - M_1 \widetilde{r}_n \le \|f - \widehat{f}_{n,\lambda}\| \le \|f - f_0\| + M_1 \widetilde{r}_n.
$$

Hence,

$$
\begin{aligned}
P_0(\|f - f_0\| \ge M_2 \widetilde{r}_n) &= \frac{\int_{\|f - f_0\| \ge M_2 \widetilde{r}_n} \exp\left(-\frac{n}{2}\|f - \widehat{f}_{n,\lambda}\|^2\right) d\Pi(f)}{\int_{S^m(\mathbb{I})} \exp\left(-\frac{n}{2}\|f - \widehat{f}_{n,\lambda}\|^2\right) d\Pi(f)} \\
&\le \frac{\int_{\|f - f_0\| \ge M_2 \widetilde{r}_n} \exp\left(-\frac{n}{2}\|f - \widehat{f}_{n,\lambda}\|^2\right) d\Pi(f)}{\int_{\|f - f_0\| \le \widetilde{r}_n} \exp\left(-\frac{n}{2}\|f - \widehat{f}_{n,\lambda}\|^2\right) d\Pi(f)} \\
&\le \exp\left(-\left((M_2 - M_1)^2/2 - (M_1 + 1)^2/2\right) n\widetilde{r}_n^2\right) \Pi(\|f - f_0\| \le \widetilde{r}_n)^{-1} \\
&\le \exp c_1^{-1}\left(-\left((M_2 - M_1)^2/2 - (M_1 + 1)^2/2 - c_0/4\right) n\widetilde{r}_n^2\right),
\end{aligned}
\tag{A.4}
$$

where the last inequality is due to Assumption A2 and the trivial fact $\widetilde{r}_n > h^{m+\frac{\psi}{2}}$ and $\widetilde{r}_n^{-\frac{2}{2m+\psi}} \le n\widetilde{r}_n^2/4$. The last inequality follows by

$$
\widetilde{r}_n \ge 2n^{-\frac{2m+\psi}{2(2m+\psi+1)}}.
$$

We can even manage $M_2$ to be large so that the quantity (A.4) is less than $\varepsilon_3$. Therefore, we get that $P_{f_0}^n(\mathcal{E}'''_n) \ge P_{f_0}^n(\mathcal{E}'_n \cap \mathcal{E}''_n) \ge 1 - \varepsilon_1/2$. Define $\mathcal{E}_n = \mathcal{E}'_n \cap \mathcal{E}''_n \cap \mathcal{E}'''_n$, then it can be seen that $P_{f_0}^n(\mathcal{E}_n) \ge 1 - \varepsilon_1$.

Using an empirical process argument (see (A.73) and (A.74) of Section A.6), it can be shown that on $\mathcal{E}_n$, for any $f \in S^m(\mathbb{I})$ satisfying $\|f - f_0\| \leq M_2 \widetilde{r}_n$,

$$|T_1(f)| \leq D_1 \times \widetilde{r}_n^3 b_{n1}, \quad |T_2(f)| \leq D_2 \times \widetilde{r}_n^2 b_{n2}, \tag{A.5}$$

where $D_1, D_2$ are constants depending on $M_1, M_2$. Since $n\widetilde{r}_n^2(\widetilde{r}_n b_{n1} + b_{n2}) = o(1)$, we choose $n$ to be large enouch so that $D_1 \times n\widetilde{r}_n^3 b_{n1} + D_2 \times n\widetilde{r}_n^2 b_{n2} \leq \varepsilon_3$.

Define

$$J_{n1} = \int_{S^m(\mathbb{I})} \exp\left(n(\ell_{n,\lambda}(f) - \ell_{n,\lambda}(\widehat{f}_{n,\lambda}))\right) d\Pi(f),$$

$$J_{n2} = \int_{S^m(\mathbb{I})} \exp\left(-\frac{n}{2}\|f - \widehat{f}_{n,\lambda}\|^2\right) d\Pi(f),$$

$$\bar{J}_{n1} = \int_{\|f-f_0\| \leq M_2 \widetilde{r}_n} \exp\left(n(\ell_{n,\lambda}(f) - \ell_{n,\lambda}(\widehat{f}_{n,\lambda}))\right) d\Pi(f),$$

$$\bar{J}_{n2} = \int_{\|f-f_0\| \leq M_2 \widetilde{r}_n} \exp\left(-\frac{n}{2}\|f - \widehat{f}_{n,\lambda}\|^2\right) d\Pi(f).$$

It is easy to see that on $\mathcal{E}_n$,

$$0 \leq \frac{J_{n1} - \bar{J}_{n1}}{J_{n1}} \leq \varepsilon_3, \quad 0 \leq \frac{J_{n2} - \bar{J}_{n2}}{J_{n2}} \leq \varepsilon_3.$$

By some algebra, it can be shown that the above inequalities lead to

$$(1 - \varepsilon_3) \cdot \frac{\bar{J}_{n2}}{\bar{J}_{n1}} \leq \frac{J_{n2}}{J_{n1}} \leq \frac{1}{1 - \varepsilon_3} \cdot \frac{\bar{J}_{n2}}{\bar{J}_{n1}}. \tag{A.6}$$

Meanwhile, on $\mathcal{E}_n$, using (A.5), Lemma A.2 and the elementary inequality $|\exp(x) - 1| \leq 2|x|$ for $|x| \leq \log 2$, we get that

$$\begin{aligned}
& |\bar{J}_{n2} - \bar{J}_{n1}| \\
\leq \ & \int_{\|f-f_0\| \leq M_2 \widetilde{r}_n} \exp\left(-\frac{n}{2}\|f - \widehat{f}_{n,\lambda}\|^2\right) \times |\exp(n(T_1(f) + T_2(f))) - 1| d\Pi(f) \\
\leq \ & 2\varepsilon_3 \bar{J}_{n2},
\end{aligned}$$

leading to that

$$\frac{1}{1 + 2\varepsilon_3} \leq \frac{\bar{J}_{n2}}{\bar{J}_{n1}} \leq \frac{1}{1 - 2\varepsilon_3}. \tag{A.7}$$

Combining (A.6) and (A.7), on $\mathcal{E}_n$,

$$\frac{1 - \varepsilon_3}{1 + 2\varepsilon_3} \leq \frac{J_{n2}}{J_{n1}} \leq \frac{1}{(1 - 2\varepsilon_3)(1 - \varepsilon_3)},$$

leading to

$$-4\varepsilon_3 \leq \frac{1 - \varepsilon_3}{1 + 2\varepsilon_3} - 1 \leq \frac{J_{n2}}{J_{n1}} - 1 \leq \frac{1}{(1 - 2\varepsilon_3)(1 - \varepsilon_3)} - 1 \leq 4\varepsilon_3 \tag{A.8}$$

For simplicity, denote $R_n(f) = n(T_1(f) + T_2(f))$. For any $S \in \mathcal{S}$, let $S' = S \cap \{f \in S^m(\mathbb{I}) : \|f - f_0\| \leq M_2 \tilde{r}_n\}$. Then on $\mathcal{E}_n$, we get that

$$|P(S|\mathbf{D}_n) - P_0(S)| \leq |P(S'|\mathbf{D}_n) - P_0(S')| + 2\varepsilon_3.$$

Moreover, it follows from Lemma A.2 and (A.8) that on $\mathcal{E}_n$,

$$\begin{aligned}
&|P(S'|\mathbf{D}_n) - P_0(S')| \\
&= \left| \int_{S'} \left( \frac{\exp(n(\ell_{n,\lambda}(f) - \ell_{n,\lambda}(\widehat{f}_{n,\lambda})))}{J_{n1}} - \frac{\exp\left(-\frac{n}{2}\|f - \widehat{f}_{n,\lambda}\|^2\right)}{J_{n2}} \right) d\Pi(f) \right| \\
&\leq \int_{S'} \exp\left(-\frac{n}{2}\|f - \widehat{f}_{n,\lambda}\|^2\right) \times \left| \frac{\exp(R_n(f))}{J_{n1}} - \frac{1}{J_{n2}} \right| d\Pi(f) \\
&\leq \int_{S'} \exp\left(-\frac{n}{2}\|f - \widehat{f}_{n,\lambda}\|^2\right) \times \frac{|\exp(R_n(f)) - 1|}{J_{n2}} d\Pi(f) \\
&\quad + \int_{S'} \exp\left(-\frac{n}{2}\|f - \widehat{f}_{n,\lambda}\|^2\right) \times \exp(R_n(f)) \times \left| \frac{1}{J_{n1}} - \frac{1}{J_{n2}} \right| d\Pi(f) \\
&\leq 2\varepsilon_3 \frac{\int_{S'} \exp\left(-\frac{n}{2}\|f - \widehat{f}_{n,\lambda}\|^2\right) d\Pi(f)}{J_{n2}} \\
&\quad + \exp(\varepsilon_3) \times \left| \frac{1}{J_{n1}} - \frac{1}{J_{n2}} \right| \times \int_{S'} \exp\left(-\frac{n}{2}\|f - \widehat{f}_{n,\lambda}\|^2\right) d\Pi(f) \\
&\leq 2\varepsilon_3 + \exp(\varepsilon_3) \times \left| \frac{J_{n2}}{J_{n1}} - 1 \right| \leq 2\varepsilon_3 + 4\varepsilon_3 \exp(\varepsilon_3) \leq \varepsilon_2/3.
\end{aligned}$$

Note that the right hand side is free of $S$. Then we get that on $\mathcal{E}_n$,

$$\sup_{S \in \mathcal{S}} |P(S|\mathbf{D}_n) - P_0(S)| \leq \varepsilon_2/3 + 2\varepsilon_3 \leq \varepsilon_2.$$

This implies that for sufficiently large $n$,

$$\begin{aligned}
&P_{f_0}^n \left( \sup_{S \in \mathcal{S}} |P(S|\mathbf{D}_n) - P_0(S)| > \varepsilon_2 \right) \\
&\leq P_{f_0}^n(\mathcal{E}_n^c) + P_{f_0}^n \left( \mathcal{E}_n, \sup_{S \in \mathcal{S}} |P(S|\mathbf{D}_n) - P_0(S)| > \varepsilon_2 \right) \\
&= P_{f_0}^n(\mathcal{E}_n^c) \leq \varepsilon_1.
\end{aligned}$$

This completes the proof. $\qquad\square$

*Proof of Lemma 3.3.* For any $f \in S^m(\mathbb{I})$, by Proposition 2.1, $f$ admits a unique series representation $f = \sum_{\nu=1}^{\infty} f_\nu \varphi_\nu$, where $f_\nu = V(f, \varphi_\nu)$ satisfies $\sum_\nu f_\nu^2 \rho_\nu < \infty$. Therefore, $T : f \mapsto \{f_\nu : \nu \geq 1\}$ defines a one-to-one map from $S^m(\mathbb{I})$ to $\mathcal{R}_m \equiv \{\{f_\nu\}_{\nu=1}^{\infty} \in \mathbb{R}^\infty : \sum_{\nu=1}^{\infty} f_\nu^2 \rho_\nu < \infty\}$.

Let $\tilde{\Pi}_\lambda$ and $\tilde{\Pi}$ be the probability measures induced by $\{w_\nu : \nu > m\}$ and $\{v_\nu : \nu > m\}$, respectively, which are both defined on $\mathbb{R}^\infty$. That is, for any subset $S \in \mathbb{R}^\infty$, $\tilde{\Pi}_\lambda(S) = P(\{w_\nu : $

$\nu > m\} \in S)$ and $\tilde{\Pi}(S) = P(\{v_\nu : \nu > m\} \in S)$. Likewise, let $\Pi'_\lambda$ and $\Pi'$ be probability measures induced by $\{w_\nu : \nu \geq 1\}$ and $\{v_\nu : \nu \geq 1\}$. It is easy to see that, for any measurable $B \subseteq \mathcal{R}_m$,

$$\Pi_\lambda(T^{-1}B) = P(G_\lambda \in T^{-1}B) = P(\{w_\nu : \nu \geq 1\} \in B) = \Pi'_\lambda(B), \text{ and}$$

$$\Pi(T^{-1}B) = P(G \in T^{-1}B) = P(\{v_\nu : \nu \geq 1\} \in B) = \Pi'(B).$$

The following result can be found in Hájek [19].

**Proposition A.2.** *The Radon-Nikodym derivative of $\tilde{\Pi}_\lambda$ w.r.t. $\tilde{\Pi}$ is*

$$\begin{aligned}
\frac{d\tilde{\Pi}_\lambda}{d\tilde{\Pi}}(\{f_\nu : \nu > m\}) &= \prod_{\nu > m}^{\infty} (1 + n\lambda\rho_\nu/\theta_\nu)^{1/2} \exp(-\frac{n\lambda}{2} f_\nu^2 \rho_\nu) \\
&= \prod_{\nu > m}^{\infty} (1 + n\lambda\rho_\nu/\theta_\nu)^{1/2} \cdot \exp\left(-\frac{n\lambda}{2} \sum_{\nu > m}^{\infty} f_\nu^2 \rho_\nu\right).
\end{aligned}$$

Note that in Proposition A.2, $\prod_{\nu > m}^{\infty} (1 + n\lambda\rho_\nu/\theta_\nu)^{1/2}$ is convergent since $\sum_{\nu > m} \rho_\nu/\theta_\nu < \infty$. Therefore, by Proposition A.2, we have

$$\begin{aligned}
&\frac{d\Pi'_\lambda}{d\Pi'}(\{f_\nu : \nu \geq 1\}) \\
&= \frac{\prod_{\nu=1}^{m} \left(\frac{2\pi\sigma_\nu^2}{1+n\lambda\sigma_\nu^2}\right)^{-1/2} \exp\left(-\frac{(1+n\lambda\sigma_\nu^2)f_\nu^2}{2\sigma_\nu^2}\right)}{\prod_{\nu=1}^{m} (2\pi\sigma_\nu^2)^{-1/2} \exp\left(-\frac{f_\nu^2}{2\sigma_\nu^2}\right)} \times \frac{d\tilde{\Pi}_\lambda}{d\tilde{\Pi}}(\{f_\nu : \nu > m\}) \\
&= \prod_{\nu=1}^{m} (1 + n\lambda\sigma_\nu^2)^{1/2} \exp\left(-\frac{n\lambda}{2} f_\nu^2\right) \times \prod_{\nu > m}^{\infty} (1 + n\lambda\rho_\nu/\theta_\nu)^{1/2} \\
&\quad \times \exp\left(-\frac{n\lambda}{2} \sum_{\nu > m}^{\infty} f_\nu^2 \rho_\nu\right) \\
&= \prod_{\nu=1}^{m} (1 + n\lambda\sigma_\nu^2)^{1/2} \prod_{\nu > m}^{\infty} (1 + n\lambda\rho_\nu/\theta_\nu)^{1/2} \times \exp\left(-\frac{n\lambda}{2} \sum_{\nu=1}^{\infty} f_\nu^2 \gamma_\nu\right).
\end{aligned}$$

Then for any measurable $S \subseteq S^m(\mathbb{I})$, by change of variable, we have

$$
\begin{aligned}
\Pi_\lambda(S) &= \Pi'_\lambda(TS) \\
&= \int_{TS} d\Pi'_\lambda(\{f_\nu : \nu \geq 1\}) \\
&= \prod_{\nu=1}^m (1 + n\lambda\sigma_\nu^2)^{1/2} \prod_{\nu>m}^\infty (1 + n\lambda\rho_\nu/\theta_\nu)^{1/2} \\
&\quad \cdot \int_{TS} \exp\left(-\frac{n\lambda}{2}\sum_{\nu=1}^\infty f_\nu^2 \gamma_\nu\right) d\Pi'(\{f_\nu : \nu \geq 1\}) \\
&= \prod_{\nu=1}^m (1 + n\lambda\sigma_\nu^2)^{1/2} \prod_{\nu>m}^\infty (1 + n\lambda\rho_\nu/\theta_\nu)^{1/2} \\
&\quad \cdot \int_{TS} \exp\left(-\frac{n\lambda}{2}J(T^{-1}(\{f_\nu : \nu \geq 1\}))\right) d(\Pi \circ T^{-1})(\{f_\nu : \nu \geq 1\}) \\
&= \prod_{\nu=1}^m (1 + n\lambda\sigma_\nu^2)^{1/2} \prod_{\nu>m}^\infty (1 + n\lambda\rho_\nu/\theta_\nu)^{1/2} \int_S \exp\left(-\frac{n\lambda}{2}J(f)\right) d\Pi(f).
\end{aligned}
$$

This completes the proof of the lemma. $\qquad\square$

The proof of Lemma 3.4 requires a concentration result (Lemma A.3). Let $\{\tilde{\varphi}_\nu : \nu \geq 1\}$ be a bounded orthonormal basis of $L^2(\mathbb{I})$ under usual $L^2$ inner product. For any $b \in [0, \beta]$, define

$$
\tilde{H}_b = \{\sum_{\nu=1}^\infty f_\nu \tilde{\varphi}_\nu : \sum_{\nu=1}^\infty f_\nu^2 \rho_\nu (\theta_\nu/\rho_\nu)^{b/\beta} < \infty\}.
$$

Then $\tilde{H}_b$ can be viewed as a version of Sobolev space with regularity $m + b/2$. Define $\tilde{G} = \sum_{\nu=1}^\infty v_\nu \tilde{\varphi}_\nu$, a centered GP, and $\tilde{f}_0 = \sum_{\nu=1}^\infty f_\nu^0 \tilde{\varphi}_\nu$. Define $\tilde{V}(f,g) = \langle f, g\rangle_{L^2} = \int_0^1 f(x)g(x)dx$, the usual $L^2$ inner product, $\tilde{J}(f) = \sum_{\nu=1}^\infty |\tilde{V}(f, \tilde{\varphi}_\nu)|^2 \rho_\nu$, a functional on $\tilde{H}_0$. For simplicity, denote $\tilde{V}(f) = \tilde{V}(f, f)$. Clearly, $\tilde{f}_0 \in \tilde{H}_\beta$. Since $\tilde{G}$ is a Gaussian process with covariance function

$$
\tilde{R}(s,t) = E\{\tilde{G}(s)\tilde{G}(t)\} = \sum_{\nu=1}^m \sigma_\nu^2 \tilde{\varphi}_\nu(s)\tilde{\varphi}_\nu(t) + \sum_{\nu>m} \theta_\nu^{-1}\tilde{\varphi}_\nu(s)\tilde{\varphi}_\nu(t),
$$

it follows by [56] that $\tilde{H}_\beta$ is the RKHS of $\tilde{G}$. For any $\tilde{H}_b$ with $0 \leq b \leq \beta$, define inner product

$$
\langle \sum_{\nu=1}^\infty f_\nu \tilde{\varphi}_\nu, \sum_{\nu=1}^\infty g_\nu \tilde{\varphi}_\nu\rangle_b = \sum_{\nu=1}^m \sigma_\nu^{-2} f_\nu g_\nu + \sum_{\nu>m} f_\nu g_\nu \rho_\nu (\theta_\nu/\rho_\nu)^{b/\beta}.
$$

Let $\|\cdot\|_b$ be the norm corresponding to the above inner product.

**Lemma A.3.** *Let $d_n$ be any positive sequence. If Condition (**S**) holds, then there exists $\omega \in \tilde{H}_\beta$ such that*

*(i).* $\tilde{V}(\omega - \tilde{f}_0) \leq \frac{1}{4}d_n^2$,

*(ii).* $\tilde{J}(\omega - \tilde{f}_0) \leq \frac{1}{4}d_n^{\frac{2(\beta-1)}{2m+\beta-1}}$,

*(iii).* $\|\omega\|_\beta^2 = O(d_n^{-\frac{2}{2m+\beta-1}})$.

*Proof of Lemma A.3.* Let $\omega = \sum_{\nu=1}^\infty \omega_\nu \tilde{\varphi}_\nu$, where $\omega_\nu = \frac{df_\nu^0}{d+(\sigma b_\nu)^\alpha}$, $\sigma = d_n^{2/(2m+\beta-1)}$, $b_\nu = \rho_\nu^{1/(2m)}$, $\alpha = m + (\beta-1)/2$, and $d > 0$ is a constant to be described. It is easy to see that for any $\nu$, $f_\nu^0 - \omega_\nu = \frac{(\sigma b_\nu)^\alpha f_\nu^0}{d+(\sigma b_\nu)^\alpha}$. Then

$$\tilde{V}(\omega - \tilde{f}_0) = \sum_{\nu=1}^\infty (f_\nu^0 - \omega_\nu)^2 = \sum_{\nu=1}^\infty \frac{|f_\nu^0|^2 (\sigma b_\nu)^{2\alpha}}{(d+(\sigma b_\nu)^\alpha)^2} \le \sigma^{2m+\beta-1} d^{-2} \sum_{\nu=1}^\infty |f_\nu^0|^2 \rho_\nu^{1+\frac{\beta-1}{2m}},$$

and

$$\begin{aligned}
\tilde{J}(\omega - \tilde{f}_0) &= \sum_{\nu=1}^\infty (f_\nu^0 - \omega_\nu)^2 \rho_\nu \\
&= \sigma^{\beta-1} \sum_{\nu=1}^\infty |f_\nu^0|^2 \rho_\nu^{1+\frac{\beta-1}{2m}} \left(d(\sigma b_\nu)^{-m} + (\sigma b_\nu)^{(\beta-1)/2}\right)^{-2} \\
&\le d^{-\frac{\beta-1}{\alpha}} \left(\left(\frac{2m}{\beta-1}\right)^{-\frac{m}{\alpha}} + \left(\frac{2m}{\beta-1}\right)^{\frac{\beta-1}{2\alpha}}\right)^{-2} \sigma^{\beta-1} \sum_{\nu=1}^\infty |f_\nu^0|^2 \rho_\nu^{1+\frac{\beta-1}{2m}}.
\end{aligned}$$

Therefore, we choose $d$ as a suitably large fixed constant such that (i) and (ii) hold.

To show (iii), observe that

$$\|\omega\|_\beta^2 = \sum_{\nu=1}^m \sigma_\nu^{-2} \omega_\nu^2 + \sum_{\nu > m} \omega_\nu^2 \theta_\nu \asymp \sum_{\nu=1}^m \sigma_\nu^{-2} |f_\nu^0|^2 + \sum_{\nu > m} \frac{d^2 |f_\nu^0|^2 \rho_\nu^{1+\frac{\beta-1}{2m}}}{(d+(\sigma b_\nu)^\alpha)^2} b_\nu = O(\sigma^{-1}).$$

The result follows by $\sigma = d_n^{2/(2m+\beta-1)}$. $\qquad\square$

*Proof of Lemma 3.4.* Let $\varepsilon \ge \lambda^{\frac{2m+\beta-1}{4m}}$. Hence, $\lambda = h^{2m} \le \varepsilon^{\frac{4m}{2m+\beta-1}}$. It follows by Lemma A.3 by replacing $d_n$ therein by $\varepsilon$, by Gaussian correlation inequality (see Theorem 1.1 of [27]), by Cameron-Martin theorem (see [6] or [25, eqn (4.18)]) and [20, Example 4.5] that

$$\begin{aligned}
&\Pi(\|f - f_0\| \le \varepsilon) \\
&= P(\|G - f_0\| \le \varepsilon) \\
&\ge P(V(G - f_0) \le \varepsilon^2/2, \lambda J(G - f_0) \le \varepsilon^2/2) \\
&\ge P(V(G - f_0) \le \varepsilon^2/2, J(G - f_0) \le \varepsilon^{\frac{2(\beta-1)}{2m+\beta-1}}/2) \\
&= P(\tilde{V}(\tilde{G} - \tilde{f}_0) \le \varepsilon^2/2, \tilde{J}(\tilde{G} - \tilde{f}_0) \le \varepsilon^{\frac{2(\beta-1)}{2m+\beta-1}}/2) \\
&\ge P(\tilde{V}(\tilde{G} - \omega) \le (1/\sqrt{2} - 1/2)^2 \varepsilon^2, \tilde{J}(\tilde{G} - \omega) \le (1/\sqrt{2} - 1/2)^2 \varepsilon^{\frac{2(\beta-1)}{2m+\beta-1}}) \\
&\ge \exp(-\|\omega\|_\beta^2/2) \times P(\tilde{V}(\tilde{G}) \le (1/\sqrt{2} - 1/2)^2 \varepsilon^2, \tilde{J}(\tilde{G}) \le (1/\sqrt{2} - 1/2)^2 \varepsilon^{\frac{2(\beta-1)}{2m+\beta-1}}) \\
&\ge \exp(-\|\omega\|_\beta^2/2) \times P(\tilde{V}(\tilde{G}) \le (1/\sqrt{2} - 1/2)^2 \varepsilon^2/2) \times P(\tilde{J}(\tilde{G}) \le (1/\sqrt{2} - 1/2)^2 \varepsilon^{\frac{2(\beta-1)}{2m+\beta-1}}/2) \\
&\ge \exp(-c_0 \varepsilon^{-\frac{2}{2m+\beta-1}}),
\end{aligned}$$

where $c_0 > 0$ is a universal constant. $\qquad\square$

*Proof of Theorem 3.5.* Let $T : f \mapsto \{f_\nu : \nu \geq 1\}$ be the one-to-one map from $S^m(\mathbb{I})$ to $\mathcal{R}_m$, as defined in the proof of Lemma 3.3.

Let $\Pi'_W$ and $\Pi'$ be the probability measures induced by $\{a_{n,\nu}\widehat{f}_\nu + b_{n,\nu}\tau_\nu v_\nu : \nu \geq 1\}$ and $\{v_\nu : \nu \geq 1\}$. Then $d\Pi'_W/d\Pi'$ equals $\lim_{N\to\infty} p_{1,N}(f_1, \ldots, f_N)/p_{2,N}(f_1, \ldots, f_N)$ (see [46, Section III]), where $p_{1,N}$ and $p_{2,N}$ are the probability densities under $f_\nu \sim a_{n,\nu}\widehat{f}_\nu + b_{n,\nu}\tau_\nu v_\nu$ and $f_\nu \sim v_\nu$, $\nu = 1, \ldots, N$, respectively. A direct evaluation leads to that

$$\frac{d\Pi'_W}{d\Pi'}(\{f_\nu : \nu \geq 1\}) = C_{n,\lambda} \exp(-\frac{n}{2}\sum_{\nu=1}^{\infty}(f_\nu - \widehat{f}_\nu)^2(1 + \lambda\gamma_\nu)),$$

where

$$C_{n,\lambda} = \prod_{\nu=1}^{\infty}(b_{n,\nu}\tau_\nu)^{-1} \exp\left(\frac{n}{2}\sum_{\nu=1}^{\infty}\frac{\tau_\nu^2(1 + \lambda\gamma_\nu)}{\tau_\nu^2 + n(1 + \lambda\gamma_\nu)}\widehat{f}_\nu^2\right).$$

Since $\sum_\nu \widehat{f}_\nu^2\gamma_\nu < \infty$ and $\beta > 1$, it is not hard to see that $C_{n,\lambda}$ is an almost surely finite constant.

For any $B \subseteq \mathcal{R}_m$, $\Pi_W(T^{-1}B) = \Pi'_W(B)$, and $\Pi(T^{-1}B) = \Pi'(B)$. To see this, note that

$$\Pi_W(T^{-1}B) = P(W \in T^{-1}B) = P(\{a_{n,\nu}\widehat{f}_\nu + b_{n,\nu}\tau_\nu v_\nu : \nu \geq 1\} \in B) = \Pi'_W(B), \text{ and}$$

$$\Pi(T^{-1}B) = P(G \in T^{-1}B) = P(\{v_\nu : \nu \geq 1\} \in B) = \Pi'(B).$$

By change of variable, for any $\Pi$-measurable $S \subseteq S^m(\mathbb{I})$,

$$
\begin{aligned}
\Pi_W(S) &= \Pi'_W(TS) \\
&= \int_{TS} d\Pi'_W(\{f_\nu : \nu \geq 1\}) \\
&= C_{n,\lambda}\int_{TS} \exp(-\frac{n}{2}\sum_{\nu=1}^{n}(f_\nu - \widehat{f}_\nu)^2(1 + \lambda\gamma_\nu))d\Pi'(f_\nu : \nu \geq 1) \\
&= C_{n,\lambda}\int_{S} \exp(-\frac{n}{2}\|f - \widehat{f}_{n,\lambda}\|^2)d\Pi(f).
\end{aligned}
$$

In particular, let $S = S^m(\mathbb{I})$ in the above equations, we get that

$$C_{n,\lambda} = \left(\int_{S^m(\mathbb{I})} \exp(-\frac{n}{2}\|f - \widehat{f}_{n,\lambda}\|^2)d\Pi(f)\right)^{-1}.$$

This proves the desired result. $\qquad\square$

### *A.2. Proofs in Section 4*

Before proving Theorem 4.1, we give some preliminary results.

**Lemma A.4.** *As $n \to \infty$, we have*

$$\frac{n\|W_n\|_2^2 - \zeta_{1,n}}{\sqrt{2\zeta_{2,n}}} \xrightarrow{d} N(0, 1),$$

*where $\zeta_{k,n} = \sum_{\nu=1}^{\infty}\frac{1}{(1 + \lambda\gamma_\nu + n^{-1}\tau_\nu^2)^k}$.*

*Proof of Lemma A.4.* Let $\eta_\nu = \tau_\nu v_\nu$. Then $\eta_\nu$ is a sequence of *iid* standard normals. Note that

$$\|W_n\|_2^2 = \sum_{\nu=1}^\infty \frac{\eta_\nu^2}{\tau_\nu^2 + n(1 + \lambda\gamma_\nu)}.$$

Let $U_n = (n\|W_n\|_2^2 - \zeta_{1,n})/\sqrt{2\zeta_{2,n}}$, then we have

$$U_n = \frac{1}{\sqrt{2\zeta_{2,n}}} \sum_{\nu=1}^\infty \frac{n(\eta_\nu^2 - 1)}{\tau_\nu^2 + n(1 + \lambda\gamma_\nu)}.$$

By straightforward calculations and Taylor's expansion of $\log(1 - x)$, it can be shown that the logarithm of the moment generating function of $U_n$ equals

$$\log E\{\exp(tU_n)\} = t^2/2 + O\left(t^3 \zeta_{2,n}^{-3/2} \zeta_{3,n}\right). \tag{A.9}$$

It can be examined that $\zeta_{2,n} \asymp n^{1/(2m+\beta)}$ and $\zeta_{3,n} \asymp n^{1/(2m+\beta)}$, so the remainder term in (A.9) is $O(n^{-1/(2(2m+\beta))}) = o(1)$. So $\lim_{n\to\infty} E\{\exp(tU_n)\} = \exp(t^2/2)$. Proof is completed. □

Define

$$R(x,y) = \sum_{\nu=1}^\infty \frac{\varphi_\nu(x)\varphi_\nu(y)}{1 + \lambda\gamma_\nu + n^{-1}\tau_\nu^2}, \quad x, y \in \mathbb{I}. \tag{A.10}$$

**Lemma A.5.** $\sup_{x,y\in\mathbb{I}} |R(x,y)| \lesssim h^{-1}$ *and* $\sup_{x,y\in\mathbb{I}} \left|\frac{\partial}{\partial x} R(x,y)\right| \lesssim h^{-2}$.

*Proof of Lemma A.5.* For any $g \in S^m(\mathbb{I})$ and $x \in \mathbb{I}$, it follows from [13, Lemmas (2.10) and (2.17)] that there exist constants $c', c'', c'''$ s.t.

$$\begin{aligned}
\left|\langle g, \frac{\partial}{\partial x} K_x\rangle\right| &= \left|\frac{\partial}{\partial x}\langle g, K_x\rangle\right| \\
&= |g'(x)| \leq c' h^{-1/2} \sqrt{\|g'\|_{L^2}^2 + h^2\|g''\|_{L^2}^2} \\
&= c' h^{-3/2} \sqrt{h^2\|g'\|_{L^2}^2 + h^4\|g''\|_{L^2}^2} \\
&\leq c' h^{-3/2} \sqrt{(\|g\|_{L^2}^2 + h^2\|g'\|_{L^2}^2) + (\|g\|_{L^2}^2 + h^4\|g''\|_{L^2}^2)} \\
&\leq c' c'' h^{-3/2} \sqrt{\|g\|_{L^2}^2 + h^{2m}\|g^{(m)}\|_{L^2}^2} \\
&\leq c''' h^{-3/2} \|g\|.
\end{aligned}$$

This implies that $\|\frac{\partial}{\partial x} K_x\| \leq c''' h^{-3/2}$. For convenience, let $R_y(\cdot) = R(\cdot, y)$. It is easy to see that

$$\begin{aligned}
\|R_y\|^2 &= \sum_{\nu=1}^\infty \frac{\varphi_\nu(y)^2}{(1 + \lambda\gamma_\nu + n^{-1}\tau_\nu^2)^2}(1 + \lambda\gamma_\nu) \\
&\leq \sum_{\nu=1}^\infty \frac{\varphi_\nu(y)^2}{1 + \lambda\gamma_\nu} = K(y,y) \leq c_K^2 h^{-1}.
\end{aligned}$$

This implies that $|R(x,y)| = |\langle R_y, K_x\rangle| \leq \|R_y\| \cdot \|K_x\| \leq c_K^2 h^{-1}$. This also leads to that, for any $x, y \in \mathbb{I}$,

$$\left|\frac{\partial}{\partial x} R(x,y)\right| = \left|\langle R_y, \frac{\partial}{\partial x} K_x\rangle\right| \leq \|R_y\| \cdot \|\frac{\partial}{\partial x} K_x\| \leq c''' c_K h^{-2}.$$

The desired result follows by the fact that both $c_k$ and $c'''$ are universal constants free of $x, y$. □

*Proof of Theorem 4.1.* By direct examinations, we can show that the Rate Conditions $(\mathbf{R})$, $n\tilde{r}_n^2(\tilde{r}_n b_{n1} + b_{n2}) = o(1)$, $nhD_n^2 = o(1)$ are all satisfied.

It is sufficient to investigate the $P_{f_0}^n$-probability of the event $\|\tilde{f}_{n,\lambda} - f_0\|_2 \le r_n(\alpha)$. To achieve this goal, we first prove the following fact:

$$|z_n(\alpha) - z_\alpha| = o_{P_{f_0}^n}(1), \tag{A.11}$$

where $z_\alpha = \Phi^{-1}(1 - \alpha)$ and $\Phi$ is the c.d.f. of $N(0,1)$, and $z_n(\alpha) = (nr_n(\alpha)^2 - \zeta_{1,n})/\sqrt{2\zeta_{2,n}}$. The proof of the theorem follows by (A.11) and a careful analysis of $f_0 - \bar{f}_{N,\lambda}$.

We first show (A.11). It follows by Theorem 3.2 that,

$$|P(R_n(\alpha)|\mathbf{D}_n) - P_0(R_n(\alpha))| \le \sup_{B \in \mathcal{B}} |P(B|\mathbf{D}_n) - P_0(B)| = o_{P_{f_0}^n}(1).$$

Together with $P(R_n(\alpha)|\mathbf{D}_n) = 1 - \alpha$, we have $|P_0(R_n(\alpha)) - (1 - \alpha)| = o_{P_{f_0}^n}(1)$. Since $W = \tilde{f}_{n,\lambda} + W_n$,

$$\begin{aligned} P_0(R_n(\alpha)) &= P(W \in R_n(\alpha)|\mathbf{D}_n) \\ &= P(\|W_n\|_2 \le r_n(\alpha)|\mathbf{D}_n) = P(U_n \le z_n(\alpha)|\mathbf{D}_n), \end{aligned}$$

and $P(U_n \le z_\alpha) \to 1 - \alpha$, where $U_n$ is defined in the proof of Lemma A.4, we get that

$$|P(U_n \le z_n(\alpha)|\mathbf{D}_n) - P(U_n \le z_\alpha)| = o_{P_{f_0}^n}(1), \tag{A.12}$$

where $U_n = (n\|W_n\|_2^2 - \zeta_{1,n})/\sqrt{2\zeta_{2,n}}$. Let $\Phi_n$ be the c.d.f. of $U_n$. Since $U_n$ is independent of the data, we have from (A.12) that

$$|\Phi_n(z_n(\alpha)) - \Phi_n(z_\alpha)| = o_{P_{f_0}^n}(1). \tag{A.13}$$

Now for any $\varepsilon > 0$, if $|z_n(\alpha) - z_\alpha| \ge \varepsilon$, then either $|\Phi_n(z_n(\alpha)) - \Phi_n(z_\alpha)| \ge \Phi_n(z_\alpha + \varepsilon) - \Phi_n(z_\alpha)$ or $|\Phi_n(z_n(\alpha)) - \Phi_n(z_\alpha)| \ge \Phi_n(z_\alpha) - \Phi_n(z_\alpha - \varepsilon)$. Since $\Phi_n$ pointwise congerges to $\Phi$, both $\Phi_n(z_\alpha + \varepsilon) - \Phi_n(z_\alpha)$ and $\Phi_n(z_\alpha) - \Phi_n(z_\alpha - \varepsilon)$ are asymptotically lower bounded by some positive numbers (possibly depending on $\varepsilon$). This implies by (A.13) that (A.11) holds.

Next we prove the theorem. Define $Rem_n = \hat{f}_{n,\lambda} - f_0 - S_{n,\lambda}(f_0)$. It follows by *Functional Bahadur Representation* ([43] or [44, Theorem 3.4]) that $\|Rem_n\| = O_{P_{f_0}^n}(D_n)$ with $D_n = a_n + b_n$. By direct examination, we have

$$\begin{aligned} \tilde{f}_{n,\lambda} - f_0 &= \sum_{\nu=1}^\infty \left(a_{n,\nu} V(\hat{f}_{n,\lambda}, \varphi_\nu) - f_\nu^0\right)\varphi_\nu \\ &= \sum_{\nu=1}^\infty \left(a_{n,\nu} V(Rem_n + f_0 + S_{n,\lambda}(f_0), \varphi_\nu) - f_\nu^0\right)\varphi_\nu \\ &= \sum_{\nu=1}^\infty a_{n,\nu} V(Rem_n, \varphi_\nu)\varphi_\nu + \sum_{\nu=1}^\infty (a_{n,\nu} - 1)f_\nu^0\varphi_\nu \\ &\quad + \sum_{\nu=1}^\infty a_{n,\nu} V(\frac{1}{n}\sum_{i=1}^n \epsilon_i K_{X_i}, \varphi_\nu)\varphi_\nu - \sum_{\nu=1}^\infty a_{n,\nu} V(\mathcal{P}_\lambda f_0, \varphi_\nu)\varphi_\nu, \end{aligned} \tag{A.14}$$

where $\epsilon_i = Y_i - \dot{A}(f_0(X_i))$. Denote the four terms in the above equation by $T_1, T_2, T_3, T_4$.

Since $a_{n,\nu} \leq 1$, it is easy to see that

$$
\begin{aligned}
\|T_1\|_2^2 &= \sum_{\nu=1}^{\infty} a_{n,\nu}^2 |V(Rem_n, \varphi_\nu)|^2 \\
&\leq \sum_{\nu=1}^{\infty} |V(Rem_n, \varphi_\nu)|^2 = \|Rem_n\|_2^2 \leq \|Rem_n\|^2 = O_{P_{f_0}^n}(D_n^2).
\end{aligned} \tag{A.15}
$$

Using $h \asymp n^{-1/(2m+\beta)}$ and a direct algebra we get that

$$
\begin{aligned}
\|T_2\|_2^2 &= \sum_{\nu=1}^{\infty} (a_{n,\nu} - 1)^2 |f_\nu^0|^2 \\
&\asymp \sum_{\nu=1}^{\infty} \left( \frac{\nu^{2m+\beta}}{\nu^{2m+\beta} + n(1 + \lambda \nu^{2m})} \right)^2 |f_\nu^0|^2 \\
&= o(n^{-\frac{2m+\beta-1}{2m+\beta}}) = o(n^{-1} h^{-1}).
\end{aligned}
$$

Meanwhile, it follows by Proposition 1 that

$$
\begin{aligned}
\|T_4\|_2^2 &= \sum_{\nu=1}^{\infty} a_{n,\nu}^2 |f_\nu^0|^2 \left( \frac{\lambda \gamma_\nu}{1 + \lambda \gamma_\nu} \right)^2 \\
&\leq \sum_{\nu=1}^{\infty} |f_\nu^0|^2 \left( \frac{\lambda \gamma_\nu}{1 + \lambda \gamma_\nu} \right)^2 \\
&\lesssim \sum_{\nu=1}^{\infty} |f_\nu^0|^2 (h\nu)^{2m+\beta-1} \frac{(h\nu)^{2m-\beta+1}}{(1 + (h\nu)^{2m})^2} \\
&= o(n^{-\frac{2m+\beta-1}{2m+\beta}}) = o(n^{-1} h^{-1}).
\end{aligned}
$$

It is easy to see that $R(x, x') = \sum_{\nu=1}^{\infty} a_{n,\nu} \frac{\varphi_\nu(x) \varphi_\nu(x')}{1 + \lambda \gamma_\nu}$ for any $x, x' \in \mathbb{I}$. Also define $R_x(\cdot) = R(x, \cdot)$. It is easy to see that $R_x \in S^m(\mathbb{I})$ for any $x \in \mathbb{I}$. Then it can be shown that $T_3 = \frac{1}{n} \sum_{i=1}^{n} \epsilon_i R_{X_i}$, leading to $\|T_3\|_2^2 = V(T_3, T_3) = \frac{1}{n^2} \sum_{i=1}^{n} \epsilon_i^2 V(R_{X_i}, R_{X_i}) + \frac{2}{n^2} \sum_{1 \leq i < k \leq n} \epsilon_i \epsilon_k V(R_{X_i}, R_{X_k})$. Define $W(n) = 2 \sum_{i < k} \epsilon_i \epsilon_k V(R_{X_i}, R_{X_k})$. Let $W_{ik} = 2 \epsilon_i \epsilon_k V(R_{X_i}, R_{X_k})$ for $1 \leq i < k \leq n$, then $W(n) = \sum_{1 \leq i < k \leq n} W_{ik}$. Note that $W(n)$ is clean in the sense of [12]. Let $\sigma^2(n) = E_{f_0}\{W(n)^2\}$ and $G_I, G_{II}, G_{IV}$ be defined as

$$
\begin{aligned}
G_I &= \sum_{1 \leq i < j \leq n} E_{f_0}\{W_{ij}^4\}, \\
G_{II} &= \sum_{1 \leq i < j < k \leq n} (E_{f_0}\{W_{ij}^2 W_{ik}^2\} + E_{f_0}\{W_{ji}^2 W_{jk}^2\} + E_{f_0}\{W_{ki}^2 W_{kj}^2\}) \\
G_{IV} &= \sum_{1 \leq i < j < k < l \leq n} (E_{f_0}\{W_{ij} W_{ik} W_{lj} W_{lk}\} + E_{f_0}\{W_{ij} W_{il} W_{kj} W_{kl}\} + E_{f_0}\{W_{ik} W_{il} W_{jk} W_{jl}\}).
\end{aligned}
$$

Since $\varphi_\nu$ are uniformly bounded, we get that

$$
\|R_x\|_2^2 = \sum_{\nu=1}^{\infty} \frac{|\varphi_\nu(x)|^2}{(1 + n^{-1} \tau_\nu^2 + \lambda \gamma_\nu)^2} \lesssim h^{-1},
$$

where "$\lesssim$" is free of $x$. This implies that $G_I = O(n^2 h^{-4})$ and $G_{II} = O(n^3 h^{-4})$. It can also be shown that for pairwise distinct $i, k, t, l$,

$$
\begin{aligned}
&E_{f_0}\{W_{ik}W_{il}W_{tk}W_{tl}\} \\
=\ & 2^4 E_{f_0}\{\epsilon_i^2\epsilon_k^2\epsilon_t^2\epsilon_l^2 V(R_{X_i}, R_{X_k})V(R_{X_i}, R_{X_l})V(R_{X_t}, R_{X_k})V(R_{X_t}, R_{X_l})\} \\
=\ & 2^4 \sum_{\nu=1}^{\infty} \frac{a_{n,\nu}^8}{(1+\lambda\gamma_\nu)^8} = O(h^{-1}),
\end{aligned}
$$

which implies that $G_{IV} = O(n^4 h^{-1})$. In the mean time, a straight algebra leads to that

$$
\begin{aligned}
\sigma^2(n) &=\ 4\binom{n}{2}\sum_{\nu=1}^{\infty}\frac{a_{n,\nu}^4}{(1+\lambda\gamma_\nu)^4} \\
&=\ 4\binom{n}{2}\sum_{\nu=1}^{\infty}\left(\frac{n}{\tau_\nu^2+n(1+\lambda\gamma_\nu)}\right)^4 = 2n(n-1)\zeta_{4,n} \asymp n^2 h^{-1}.
\end{aligned}
$$

Since $nh^2 \asymp n^{1-2/(2m+\beta)} \to \infty$, we get that $G_I, G_{II}$ and $G_{IV}$ are all of order $o(\sigma^4(n))$. Then it follows by [12] that as $n \to \infty$,

$$
\frac{W(n)}{n\sqrt{2\zeta_{4,n}}} \xrightarrow{d} N(0,1).
$$

Since $\zeta_{4,n} \asymp h^{-1}$, the above equation leads to that $W(n)/n = O_{P_{f_0}^n}(h^{-1/2})$. It follows by direct examination that $Var_{f_0}\{\sum_{i=1}^{n}\epsilon_i^2 V(R_{X_i}, R_{X_i})\} \le nE_{f_0}\{\epsilon_i^4\|R_{X_i}\|_2^4\} = O(nh^{-2})$, leading to that

$$
\begin{aligned}
\sum_{i=1}^{n}\epsilon_i^2 V(R_{X_i}, R_{X_i}) &=\ E_{f_0}\{\sum_{i=1}^{n}\epsilon_i^2 V(R_{X_i}, R_{X_i})\} + O_{P_{f_0}^n}(n^{1/2}h^{-1}) \\
&=\ n\zeta_{2,n} + O_{P_{f_0}^n}(n^{1/2}h^{-1}).
\end{aligned}
$$

Therefore, it follows by condition $nhD_n^2 = o(1)$ and the above analysis on $T_1, T_2, T_3, T_4$ that

$$
\begin{aligned}
nh\|\widetilde{f}_{n,\lambda} - f_0\|_2^2 &=\ nh\|T_3\|_2^2 + O_{P_{f_0}^n}(nhD_n^2) + o_{P_{f_0}^n}(1) \\
&=\ h\zeta_{2,n} + o_{P_{f_0}^n}(1).
\end{aligned}
\tag{A.16}
$$

In the end, note from (A.11) and $\zeta_{k,n} \asymp n^{1/(2m+\beta)}$ (see the proof of Lemma A.4) that $nr_n(\alpha)^2 = \zeta_{1,n} + \sqrt{2\zeta_{2,n}}z_\alpha + o_{P_{f_0}^n}(\sqrt{\zeta_{2,n}})$. Therefore, $nhr_n(\alpha)^2 = h\zeta_{1,n}(1 + o_{P_{f_0}^n}(1))$. Since $\liminf_{n\to\infty}(h\zeta_{1,n} - h\zeta_{2,n}) > 0$, we get that, with $P_{f_0}^n$-probability approaching one, $\|\widetilde{f}_{n,\lambda} - f_0\|_2 \le r_n(\alpha)$, i.e., $f_0 \in R_n(\alpha)$. Proof is completed. $\qquad\square$

Before proving Theorem 4.2, let us present two preliminary lemmas.

**Lemma A.6.** *As* $n \to \infty$,

$$
n\|W_n\|_\omega^2 \xrightarrow{d} \sum_{\nu=1}^{\infty}\omega_\nu\eta_\nu^2,
$$

*where* $\eta_\nu$ *are independent standard normal random variables.*

*Proof of Lemma A.6.* The proof follows by moment generating function approach and direct calculations, as in the proof of Lemma A.4. □

*Proof of Theorem 4.2.* By direct examinations, one can show that Rate Conditions (**R**), $n\tilde{r}_n^2(\tilde{r}_n b_{n1} + b_{n2}) = o(1)$, and $nD_n^2 = o(1)$ are all satisfied.

We first have the following fact:

$$|\sqrt{n} r_{\omega,n}(\alpha) - \sqrt{c_\alpha}| = o_{P_{f_0}^n}(1), \tag{A.17}$$

where $c_\alpha > 0$ satisfies $P(\sum_{\nu=1}^{\infty} \omega_\nu \eta_\nu^2 \le c_\alpha) = 1 - \alpha$ with $\eta_\nu$ being independent standard normal random variables. It follows from (A.17) that $n r_{\omega,n}(\alpha)^2 = c_\alpha + o_{P_{f_0}^n}(1)$. The proof of (A.17) is similar to the proof of (A.11) and is omitted.

Let $T_1, T_2, T_3, T_4$ be items defined in (A.14). It follows from the proof of Theorem 4.1 that $\|T_1\|_\omega^2 \le \|T_1\|_2^2 = O_{P_{f_0}^n}(D_n^2)$. So, $n\|T_1\|_\omega^2 = O_{P_{f_0}^n}(nD_n^2) = o_{P_{f_0}^n}(1)$ due to the condition $nD_n^2 = o(1)$.

It follows by condition $h \asymp n^{-1/(2m+\beta)}$, dominated convergence theorem and direct examinations that

$$
\begin{aligned}
\|T_2\|_\omega^2 &= \sum_{\nu=1}^{\infty} \omega_\nu (a_{n,\nu} - 1)^2 |f_\nu^0|^2 \\
&\asymp n^{-2} \sum_{\nu=1}^{\infty} \omega_\nu \frac{\nu^{2m+\beta+1}}{(1 + (h\nu)^{2m} + (h\nu)^{2m+\beta})^2} \times \nu^{2m+\beta-1} |f_\nu^0|^2 \\
&\lesssim n^{-1} \sum_{\nu=1}^{\infty} \frac{(h\nu)^{2m+\beta+1}}{(1 + (h\nu)^{2m} + (h\nu)^{2m+\beta})^2} \times \nu^{2m+\beta-1} |f_\nu^0|^2 = o(n^{-1}),
\end{aligned}
$$

and

$$
\begin{aligned}
\|T_4\|_\omega^2 &= \sum_{\nu=1}^{\infty} \omega_\nu a_{n,\nu}^2 \left( \frac{\lambda \gamma_\nu}{1 + \lambda \gamma_\nu} \right)^2 |f_\nu^0|^2 \\
&\lesssim \sum_{\nu=1}^{\infty} \omega_\nu \frac{(h\nu)^{2m-\beta+1}}{(1 + (h\nu)^{2m} + (h\nu)^{2m+\beta})^2} \times |f_\nu^0|^2 (h\nu)^{2m+\beta-1} \\
&\lesssim h^{2m+\beta} \sum_{\nu=1}^{\infty} \frac{(h\nu)^{2m-\beta}}{(1 + (h\nu)^{2m} + (h\nu)^{2m+\beta})^2} \times |f_\nu^0|^2 \nu^{2m+\beta-1} = o(n^{-1}).
\end{aligned}
$$

Next we handle $T_3$. By proof of Theorem 4.1, we have $T_3 = n^{-1} \sum_{i=1}^{n} \epsilon_i R_{X_i}$, where $\epsilon_i = Y_i - \dot{A}(f_0(X_i))$, which implies $n\|T_3\|_\omega^2 = n^{-1} \|\sum_{i=1}^{n} \epsilon_i R_{X_i}\|_\omega^2$.

Since $E_{f_0}\{\exp(|\epsilon|/C_0)\} \le C_1$, we can choose a constant $L > C_0$ such that $P_{f_0}^n(\mathcal{E}_n) \to 1$, where $\mathcal{E}_n = \{\max_{1 \le i \le n} |\epsilon_i| \le b_n \equiv L \log n\}$. We can even choose the above $L$ to be properly large so that the following rate condition holds:

$$h^{-1} n^{1/2} \exp(-b_n/(2C_0)) = o(1), \quad h^{-2} \exp(-b_n/(2C_0)) = o(1). \tag{A.18}$$

Define $H_n(\cdot) = n^{-1/2} \sum_{i=1}^{n} \epsilon_i R_{X_i}(\cdot)$ and $H_n^b(\cdot) = n^{-1/2} \sum_{i=1}^{n} \epsilon_i I(|\epsilon_i| \leq b_n) R_{X_i}(\cdot)$. Write $H_n = H_n - H_n^b - E_{f_0}\{H_n - H_n^b\} + H_n^b - E_{f_0}\{H_n^b\}$. Clearly, on $\mathcal{E}_n$, $H_n = H_n^b$, and hence,

$$
\begin{aligned}
|H_n(z) - H_n^b(z) - E_{f_0}\{H_n(z) - H_n^b(z)\}| &= |E_{f_0}\{H_n(z) - H_n^b(z)\}| \\
&= n^{1/2}|E_{f_0}\{\epsilon I(|\epsilon| > b_n) R_X(z)\}| \\
&\lesssim n^{1/2} h^{-1} E_{f_0}\{\epsilon^2\}^{1/2} P_{f_0}^n(|\epsilon| > b_n)^{1/2} \\
&\lesssim n^{1/2} h^{-1} \exp(-b_n/(2C_0)) = o(1),
\end{aligned}
$$

where the last $o(1)$-term follows by (A.18) and is free of the argument $z$. Thus,

$$
\sup_{z \in \mathbb{I}} |H_n(z) - H_n^b(z) - E_{f_0}\{H_n(z) - H_n^b(z)\}| = o_{P_{f_0}^n}(1). \tag{A.19}
$$

Define $\mathcal{R}_n = H_n^b - E_{f_0}\{H_n^b\}$ and $Z_n(e,x) = n^{1/2}(P_n(e,x) - P(e,x))$, where $P_n(e,x)$ is the empirical distribution of $(\epsilon, X)$ and $P(e,x)$ is the population distribution of $(\epsilon, X)$ under $P_{f_0}^n$-probability. It follows by Theorem 1 of [52] that

$$
\sup_{e \in \mathbb{R}, x \in \mathbb{I}} |Z_n(e,x) - W(t(e,x))| = O_{P_{f_0}^n}(n^{-1/2}(\log n)^2), \tag{A.20}
$$

where $W(\cdot, \cdot)$ is Brownian bridge indexed on $\mathbb{I}^2$, $t(e,x) = (F_1(x), F_2(e|x))$, $F_1$ is the marginal distribution of $X$ and $F_2$ is the conditional distribution of $\epsilon$ given $X$ both under $P_{f_0}^n$-probability. It can be seen that $\mathcal{R}_n(z) = \int_0^1 \int_{-b_n}^{b_n} e R_x(z) dZ_n(e,x)$. Define $\mathcal{R}_n^0(z) = \int_0^1 \int_{-b_n}^{b_n} e R_x(z) dW(t(e,x))$. Write $dU_n(x) = \int_{-b_n}^{b_n} e dZ_n(e,x)$, $dU_n^0(x) = \int_{-b_n}^{b_n} e dW(t(e,x))$. It follows from integration by parts where all quadratic variation terms are zero that

$$
U_n(x) = \int_{-b_n}^{b_n} e d_e Z_n(e,x) = Z_n(e,x) e|_{e=-b_n}^{b_n} - \int_{-b_n}^{b_n} Z_n(e,x) de,
$$

$$
U_n^0(x) = \int_{-b_n}^{b_n} e d_e W(t(e,x)) = W(t(e,x)) e|_{e=-b_n}^{b_n} - \int_{-b_n}^{b_n} W(t(e,x)) de,
$$

and hence, it follows by (A.20) that $\sup_{x \in \mathbb{I}} |U_n(x) - U_n^0(x)| = O_{P_{f_0}^n}(b_n n^{-1/2}(\log n)^2)$. It follows from integration by parts again and $\sup_{x,y \in \mathbb{I}} |\frac{\partial}{\partial x} R(x,y)| = O(h^{-2})$ (Lemma A.5) that

$$
\mathcal{R}_n(z) = \int_0^1 R_x(z) dU_n(x) = U_n(x) R(x,z)|_{x=0}^1 - \int_0^1 U_n(x) \frac{\partial}{\partial x} R(x,z) dx,
$$

$$
\mathcal{R}_n^0(z) = \int_0^1 R_x(z) dU_n^0(x) = U_n^0(x) R(x,z)|_{x=0}^1 - \int_0^1 U_n^0(x) \frac{\partial}{\partial x} R(x,z) dx,
$$

and hence,

$$
\sup_{z \in \mathbb{I}} |\mathcal{R}_n(z) - \mathcal{R}_n^0(z)| = O_{P_{f_0}^n}(h^{-2} b_n n^{-1/2}(\log n)^2) = o_{P_{f_0}^n}(1), \tag{A.21}
$$

due to the fact $2m + \beta > 4$, hence $h^{-2} n^{-1/2} b_n (\log n)^2 = O(n^{-1/2 + 2/(2m+\beta)}(\log n)^3) = o(1)$.

Next we handle the term $\mathcal{R}_n^0$. Write $W(s,t) = B(s,t) - stB(1,1)$, where $B(s,t)$ is standard Brownian motion indexed on $\mathbb{I}^2$. Define $\bar{\mathcal{R}}_n^0(z) = \int_0^1 \int_{-b_n}^{b_n} e R_x(z) dB(t(e,x))$. Let $F(e,x) = F_1(x)F_2(e,x)$ be the joint distribution of $(\epsilon, X)$. It is easy to see that

$$
\begin{aligned}
|\bar{\mathcal{R}}_n^0(z) - \mathcal{R}_n^0(z)| &= |B(1,1)| \cdot |\int_0^1 \int_{-b_n}^{b_n} e R_x(z) dF(e,x)| \\
&= |B(1,1)| \cdot |E_{f_0}\{\epsilon I(|\epsilon| \le b_n) R_X(z)\}| \\
&= |B(1,1)| \cdot |E_{f_0}\{\epsilon^2 I(|\epsilon| > b_n) R_X(z)\}| \\
&= O_{P_{f_0}^n}(h^{-1}\exp(-b_n/(2C_0))) = o_{P_{f_0}^n}(1),
\end{aligned}
$$

where the last equality follows by (A.18). Therefore, we have shown that

$$
\sup_{z \in \mathbb{I}} |\bar{\mathcal{R}}_n^0(z) - \mathcal{R}_n^0(z)| = o_{P_{f_0}^n}(1). \tag{A.22}
$$

By (A.19), (A.21) and (A.22) that

$$
n\|T_3\|_\omega^2 = \|H_n\|_\omega^2 = \|\bar{\mathcal{R}}_n^0\|_\omega^2 + o_{P_{f_0}^n}(1). \tag{A.23}
$$

Define $\widetilde{\mathcal{R}}(z) = \int_0^1 \int_{-\infty}^{\infty} e R_x(z) dB(t(e,x))$. Let $\Delta(z) = \widetilde{\mathcal{R}}(z) - \bar{\mathcal{R}}_n^0(z)$. Then

$$
\Delta(z) = \int_0^1 \int_{|e|>b_n} e R_x(z) dB(t(e,x)).
$$

For each $z$, $\Delta(z)$ is a zero-mean Gaussian random variable with variance

$$
\begin{aligned}
E_{f_0}\{\Delta(z)^2\} &= \int_0^1 \int_{|e|>b_n} e^2 R_x(z)^2 dF(e,x) \\
&\lesssim h^{-2} E_{f_0}\{\epsilon^2 I(|\epsilon| > b_n)\} = O(h^{-2}\exp(-b_n/(2C_0))) = o(1),
\end{aligned}
$$

where the last $o(1)$-term follows from (A.18) and is free of the argument $z$. Therefore,

$$
E_{f_0}\{\|\widetilde{\mathcal{R}} - \bar{\mathcal{R}}_n^0\|_\omega^2\} \lesssim E_{f_0}\{\|\Delta\|_{L^2}^2\} = \int_0^1 E_{f_0}\{\Delta(z)^2\} dz = o(1),
$$

implying that $\|\widetilde{\mathcal{R}} - \bar{\mathcal{R}}_n^0\|_\omega = o_{P_{f_0}^n}(1)$. Therefore, it follows by (A.23) that

$$
n\|T_3\|_\omega^2 = \|\widetilde{\mathcal{R}}\|_\omega^2 + o_{P_{f_0}^n}(1). \tag{A.24}
$$

It follows from the definition of $R(\cdot, \cdot)$ that

$$
\|\widetilde{\mathcal{R}}\|_\omega^2 = \sum_{\nu=1}^{\infty} \frac{\omega_\nu \tilde{\eta}_\nu^2}{(1 + \lambda\gamma_\nu + n^{-1}\tau_\nu^2)^2} \xrightarrow{d} \sum_{\nu=1}^{\infty} \omega_\nu \tilde{\eta}_\nu^2,
$$

where $\tilde{\eta}_\nu = \int_0^1 \int_{-\infty}^{\infty} e \varphi_\nu(x) dB(t(e,x))$. It is easy to see that for any $\nu, \mu$,

$$
E_{f_0}\{\tilde{\eta}_\nu \tilde{\eta}_\mu\} = E_{f_0}\{\epsilon^2 \varphi_\nu(X)\varphi_\mu(X)\} = E_{f_0}\{B(X)\varphi_\nu(X)\varphi_\mu(X)\} = V(\varphi_\nu, \varphi_\mu) = \delta_{\nu\mu},
$$

that is, $\tilde{\eta}_\nu$ are *iid* standard normal random variables. Combined with the above analysis of terms $T_1, T_2, T_3, T_4$, we have shown that as $n \to \infty$, $n\|f_0 - \tilde{f}_{n,\lambda}\|_\omega^2 \xrightarrow{d} \sum_{\nu=1}^{\infty} \omega_\nu \tilde{\eta}_\nu^2$. This implies that as $n \to \infty$, $P_{f_0}^n(f_0 \in R_n^\omega(\alpha)) = P_{f_0}^n(n\|f_0 - \tilde{f}_{n,\lambda}\|_\omega^2 \le c_\alpha) \to 1 - \alpha$. The proof is completed. $\qquad \square$

*Proof of Theorem 4.4.* Recall in the proof of Theorem 4.2 we show that Rate Conditions (**R**), $n\widetilde{r}_n^2(\widetilde{r}_n b_{n1} + b_{n2}) = o(1)$ and $nD_n^2 = o(1)$ are all satisfied.

It is easy to see that

$$F(W_n) \overset{d}{=} N(0, \theta_{1,n}^2). \tag{A.25}$$

Define $R_n^F(\alpha) = \{f \in S^m(\mathbb{I}) : |F(f) - F(\widetilde{f}_{n,\lambda})| \le r_{F,n}(\alpha)\}$. It follows by Theorem 3.2 that $|1 - \alpha - P_0(R_n^F(\alpha))| = o_{P_{f_0}^n}(1)$. It is easy to see that

$$P_0(R_n^F(\alpha)) = P(|F(W_n)| \le r_{F,n}(\alpha)|\mathbf{D}_n) = 2\Phi(r_{F,n}(\alpha)/\theta_{1,n}) - 1,$$

which leads to

$$|r_{F,n}(\alpha)/\theta_{1,n} - z_{\alpha/2}| = o_{P_{f_0}^n}(1). \tag{A.26}$$

Consider the decomposition (A.14) with $T_1, T_2, T_3, T_4$ being defined therein. It follows by (A.15) and rate condition $nD_n^2 = o(1)$ that $n\|T_1\|^2 = O_{P_{f_0}}(nD_n^2) = o_{P_{f_0}^n}(1)$. Meanwhile, it follows by Condition (**S'**), $n^{-1} \asymp h^{2m+\beta}$ and $\lambda = h^{2m}$ and direct examinations that

$$
\begin{aligned}
n\|T_2\|^2 &= n\sum_{\nu=1}^{\infty}(a_{n,\nu}-1)^2|f_\nu^0|^2(1+\lambda\gamma_\nu) \\
&\asymp n\sum_{\nu=1}^{\infty}\left(\frac{\nu^{2m+\beta}}{\nu^{2m+\beta}+n(1+\lambda\nu^{2m})}\right)^2|f_\nu^0|^2(1+\lambda\nu^{2m}) \\
&\asymp \sum_{\nu=1}^{\infty}\frac{(h\nu)^{2m+\beta}+(h\nu)^{4m+\beta}}{(1+(h\nu)^{2m}+(h\nu)^{2m+\beta})^2} \times |f_\nu^0|^2\nu^{2m+\beta} = o(1),
\end{aligned}
$$

and

$$
\begin{aligned}
n\|T_4\|^2 &= n\sum_{\nu=1}^{\infty}a_{n,\nu}^2\left(\frac{\lambda\gamma_\nu}{1+\lambda\gamma_\nu}\right)^2|f_\nu^0|^2(1+\lambda\gamma_\nu) \\
&\asymp \sum_{\nu=1}^{\infty}\frac{(h\nu)^{2m-\beta}}{1+(h\nu)^{2m}} \times |f_\nu^0|^2\nu^{2m+\beta} = o(1).
\end{aligned}
$$

Therefore, $\|\widetilde{f}_{n,\lambda} - f_0 - T_3\| = \|T_1 + T_2 + T_4\| = o_{P_{f_0}^n}(n^{-1/2})$. If follows from (4.3) that $|F(\widetilde{f}_{n,\lambda} - f_0) - F(T_3)| = o_{P_{f_0}^n}(h^{-r/2}n^{-1/2})$.

Recall $F(T_3) = \frac{1}{n}\sum_{i=1}^{n}\epsilon_i F(R_{X_i})$, where the kernel $R_x$ is defined in (A.10). We will derive the asymptotic distribution for $F(T_3)$. Let $s_n^2 = Var_{f_0}(\sum_{i=1}^{n}\epsilon_i F(R_{X_i}))$. It is easy to show that

$$s_n^2 = n^3\sum_{\nu=1}^{\infty}\frac{F(\varphi_\nu)^2}{(\tau_\nu^2+n(1+\lambda\gamma_\nu))^2} = n^3\theta_{2,n}^2 \asymp nh^{-r}.$$

By (4.3) and $\|R_x\| \le c_K h^{-1/2}$ (see proof of Lemma A.5), we get $|F(R_x)| \le \kappa h^{-r/2}\|R_x\| \le \kappa c_K h^{-(1+r)/2}$. Meanwhile,

$$E_{f_0}\{\epsilon^2 F(R_X)^2\} = n^2\sum_{\nu=1}^{\infty}\frac{F(\varphi_\nu)^2}{(\tau_\nu^2+n(1+\lambda\gamma_\nu))^2} = n^2\theta_{2,n}^2 \asymp h^{-r}. \tag{A.27}$$

By Assumption A1 there exists a constant $M_4$ s.t. $E_{f_0}\{\epsilon^4|X\} \le M_4$ a.s. Then for any $\delta > 0$,

$$\frac{1}{s_n^2} \sum_{i=1}^n E_{f_0}\{\epsilon_i^2 F(R_{X_i})^2 I(|\epsilon_i F(R_{X_i})| \ge \delta s_n)\}$$

$$\le \frac{n}{s_n^2}(\delta s_n)^{-2} E_{f_0}\{\epsilon^4 F(R_X)^4\}$$

$$\lesssim \frac{n}{s_n^2}(\delta s_n)^{-2} h^{-(1+r)} E_{f_0}\{\epsilon^2 F(R_X)^2\} \lesssim \delta^{-2} n^{-1} h^{-1} = o(1),$$

where the last $o(1)$-term follows by $h \asymp h_*$ and $2m + \beta > 1$. By Lindeberg's central limit theorem, as $n \to \infty$,

$$\frac{F(T_3)}{\sqrt{n}\theta_{2,n}} = \frac{1}{s_n} \sum_{i=1}^n \epsilon_i F(R_{X_i}) \xrightarrow{d} N(0,1). \tag{A.28}$$

By condition $n^2\theta_{2,n}^2 \asymp h^{-r}$, we have

$$\left|\frac{F(\widetilde{f}_{n,\lambda} - f_0 - T_3)}{\sqrt{n}\theta_{2,n}}\right| = o_{P_{f_0}^n}\left(\frac{h^{-r/2}n^{-1/2}}{\sqrt{n}\theta_{2,n}}\right) = o_{P_{f_0}^n}(1).$$

It follows by (A.26) that

$$\frac{r_{F,n}(\alpha)}{\sqrt{n}\theta_{2,n}} = \frac{\theta_{1,n}}{\sqrt{n}\theta_{2,n}} \times z_{\alpha/2}(1 + o_{P_{f_0}^n}(1)).$$

It can be easily seen that

$$\frac{\theta_{1,n}^2}{n\theta_{2,n}^2} = \frac{\sum_{\nu=1}^\infty \frac{F(\varphi_\nu)^2}{1+\lambda\gamma_\nu+n^{-1}\tau_\nu^2}}{\sum_{\nu=1}^\infty \frac{F(\varphi_\nu)^2}{(1+\lambda\gamma_\nu+n^{-1}\tau_\nu^2)^2}} \ge 1,$$

together with (A.28) we get that

$$P_{f_0}^n(|F(f_0) - F(\widetilde{f}_{n,\lambda})| \le r_{F,n}(\alpha))$$

$$= P_{f_0}^n\left(\left|\frac{F(\widetilde{f}_{n,\lambda} - f_0 - T_3)}{\sqrt{n}\theta_{2,n}} + \frac{F(T_3)}{\sqrt{n}\theta_{2,n}}\right| \le \frac{r_{F,n}(\alpha)}{\sqrt{n}\theta_{2,n}}\right)$$

$$\ge P_{f_0}^n\left(\left|\frac{F(\widetilde{f}_{n,\lambda} - f_0 - T_3)}{\sqrt{n}\theta_{2,n}} + \frac{F(T_3)}{\sqrt{n}\theta_{2,n}}\right| \le z_{\alpha/2}(1 + o_{P_{f_0}}(1))\right) \longrightarrow 1 - \alpha. \tag{A.29}$$

Notice that when $0 < \sum_{\nu=1}^\infty F(\varphi_\nu)^2 < \infty$, $\frac{\theta_{1,n}^2}{n\theta_{2,n}^2} \to 1$, leading to that the probability in (A.29) approaches exactly $1 - \alpha$. Proof is completed. □

*Proof of Proposition 4.5.* Under the setup of Proposition 4.5, it follows from [44] that $\ddot{A}(\cdot) \equiv 1$, and hence, (2.5) becomes the following uniform free beam problem:

$$\varphi_\nu^{(4)}(\cdot) = \rho_\nu \varphi_\nu(\cdot), \ \varphi_\nu^{(j)}(0) = \varphi_\nu^{(j)}(1) = 0, \ j = 2, 3. \tag{A.30}$$

The eigenvalues satisfy $\lim_{\nu\to\infty} \rho_\nu/(\pi\nu)^4 = 1$; see [21, Problem 3.10]. The normalized solutions to (A.30) are

$$\varphi_1(z) = 1, \ \varphi_2(z) = \sqrt{3}(2z - 1), \tag{A.31}$$

$$\varphi_{2k+1}(z) = \frac{\sin(\gamma_{2k+1}(z-1/2))}{\sin(\gamma_{2k+1}/2)} + \frac{\sinh(\gamma_{2k+1}(z-1/2))}{\sinh(\gamma_{2k+1}/2)}, \ k \geq 1. \tag{A.32}$$

$$\varphi_{2k+2}(z) = \frac{\cos(\gamma_{2k+2}(z-1/2))}{\cos(\gamma_{2k+2}/2)} + \frac{\cosh(\gamma_{2k+2}(z-1/2))}{\cosh(\gamma_{2k+2}/2)}, \ k \geq 1, \tag{A.33}$$

where $\gamma_\nu = \rho_\nu^{1/4}$ satisfying $\cos(\gamma_\nu)\cosh(\gamma_\nu) = 1$; see [5], page 295–296].

*Proof of (i).* By direct examinations, it can be shown that when $\nu \geq 3$ is odd, $\cos(x)\cosh(x) = 1$ has a unique solution in $((\nu+1/2)\pi, (\nu+1)\pi)$, that is, $\gamma_\nu \in ((\nu+1/2)\pi, (\nu+1)\pi)$; when $\nu \geq 3$ is even, $\cos(x)\cosh(x) = 1$ has a unique solution in $(\nu\pi, (\nu+1/2)\pi)$, that is, $\gamma_\nu \in (\nu\pi, (\nu+1/2)\pi)$. Consequently, for any $k \geq 1$, $0 < \gamma_{2k+2} - \gamma_{2k+1} < \pi$.

Let $\delta_0$ be constant such that $0 < \delta_0 < \pi/2 - \pi|z-1/2|$, and $d_0 = \min\{\sin^2(\delta_0), \cos^2(\delta_0 + \pi|z-1/2|)\}$. Clearly, $d_0 > 0$ is a constant. It is easy to see that when $k \to \infty$,

$$\frac{\sinh(\gamma_{2k+1}(z-1/2))}{\sinh(\gamma_{2k+1}/2)} \to 0, \ \text{and} \ \frac{\cosh(\gamma_{2k+2}(z-1/2))}{\cosh(\gamma_{2k+2}/2)} \to 0.$$

Then for arbitrarily small $\varepsilon \in (0, d_0/8)$, there exists $N$ s.t. for any $k \geq N$,

$$\varphi_{2k+1}(z)^2 \geq \frac{1}{2}\sin^2(\gamma_{2k+1}(z-1/2)) - \varepsilon \ \text{ and } \ \varphi_{2k+2}(z)^2 \geq \frac{1}{2}\cos^2(\gamma_{2k+2}(z-1/2)) - \varepsilon.$$

Let $\phi_k' = (\gamma_{2k+2} - \gamma_{2k+1})(z-1/2)$. Then $|\phi_k'| \leq \pi|z-1/2| < \pi/2$. There exists an integer $l_k$ s.t. $\gamma_{2k+1}(z-1/2) = \phi_k + l_k\pi$, where $\phi_k \in [0, \pi)$. Then, $\sin^2(\gamma_{2k+1}(z-1/2)) = \sin^2(\phi_k)$ and $\cos^2(\gamma_{2k+2}(z-1/2)) = \cos^2(\phi_k + \phi_k')$.

If $0 \leq \phi_k \leq \delta_0$, then it can be seen that

$$-\pi|z-1/2| \leq \phi_k' \leq \phi_k + \phi_k' \leq \delta_0 + \phi_k' \leq \delta_0 + \pi|z-1/2|.$$

Therefore, $\cos^2(\phi_k + \phi_k') \geq \cos^2(\delta_0 + \pi|z-1/2|)$. If $\delta_0 < \phi_k < \pi - \delta_0$, then $\sin^2(\phi_k) \geq \sin^2(\delta_0)$. If $\pi - \delta \leq \phi_k < \pi$, then it can be seen that

$$\pi - \delta_0 - \pi|z-1/2| \leq \phi_k + \phi_k' < \pi + \pi|z-1/2|.$$

Therefore, $\cos^2(\phi_k + \phi_k') \geq \cos^2(\delta_0 + \pi|z-1/2|)$. Consequently, for any $k \geq N$,

$$\begin{aligned}
\varphi_{2k+1}(z)^2 + \varphi_{2k+2}(z)^2 &\geq \frac{1}{2}(\sin^2(\gamma_{2k+1}(z-1/2)) + \cos^2(\gamma_{2k+2}(z-1/2))) - 2\varepsilon \\
&\geq \frac{1}{2}\min\{\sin^2(\delta_0), \cos^2(\delta_0 + \pi|z-1/2|)\} - 2\varepsilon \geq d_0/4.
\end{aligned}$$

Then we have

$$
\begin{aligned}
&\sum_{\nu > 2} \frac{h\varphi_\nu(z)^2}{(1 + \lambda\rho_\nu + (\lambda\rho_\nu)^{1+\beta/4})^j} \\
&= \sum_{k \geq 1} \frac{h\varphi_{2k+1}(z)^2}{(1 + \lambda\rho_{2k+1} + (\lambda\rho_{2k+1})^{1+\beta/4})^j} + \sum_{k \geq 1} \frac{h\varphi_{2k+2}(z)^2}{(1 + \lambda\rho_{2k+2} + (\lambda\rho_{2k+2})^{1+\beta/4})^j} \\
&\geq \sum_{k \geq 1} \frac{h\varphi_{2k+1}(z)^2 + h\varphi_{2k+2}(z)^2}{(1 + \lambda\rho_{2k+2} + (\lambda\rho_{2k+2})^{1+\beta/4})^j} \\
&\geq \sum_{k \geq N} \frac{hd_0/4}{(1 + \lambda\rho_{2k+2} + (\lambda\rho_{2k+2})^{1+\beta/4})^j} \\
&\gtrsim \sum_{k \geq N} \frac{h}{(1 + (k\pi h)^4 + (k\pi h)^{4+\beta})^j} \\
&\geq \int_N^\infty \frac{h}{(1 + (\pi h x)^4 + (\pi h x)^{4+\beta})^j} dx \\
&= \frac{1}{\pi} \int_{\pi N h}^\infty \frac{1}{(1 + x^4 + x^{4+\beta})^j} dx \xrightarrow{h \to 0} \frac{1}{\pi} \int_0^\infty \frac{1}{(1 + x^4 + x^{4+\beta})^j} dx > 0.
\end{aligned}
$$

This shows that condition (4.5) holds for $r = 1$.

*Proof of (ii).* Write $\omega = \sum_\nu \omega_\nu \varphi_\nu$ where $\omega_\nu$ is a square-summable real sequence. Then $F_\omega(\varphi_\nu) = \int_0^1 \omega(z)\varphi_\nu(z)dz = \omega_\nu$. Therefore, $\sum_\nu F_\omega(\varphi_\nu)^2 = \sum_\nu \omega_\nu^2 < \infty$. Meanwhile, since $\omega \neq 0$, $\sum_{\nu=1}^\infty F_\omega(\varphi_\nu)^2 > 0$. Consequently, for $j = 1, 2$, it follows by dominated convergence theorem that as $n \to \infty$,

$$
\sum_{\nu=1}^m \frac{F_\omega(\varphi_\nu)^2}{(1 + \lambda + n^{-1}\sigma_\nu^{-2})^j} + \sum_{\nu > m} \frac{F_\omega(\varphi_\nu)^2}{(1 + \lambda\rho_\nu + (\lambda\rho_\nu)^{1+\beta/(2m)})^j} \to \sum_{\nu=1}^\infty F_\omega(\varphi_\nu)^2 > 0.
$$

Hence (4.5) holds for $r = 0$. $\qquad\square$

### A.3. Proofs in Section 5

*Proof of Theorem 5.1.* It follows by Theorem 3.2 that

$$
\begin{aligned}
&\sup_{t \geq 0} |P(\sqrt{nh}\|f - \widehat{f}_{n,\lambda}\|_2 \leq t | \mathbf{D}_n) - P_0(\sqrt{nh}\|f - \widehat{f}_{n,\lambda}\|_2 \leq t)| \\
&\leq \sup_{B \in \mathcal{B}} |P(B|\mathbf{D}_n) - P_0(B)| = o_{P_{f_0}^n}(1).
\end{aligned}
$$

It is sufficient to prove

$$
\sup_{t \geq 0} |P_0(\sqrt{nh}\|f - \widehat{f}_{n,\lambda}\|_2 \leq t) - P_\star(\sqrt{nh}\|f - \widehat{f}_{n,\lambda}\|_2 \leq t)| = o_{P_{f_0}^n}(1). \tag{A.34}
$$

Define $\Delta f = \widetilde{f}_{n,\lambda} - \widehat{f}_{n,\lambda}$ and $\Delta W = W_n - W_n^\star$. Define $Rem_n = \widehat{f}_{n,\lambda} - f_0 - S_{n,\lambda}(f_0)$, where $S_{n,\lambda}(f_0) = \frac{1}{n}\sum_{i=1}^n \epsilon_i K_{X_i} - \mathcal{P}_\lambda f_0$, $\epsilon_i = Y_i - \dot{A}(f_0(X_i))$, $K = K^{f_0}$ and $\mathcal{P}_\lambda = \mathcal{P}_\lambda^{f_0}$. See [45] for definition of $K^{f_0}$ and $\mathcal{P}_\lambda^{f_0}$. It follows by [44, Theorem 3.4]) or [43] that $\|Rem_n\| = O_{P_{f_0}^n}(D_n)$ with $D_n = a_n + b_n$. Then

$$
\begin{aligned}
\widehat{f}_\nu &= V(\widehat{f}_{n,\lambda}, \varphi_\nu) = V(f_0 + S_{n,\lambda}(f_0) + Rem_n, \varphi_\nu) \\
&= f_\nu^0 + V(S_{n,\lambda}(f_0), \varphi_\nu) + V(Rem_n, \varphi_\nu) \\
&= \frac{1}{1 + \lambda\gamma_\nu} f_\nu^0 + V(\frac{1}{n}\sum_{i=1}^n \epsilon_i K_{X_i}, \varphi_\nu) + V(Rem_n, \varphi_\nu).
\end{aligned}
$$

It holds that $\Delta f = \sum_\nu (a_{n,\nu} - 1)\widehat{f}_\nu \varphi_\nu$, which leads to

$$
\begin{aligned}
&\|\Delta f\|_2^2 \\
={}& \sum_\nu (a_{n,\nu} - 1)^2 \widehat{f}_\nu^2 \\
={}& O_{P_{f_0}^n}(n^{-2}) + \sum_{\nu > m} \left(\frac{\tau_\nu^2}{\tau_\nu^2 + n(1 + \lambda\gamma_\nu)}\right)^2 \widehat{f}_\nu^2 \\
\lesssim{}& O_{P_{f_0}^n}(n^{-2}) + \sum_{\nu > m} \left(\frac{n^{-1}\tau_\nu^2}{1 + \lambda\gamma_\nu + n^{-1}\tau_\nu^2}\right)^2 \left(\frac{1}{(1 + \lambda\gamma_\nu)^2}|f_\nu^0|^2 + |V(n^{-1}\sum_{i=1}^n \epsilon_i K_{X_i}, \varphi_\nu)|^2 \right. \\
&\left. + |V(Rem_n, \varphi_\nu)|^2\right) \\
\equiv{}& O_{P_{f_0}^n}(n^{-2}) + I + II + III.
\end{aligned}
$$

Next we will analyze the terms $I, II, III$.

$$
\begin{aligned}
I &\leq \sum_{\nu > m} \frac{(h_* \nu)^{4m+2\beta}}{(1 + (h\nu)^{2m})^4}|f_\nu^0|^2 \\
&= (h_* h^{-1})^{4m+2\beta} \sum_{\nu > m} \frac{(h\nu)^{4m+2\beta}}{(1 + (h\nu)^{2m})^4}|f_\nu^0|^2 \\
&= (h_* h^{-1})^{4m+2\beta} h^{2m+\beta-1} \sum_{\nu > m} \frac{(h\nu)^{2m+\beta+1}}{(1 + (h\nu)^{2m})^4}|f_\nu^0|^2 \nu^{2m+\beta-1} \\
&= o((h_* h^{-1})^{4m+2\beta} h^{2m+\beta-1}) = o(n^{-2}h^{-(2m+\beta+1)}) = o(n^{-1}h^{-1/2}), \\
&\quad (\text{recall } h = n^{-a} \text{ with } a < \frac{2}{4m+2\beta+1})
\end{aligned}
$$

By [28] we have $E_{f_0}\{\epsilon^2|X\} = \ddot{A}(f_0(X))$, leading to

$$
\delta_{\nu\mu} = E\{\ddot{A}(f_0(X))\varphi_\nu(X)\varphi_\mu(X)\} = E_{f_0}\{\epsilon^2 \varphi_\nu(X)\varphi_\mu(X)\}.
$$

Define $A_{n,\nu} = \left(\frac{n^{-1}\tau_\nu^2}{1+\lambda\gamma_\nu+n^{-1}\tau_\nu^2}\right)^2 \frac{1}{(1+\lambda\gamma_\nu)^2}$. Then

$$
\begin{aligned}
II &= \sum_{\nu>m} A_{n,\nu} \left(\frac{1}{n}\sum_{i=1}^n \epsilon_i\varphi_\nu(X_i)\right)^2 \\
&= n^{-2}\sum_{i,j=1}^n \epsilon_i\epsilon_j \sum_{\nu>m} A_{n,\nu}\varphi_\nu(X_i)\varphi_\nu(X_j) \\
&= n^{-2}\left(\sum_{i=1}^n \epsilon_i^2\sum_{\nu>m}A_{n,\nu}\varphi_\nu(X_i)^2 + 2\sum_{i<j}\epsilon_i\epsilon_j\sum_{\nu>m}A_{n,\nu}\varphi_\nu(X_i)\varphi_\nu(X_j)\right).
\end{aligned}
\tag{A.35}
$$

Through direct examinations similar to the proof of Theorem 4.1, we have the following

$$
\begin{aligned}
&\sum_{i=1}^n \left(\epsilon_i^2\sum_{\nu>m}A_{n,\nu}\varphi_\nu(X_i)^2 - E_{f_0}\{\epsilon_i^2\sum_{\nu>m}A_{n,\nu}\varphi_\nu(X_i)^2\}\right) = O_{P_{f_0}^n}(n^{-3/2}h^{-(4m+2\beta+1)}), \\
&\sum_{i<j}\epsilon_i\epsilon_j\sum_{\nu>m}A_{n,\nu}\varphi_\nu(X_i)\varphi_\nu(X_j) = O_{P_{f_0}^n}(n^{-1}h^{-(8m+4\beta+1)/2}), \\
&E_{f_0}\{\epsilon_i^2\sum_{\nu>m}A_{n,\nu}\varphi_\nu(X_i)^2\} = \sum_{\nu>m}A_{n,\nu} \equiv c_h''/h.
\end{aligned}
\tag{A.36}
$$

Since $c_h'' \lesssim n^{-2}h^{-(4m+2\beta)}$, we have $II = O_{P_{f_0}^n}(n^{-3}h^{-(4m+2\beta+1)}) = o(n^{-1})$, implying that $II = o_{P_{f_0}^n}(n^{-1})$. Meanwhile,

$$
\begin{aligned}
III &\le n^{-2}\sum_{\nu>m}\frac{\nu^{4m+2\beta}}{(1+(h\nu)^{2m})^2}|V(Rem_n,\varphi_\nu)|^2 \\
&= n^{-2}h^{-(4m+2\beta)}\sum_{\nu>m}\frac{(h\nu)^{2m+2\beta}}{(1+(h\nu)^{2m})^2}|V(Rem_n,\varphi_\nu)|^2(h\nu)^{2m} \\
&\lesssim n^{-2}h^{-(4m+2\beta)}\|Rem_n\|^2 = O_{P_{f_0}^n}(n^{-2}h^{-(4m+2\beta)}D_n^2) = o_{P_{f_0}^n}(n^{-1}h^{-1/2}),
\end{aligned}
\tag{A.37}
$$

where the last equation follows by $D_n^2 = o(nh^{4m+2\beta-1/2})$ which can be verified by $a < \frac{2}{4m+2\beta+1}$ and direct calculations. In summary, $\|\Delta f\|_2^2 = o_{P_{f_0}^n}(n^{-1}h^{-1/2})$.

Meanwhile, by $\beta < m - 1/2$,

$$
\begin{aligned}
E\|\Delta W\|_2^2 &= O(n^{-3}) + \frac{1}{n}\sum_{\nu>m}\left(\frac{1}{\sqrt{1+\lambda\gamma_\nu+n^{-1}\tau_\nu^2}} - \frac{1}{\sqrt{1+\lambda\gamma_\nu}}\right)^2 \\
&\lesssim O(n^{-3}) + \frac{1}{n}\sum_{\nu>m}\frac{(h_*\nu)^{4m+2\beta}}{(1+(h\nu)^{2m})^3} \\
&= O(n^{-3}) + O(n^{-3}h^{-(4m+2\beta+1)}) = o(n^{-1}),
\end{aligned}
$$

implying that $\|\Delta W\|_2^2 = o_P(n^{-1})$. It follows by analysis of the terms $I, II, III$ and direct calcu-

lations that

$$
\begin{aligned}
& E\left(|\langle \Delta f, W_n^\star \rangle_2|^2 | \mathbf{D}_n\right) \\
= & \sum_\nu \frac{(a_{n,\nu}-1)^2}{n(1+\lambda\gamma_\nu)} \widehat{f}_\nu^2 \\
\lesssim & \ O_{P_{f_0}^n}(n^{-3}) + \sum_{\nu > m} \frac{1}{n(1+\lambda)} \left(\frac{n^{-1}\tau_\nu^2}{1+\lambda\gamma_\nu+n^{-1}\tau_\nu^2}\right)^2 \left(\frac{|f_\nu^0|^2}{(1+\lambda\gamma_\nu)^2} + |V(\frac{1}{n}\sum_{i=1}^n \epsilon_i K_{X_i}, \varphi_\nu)|^2\right. \\
& \left. + |V(Rem_n, \varphi_\nu)|^2\right) \\
= & \ O_{P_{f_0}^n}(n^{-3}) + n^{-3}h^{-(2m+\beta+1)} \sum_{\nu > m} \frac{(h\nu)^{2m+\beta+1}}{(1+(h\nu)^{2m})^5} |f_\nu^0|^2 \nu^{2m+\beta-1} \\
& + n^{-4}h^{-(4m+2\beta)} \sum_{\nu > m} \frac{(h\nu)^{4m+2\beta}}{(1+(h\nu)^{2m})^5} O_{P_{f_0}^n}(1) \\
& + n^{-3}h^{-(4m+2\beta)} \sum_{\nu > m} \frac{(h\nu)^{4m+2\beta}}{(1+(h\nu)^{2m})^3} |V(Rem_n, \varphi_\nu)|^2 \\
= & \ O_{P_{f_0}^n}(n^{-3}) + o(n^{-3}h^{-(2m+\beta+1)}) + O_{P_{f_0}^n}(n^{-4}h^{-(4m+2\beta+1)}) + O_{P_{f_0}^n}(n^{-3}h^{-(4m+2\beta)}D_n^2) \\
= & \ o_{P_{f_0}^n}(n^{-2}h^{-1}),
\end{aligned}
$$

(A.38)

where the last equation follows by $D_n^2 = o(nh^{4m+2\beta-1})$. This implies that, with $P_{f_0}^n$-probability approaching one,

$$
P\left(|\langle \Delta f, W_n^\star \rangle_2| \geq \varepsilon n^{-1}h^{-1/2}|\mathbf{D}_n\right) \leq \varepsilon,
$$

(A.39)

for arbitrary $\varepsilon > 0$.

Define $U_n^\star = (nh\|W_n^\star\|_2^2 - c_h)/\sqrt{2hd_h}$, where $c_h = \sum_{\nu \geq 1} \frac{h}{1+\lambda\gamma_\nu}$, $d_h = \sum_{\nu \geq 1} \frac{h}{(1+\lambda\gamma_\nu)^2}$. Similar to the proof of Lemma A.4, $U_n^\star \xrightarrow{d} Z \sim N(0,1)$. So $\|W_n^\star\|_2^2 = O_P((nh)^{-1})$. Therefore, we have proved the following: with $P_{f_0}^n$-probability or $P$-probability approaching one,

$$
\begin{aligned}
\|\Delta W\|_2^2 & \leq \varepsilon n^{-1}, \\
\|\Delta f\|_2^2 & \leq \varepsilon n^{-1}h^{-1/2}, \\
|\langle \Delta W, \Delta f \rangle_2| & \leq \|\Delta W\|_2 \|\Delta f\|_2 \leq \varepsilon n^{-1}h^{-1/4}, \\
|\langle \Delta W, W_n^\star \rangle_2| & \leq \|\Delta W\|_2 \|W_n^\star\|_2 \leq \varepsilon n^{-1}h^{-1/2},
\end{aligned}
$$

(A.40)

Note that

$$
\begin{aligned}
\|W - \widehat{f}_{n,\lambda}\|_2^2 & = \|W_n^\star + \Delta f + \Delta W\|_2^2 \\
& = \|W_n^\star\|_2^2 + \|\Delta W\|_2^2 + \|\Delta f\|_2^2 + 2\langle \Delta W, \Delta f \rangle_2 + 2\langle \Delta W, W_n^\star \rangle_2 + 2\langle \Delta f, W_n^\star \rangle_2.
\end{aligned}
$$

Together with (A.39) and (A.40), we have that, with $P_{f_0}^n$-probability approaching one,

$$
\begin{aligned}
P\left(\|W - \widehat{f}_{n,\lambda}\|_2^2 \le \frac{t^2}{nh} | \mathbf{D}_n\right) &\le P\left(\|W_n^\star\|_2^2 \le \frac{t^2}{nh} + 8\varepsilon n^{-1}h^{-1/2} | \mathbf{D}_n\right) + \varepsilon \\
&= P\left(U_n^\star \le \frac{t^2 - c_h}{\sqrt{2hd_h}} + \frac{8\varepsilon}{\sqrt{2d_h}}\right) + \varepsilon,
\end{aligned}
$$

and

$$
\begin{aligned}
P\left(\|W - \widehat{f}_{n,\lambda}\|_2^2 \le \frac{t^2}{nh} | \mathbf{D}_n\right) &\ge P\left(\|W_n^\star\|_2^2 \le \frac{t^2}{nh} - 8\varepsilon n^{-1}h^{-1/2} | \mathbf{D}_n\right) - \varepsilon \\
&= P\left(U_n^\star \le \frac{t^2 - c_h}{\sqrt{2hd_h}} - \frac{8\varepsilon}{\sqrt{2d_h}}\right) - \varepsilon.
\end{aligned}
$$

Since

$$
P\left(\|W^\star - \widehat{f}_{n,\lambda}\|_2^2 \le \frac{t^2}{nh} | \mathbf{D}_n\right) = P\left(U_n^\star \le \frac{t^2 - c_h}{\sqrt{2hd_h}}\right),
$$

we have

$$
\begin{aligned}
&\sup_{t \ge 0} |P_0\left(\sqrt{nh}\|f - \widehat{f}_{n,\lambda}\|_2 \le t | \mathbf{D}_n\right) - P_\star\left(\sqrt{nh}\|f - \widehat{f}_{n,\lambda}\|_2 \le t | \mathbf{D}_n\right)| \\
&= \sup_{t \ge 0} |P\left(\|W - \widehat{f}_{n,\lambda}\|_2^2 \le \frac{t^2}{nh} | \mathbf{D}_n\right) - P\left(\|W^\star - \widehat{f}_{n,\lambda}\|_2^2 \le \frac{t^2}{nh} | \mathbf{D}_n\right)| \\
&\le \sup_{t \ge 0} P\left(\frac{t^2 - c_h}{\sqrt{2hd_h}} - \frac{8\varepsilon}{\sqrt{2d_h}} \le U_n^\star \le \frac{t^2 - c_h}{\sqrt{2hd_h}} + \frac{8\varepsilon}{\sqrt{2d_h}}\right) + \varepsilon. \quad (A.41)
\end{aligned}
$$

It follows by Polya's theorem ([11]), the cdf of $U_n^\star$ uniformly converges to $\Phi(\cdot)$, the cdf of $Z \sim N(0,1)$. We have that, as $n \to \infty$,

$$
\begin{aligned}
\sup_{t \ge 0} \Big| &P\left(\frac{t^2 - c_h}{\sqrt{2hd_h}} - \frac{8\varepsilon}{\sqrt{2d_h}} \le U_n^\star \le \frac{t^2 - c_h}{\sqrt{2hd_h}} + \frac{8\varepsilon}{\sqrt{2d_h}}\right) \\
&- \left(\Phi\left(\frac{t^2 - c_h}{\sqrt{2hd_h}} + \frac{8\varepsilon}{\sqrt{2d_h}}\right) - \Phi\left(\frac{t^2 - c_h}{\sqrt{2hd_h}} - \frac{8\varepsilon}{\sqrt{2d_h}}\right)\right) \Big| \le \varepsilon.
\end{aligned}
$$

Since

$$
\sup_{t \ge 0} \left| \Phi\left(\frac{t^2 - c_h}{\sqrt{2hd_h}} + \frac{8\varepsilon}{\sqrt{2d_h}}\right) - \Phi\left(\frac{t^2 - c_h}{\sqrt{2hd_h}} - \frac{8\varepsilon}{\sqrt{2d_h}}\right) \right| \le c\varepsilon,
$$

where $c \ge \frac{16}{\sqrt{2d_h}}$ is a constant. Therefore, with $n \to \infty$,

$$
\sup_{t \ge 0} P\left(\frac{t^2 - c_h}{\sqrt{2hd_h}} - \frac{8\varepsilon}{\sqrt{2d_h}} \le U_n^\star \le \frac{t^2 - c_h}{\sqrt{2hd_h}} + \frac{8\varepsilon}{\sqrt{2d_h}}\right) \le (1 + c)\varepsilon. \quad (A.42)
$$

The result (5.2) then follows from (A.41) and (A.42).

To finish the rest of the proof, we show that (5.2) fails when $\frac{1}{2m+\beta+1/2} \le a \le \frac{8m+4\beta+2}{(8m+4\beta+1)(2m+\beta)}$. Define $c'_h = \sum_\nu \frac{h}{1 + \lambda\gamma_\nu + n^{-1}\tau_\nu^2}$ and $d'_h = \sum_\nu \frac{h}{(1 + \lambda\gamma_\nu + n^{-1}\tau_\nu^2)^2}$. Let $f_0 = 0$, then it holds that $I = 0$. The proof relies on the following decomposition

$$
\|W - \widehat{f}_{n,\lambda}\|_2^2 = \|W_n\|_2^2 + \|\Delta f\|_2^2 + 2\langle \Delta f, W_n \rangle_2.
$$

Similar to the proof of (A.38) and (A.39), one can show that

$$E\left(|\langle \Delta f, W_n \rangle_2|^2 | \mathbf{D}_n\right) = O_{P_{f_0}^n}(n^{-4}h^{-(4m+2\beta+1)}),$$

and hence, with $P_{f_0}^n$-probability approaching one,

$$P\left(|\langle \Delta f, W_n \rangle_2| \geq C_\varepsilon n^{-2}h^{-(4m+2\beta+1)/2} | \mathbf{D}_n\right) \leq \varepsilon, \tag{A.43}$$

where $C_\varepsilon > 0$ is a (possibly large) constant. By definitions of $c_h, c'_h, c''_h$ we have

$$c_h - c'_h - c''_h = \frac{h}{n} \sum_{\nu > m} \frac{\tau_\nu^2(2\lambda\gamma_\nu + \lambda^2\gamma_\nu^2 + n^{-1}\tau_\nu^2 + (n^{-1}\tau_\nu^2)(\lambda\gamma_\nu))}{(1 + \lambda\gamma_\nu + n^{-1}\tau_\nu^2)^2(1 + \lambda\gamma_\nu)^2}. \tag{A.44}$$

Keep in mind that, since $\frac{2}{4m+2\beta+1} \leq a \leq \frac{1}{2m+1}$, it holds that $nhD_n^2 = o(1)$. The proof proceeds in two cases.

**Case 1**: $\frac{2}{4m+2\beta+1} \leq a \leq \frac{1}{2m+\beta}$. In this case, it can be verified that $c''_h \asymp n^{-2}h^{-(4m+2\beta)}$. And it follows from (A.35), (A.36) and (A.37) that

$$\|\Delta f\|_2^2 = \frac{c''_h}{nh}(1 + o_{P_{f_0}^n}(1)).$$

By $a \leq \frac{1}{2m+\beta}$ and (A.44) we have

$$c_h - c'_h - c''_h \asymp n^{-1}h^{-(2m+\beta)}, \text{ and } c''_h \lesssim c_h - c'_h - c''_h.$$

Therefore,

$$\begin{aligned}
&P\left(\|W_n\|_2^2 + \|\Delta f\|_2^2 + 2\langle \Delta f, W_n \rangle_2 \leq \frac{c_h}{nh} | \mathbf{D}_n\right) \\
\geq \quad &P\left(U_n \leq \frac{c_h - c'_h - nh\|\Delta f\|_2 - 2nh\langle \Delta f, W_n \rangle_2}{\sqrt{2hd'_h}} | \mathbf{D}_n\right) \\
\geq \quad &P\left(U_n \leq \frac{(c_h - c'_h - c''_h)(1 + o_{P_{f_0}^n}(1)) - C_\varepsilon n^{-1}h^{-(4m+2\beta-1)/2}}{\sqrt{2hd'_h}} | \mathbf{D}_n\right) - \varepsilon \\
\geq \quad &P\left(U_n \leq \frac{(c_h - c'_h - c''_h)/2 - C_\varepsilon n^{-1}h^{-(4m+2\beta-1)/2}}{\sqrt{2hd'_h}}\right) - \varepsilon \\
&\text{(with } P_{f_0}^n\text{-probability approaching one)} \\
\geq \quad &1/2 + c - \varepsilon,
\end{aligned}$$

where $c > 0$ is a constant. The existence of $c$ in the last inequality follows by $n^{-1}h^{-(4m+2\beta-1)/2} = o(c_h - c'_h - c''_h)$, $U_n \xrightarrow{d} N(0,1)$, and the fact that $(c_h - c'_h - c''_h)/\sqrt{h} \asymp n^{-1}h^{-(2m+\beta+1/2)}$ which is greater than some fixed constant. Note that we should select $\varepsilon > 0$ to be small enough (at least smaller than $c/2$). Since $P(\sqrt{nh}\|W^\star - \widehat{f}_{n,\lambda}\|_2 \leq \sqrt{c_h} | \mathbf{D}_n) = P(U_n^\star \leq 0) \to 1/2$, the above analysis implies that, with $P_{f_0}^n$-probability approaching one, the left side of (5.2) is greater than a positive constant. So (5.2) does not hold.

**Case 2:** $\frac{1}{2m+\beta} < a < \frac{8m+4\beta+2}{(8m+4\beta+1)(2m+\beta)}$. In this case, $h_*/h \to \infty$. It follows from (A.35), (A.36), (A.37), (A.43), (A.44) that the following consequences immediately hold:

$$
\begin{aligned}
c_h'' &= \sum_{\nu > m} \left( \frac{n^{-1}\tau_\nu^2}{1 + \lambda + n^{-1}\tau_\nu^2} \right)^2 \frac{h}{(1+\lambda\gamma_\nu)^2} \gtrsim h/h_*, \\
\|\Delta f\|_2^2 &= \frac{c_h''}{nh}(1 + o_{P_{f_0}^n}(1)), \\
\langle \Delta f, W_n \rangle_2 &= o_{P_{f_0}^n}(\|\Delta f\|_2^2), \\
c_h - c_h' - c_h'' &\gtrsim h/h_*, \\
\frac{c_h - c_h' - c_h''}{\sqrt{h}} &\to \infty, \text{ as } n \to \infty.
\end{aligned}
$$

The proof of the above assertions rely on the specified range of $a$ and direct calculations. Therefore,

$$
\begin{aligned}
&P\left( \|W_n\|_2^2 + \|\Delta f\|_2^2 + 2\langle \Delta f, W_n \rangle_2 \le \frac{c_h}{nh} | \mathbf{D}_n \right) \\
={}& P\left( \|W_n\|_2^2 + \|\Delta f\|_2^2(1 + o_{P_{f_0}^n}) \le \frac{c_h}{nh} | \mathbf{D}_n \right) \\
={}& P\left( \|W_n\|_2^2 + \frac{c_h''}{nh}(1 + o_{P_{f_0}^n}(1)) \le \frac{c_h}{nh} | \mathbf{D}_n \right) \\
={}& P\left( U_n \le \frac{c_h - c_h' - c_h''(1 + o_{P_{f_0}^n}(1))}{\sqrt{2hd_h'}} \right) \to 1,
\end{aligned}
$$

where the last limit follows by (A.43). This would violate (5.2) based on arguments in **Case 1**. Proof is completed. $\qquad\square$

*Proof of Corollary 5.2.* Recall that $Rem_n = \widehat{f}_{n,\lambda} - f_0 - S_{n,\lambda}(f_0)$ satisfies $\|Rem_n\|_2 = O_{P_{f_0}^n}(D_n) = o_{P_{D_n}}$ with $nhD_n^2 = o(1)$. By the proof of Theorem 4.1, we have

$$
\begin{aligned}
&\|\widehat{f}_{n,\lambda} - f_0\|_2^2 \\
={}& \|S_{n,\lambda}(f_0)\|_2^2 + O_{P_{f_0}^n}(D_n^2) \\
={}& \frac{1}{n^2} \sum_{i=1}^n \epsilon_i^2 \langle K_{X_i}, K_{X_i} \rangle_2 + \frac{2}{n^2} \sum_{i<j} \epsilon_i \epsilon_j \langle K_{X_i}, K_{X_j} \rangle_2 + \frac{2}{n} \sum_{i=1}^n \epsilon_i \langle K_{X_i}, \mathcal{P}_\lambda f_0 \rangle_2 \\
&+ \|\mathcal{P}_\lambda f_0\|_2^2 + O_{P_{f_0}^n}(D_n^2) \\
={}& \frac{1}{nh} \int_0^\infty \frac{1}{(1+x^{2m})^2} dx + o(h^{2m+\beta-1}) + O_{P_{f_0}^n}(n^{-3/2}h^{-1} + n^{-1}h^{-1/2} + n^{-1/2}h^{(2m+\beta-1)/2}) \\
\le{}& (r_n^\star)^2 (1 + o_{P_{f_0}^n}(1)).
\end{aligned}
$$

Therefore, for any $\varepsilon > 0$, with $P_{f_0}^n$-probability approaching one,

$$
\|\widehat{f}_{n,\lambda} - f_0\|_2 \le (1+\varepsilon) r_n^\star.
$$

This implies that $P_{f_0}^n(f_0 \in R_n^\star(\varepsilon)) \to 1$ as $n \to \infty$.

Next we examine the posterior coverage of $R_n^\star(\varepsilon)$. Since $a$ satisfies (5.1), by Theorem 5.2,

$$P\left(R_n^\star(\varepsilon)|\mathbf{D}_n\right) - P_\star(R_n^\star(\varepsilon)) = o_{P_{f_0}^n}(1).$$

Since $\frac{1}{nh} \ll (r_n^\star)^2$, we have

$$P_\star(R_n^\star(\varepsilon)) = P\left(U_n^\star \leq \frac{(1+\varepsilon)nh(r_n^\star)^2 - c_h}{\sqrt{2hd_h}}\right) \to 1.$$

This shows that $P\left(R_n^\star(\varepsilon)|\mathbf{D}_n\right) = 1 + o_{P_{f_0}^n}(1)$, completing the proofs. $\qquad\square$

### A.4. $L^2$-diameter of $R_n^\omega(\alpha)$

Without additional restrictions, the $L^2$-diameter of $R_n^\omega(\alpha)$ in (4.2) is infinity. To see this, consider $f = \widetilde{f}_{n,\lambda} + \sum_{\nu=1}^{N} f_\nu \varphi_\nu$, where $f_\nu^2 = \frac{r_{\omega,n}(\alpha)^2}{N\omega_\nu}$ for $1 \leq \nu \leq N$. Then $f \in R_n^\omega(\alpha)$ since $\sum_{\nu=1}^{N} \omega_\nu f_\nu^2 = r_{\omega,n}(\alpha)^2$. However,

$$\|f\|_2^2 = \sum_{\nu=1}^{N} f_\nu^2 = \frac{r_{\omega,n}(\alpha)^2}{N} \sum_{\nu=1}^{N} \omega_\nu^{-1} \geq \frac{r_{\omega,n}(\alpha)^2}{N} \sum_{\nu=1}^{N} \nu = \frac{r_{\omega,n}(\alpha)^2(N+1)}{2}.$$

Letting $N \to \infty$, we can see that $\|f\|_2^2 \to \infty$. Therefore, the $L^2$-diameter of $R_n^\omega(\alpha)$ is infinity.

Next we investigate the $L^2$-diameter of $R_n^{\star\omega}(\alpha)$. For any $g, f \in R_n^{\star\omega}(\alpha)$, let $u = g - f \equiv \sum_{\nu=1}^{\infty} u_\nu \varphi_\nu$, and choose $J_n \sim n^{1/(2m+\beta)}$. It follows by Remark 4.1 in the revised manuscript that $r_{\omega,n}(\alpha) = O_{P_{f_0}}(n^{-1/2})$, and hence, $\|u\|_\omega \leq 2r_{\omega,n}(\alpha) = O_{P_{f_0}}(n^{-1/2})$. Then

$$
\begin{aligned}
\|u\|_2^2 &= \sum_{1 \leq \nu \leq J_n} u_\nu^2 + \sum_{\nu > J_n} u_\nu^2 \\
&= \sum_{1 \leq \nu \leq J_n} \omega_\nu u_\nu^2 \omega^{-1} + \sum_{\nu \geq J_n} \rho_\nu^{1+\frac{\beta-1}{2m}} u_\nu^2 \rho_\nu^{-(1+\frac{\beta-1}{2m})} \\
&\leq 4J_n \log(2J_n) r_{\omega,n}(\alpha)^2 + 4M J_n^{-(2m+\beta-1)} \\
&= O_{P_{f_0}}(n^{-\frac{2m+\beta-1}{2m+\beta}} \log n),
\end{aligned}
$$

indicating that the $L^2$-diameter of $R_n^{\star\omega}(\alpha)$ is $O_{P_{f_0}}(n^{-\frac{2m+\beta-1}{2(2m+\beta)}}\sqrt{\log n})$.

## Supplementary Document: Part II

We first establish a theoretical foundation including a Bayesian RKHS framework in Section A.5, and then prove Proposition A.1 in Section A.6.

### A.5. Some Preliminary Results

In this section, let us introduce some technical preliminaries. Using (2.7), for any $g = \sum_\nu g_\nu \varphi_\nu$, $\tilde{g} = \sum_\nu \tilde{g}_\nu \varphi_\nu \in S^m(\mathbb{I})$, we have $J(g, \tilde{g}) = \sum_{\nu \geq 1} g_\nu \tilde{g}_\nu \gamma_\nu$. It therefore holds that

$$J(\varphi_\nu, \varphi_\mu) = \gamma_\nu \delta_{\nu\mu}, \ \nu, \mu \geq 1. \tag{A.45}$$

This shows that

$$\|g\|_{U,V}^2 = \sum_{\nu \geq 1} g_\nu^2 (1 + \rho_\nu), \ J(g) = \sum_{\nu \geq 1} g_\nu^2 \gamma_\nu.$$

Since $\gamma_\nu \asymp 1 + \rho_\nu$, we can see that the $\|\cdot\|_{U,V}$-norm and $J^{1/2}$-norm are equivalent. By Sobolev embedding theorem ([1]) which implies that the supremum norm is "weaker" than the $\|\cdot\|_{U,V}$-norm, there exists an absolute constant $C_3 > 0$ s.t. for any $g \in S^m(\mathbb{I})$,

$$\|g\|_\infty \leq C_3 \sqrt{J(g)}. \tag{A.46}$$

For any $f, g, \tilde{g} \in S^m(\mathbb{I})$, define $V_f(g, \tilde{g}) = E\{\ddot{A}(f(X))g(X)\tilde{g}(X)\}$. In particular, $V_{f_0}(\cdot, \cdot) = V(\cdot, \cdot)$. Let $(\varphi_{f,\nu}, \rho_{f,\nu})$ be the eigen-system corresponding to the following ODE:

$$(-1)^m \varphi_{f,\nu}^{(2m)}(\cdot) = \rho_{f,\nu} \ddot{A}(f(\cdot))\pi(\cdot)\varphi_{f,\nu}(\cdot),$$
$$\varphi_{f,\nu}^{(j)}(0) = \varphi_{f,\nu}^{(j)}(1) = 0, j = m, m+1, \ldots, 2m-1. \tag{A.47}$$

It follows from [44, Proposition 2.2] that $(\varphi_{f,\nu}, \rho_{f,\nu})$ satisfy the properties stated in Proposition 2.1 with $V$ therein replaced by $V_f$. Let $\gamma_{f,\nu} = 1$ if $\nu = 1, 2, \ldots, m$; $= \rho_{f,\nu}$ if $\nu > m$. For any $g, \tilde{g} \in S^m(\mathbb{I})$ with $g = \sum_\nu g_\nu \varphi_{f,\nu}$ and $\tilde{g} = \sum_\nu \tilde{g}_\nu \varphi_{f,\nu}$, define $J_f(g, \tilde{g}) = \sum_\nu g_\nu \tilde{g}_\nu \gamma_{f,\nu}$. Define an inner product

$$\langle g, \tilde{g} \rangle_f = V_f(g, \tilde{g}) + \lambda J_f(g, \tilde{g}), \ g \in S^m(\mathbb{I}),$$

and let $\|\cdot\|_f$ be the corresponding norm. Let $\mathcal{P}_\lambda^f$ be a self-adjoint positive-definite operator from $S^m(\mathbb{I})$ to itself s.t. $\langle \mathcal{P}_\lambda^f g, \tilde{g} \rangle_f = \lambda J_f(g, \tilde{g})$ for any $g, \tilde{g} \in S^m(\mathbb{I})$. For convenience, define $\mathcal{P}_\lambda = \mathcal{P}_\lambda^{f_0}$. In particular,

$$J_{f_0}(g, \tilde{g}) = J(g, \tilde{g}), \ \langle g, \tilde{g} \rangle_{f_0} = \langle g, \tilde{g} \rangle, \ \|g\|_{f_0} = \|g\|.$$

For any constant $C$ with $C > \|f_0\|_\infty$, let $C_0, C_1, C_2$ be positive constants satisfying Assumption A1. Since $1/C_2 \leq \ddot{A}(z) \leq C_2$ if $|z| \leq 2C$ (Assumption A1), we get that for any $f \in \mathcal{F}(C)$ and $g \in S^m(\mathbb{I})$, (leading to that $C_2^{-1} \leq \ddot{A}(f(X)) \leq C_2$ a.s.)

$$C_2^{-2} V(g, g) \leq V_f(g, g) \leq C_2^2 V(g, g), \tag{A.48}$$

that is, $V_f$ is *uniformly equivalent* to $V$ for $f \in \mathcal{F}(C)$. This leads to

$$C_2^{-2} \frac{V(g)}{V(g) + U(g)} \leq \frac{V_f(g)}{V_f(g) + U(g)} \leq C_2^2 \frac{V(g)}{V(g) + U(g)}.$$

It follows from (A.48) and mapping principle (see [60, Theorem 5.3]) that

$$C_2^{-2} \rho_\nu \leq \rho_{f,\nu} \leq C_2^2 \rho_\nu, \text{ for any } \nu > m \text{ and } f \in \mathcal{F}(C).$$

The following lemma says that the norms $\|\cdot\|$ and $\|\cdot\|_f$ are equivalent.

**Lemma A.7.** *If* $0 < \lambda \leq \frac{1}{2C_2^2}$, *then for any* $f \in \mathcal{F}(C)$ *and* $g \in S^m(\mathbb{I})$,

$$\frac{1}{\sqrt{2}C_2} \|g\| \leq \|g\|_f \leq \sqrt{2}C_2 \|g\|,$$

$$\left(1 + \frac{C_2^2}{\rho_{m+1}}\right)^{-1} C_2^{-2} J(g) \leq J_f(g) \leq \left(1 + \frac{1}{\rho_{m+1}}\right) C_2^2 J(g).$$

*Proof of Lemma A.7.* For any $g \in S^m(\mathbb{I})$ with $g = \sum_\nu g_\nu \varphi_{f,\nu}$, we have

$$V_f(g) = \sum_{\nu \geq 1} g_\nu^2, \ U(g) = \sum_{\nu > m} g_\nu^2 \rho_{f,\nu}, \ J_f(g) = \sum_{\nu=1}^m g_\nu^2 + \sum_{\nu > m} g_\nu^2 \rho_{f,\nu}.$$

So, $J_f(g) \leq V_f(g) + U(g)$ and $U(g) \leq J_f(g)$. Therefore, it follows by (A.48) that

$$
\begin{aligned}
\|g\|_f^2 &= V_f(g) + \lambda J_f(g) \\
&\leq (1 + \lambda) V_f(g) + \lambda U(g) \\
&\leq (1 + \lambda) C_2^2 V(g) + \lambda J(g) \leq (1 + \lambda) C_2^2 (V(g) + \lambda J(g)) \leq 2C_2^2 \|g\|^2,
\end{aligned}
$$

where the last inequality is because $\lambda \leq \frac{1}{2C_2^2} < 1$.

On the other hand,

$$
\begin{aligned}
\|g\|_f^2 &= V_f(g) + \lambda J_f(g) \\
&\geq C_2^{-2} V(g) + \lambda U(g) \\
&\geq C_2^{-2} V(g) + \lambda (J(g) - V(g)) \\
&= (C_2^{-2} - \lambda) V(g) + \lambda J(g) \geq \frac{1}{2C_2^2} (V(g) + \lambda J(g)) = \frac{1}{2C_2^2} \|g\|^2.
\end{aligned}
$$

Meanwhile, $J_f(g) \leq V_f(g) + U(g) \leq C_2^2 V(g) + J(g)$. It can be shown that $V(g) + U(g) \leq (1 + 1/\rho_{m+1}) J(g)$. To see this, write $g = \sum_\nu g_\nu \varphi_\nu$. Then it follows by $1 + \rho_\nu \leq (1 + 1/\rho_{m+1}) \gamma_\nu$ that

$$V(g) + U(g) = \sum_\nu g_\nu^2 (1 + \rho_\nu) \leq (1 + 1/\rho_{m+1}) \sum_\nu g_\nu^2 \gamma_\nu = (1 + 1/\rho_{m+1}) J(g).$$

So $J_f(g) \leq (1 + 1/\rho_{m+1}) C_2^2 J(g)$.

Similarly, we have that $J(g) \leq V(g) + U(g) \leq C_2^2 V_f(g) + U(g)$. Write $g = \sum_\nu g_\nu \varphi_{f,\nu}$. Since $C_2^2 \rho_\nu \geq \rho_{f,\nu} \geq C_2^{-2} \rho_\nu \geq C_2^{-2} \rho_{m+1}$ for $\nu > m$, we have $1 + \rho_{f,\nu} \leq (1 + C_2^2/\rho_{m+1})\gamma_{f,\nu}$. So

$$
\begin{aligned}
V_f(g) + U(g) &= \sum_\nu g_\nu^2 (1 + \rho_{f,\nu}) \\
&\leq (1 + C_2^2/\rho_{m+1}) \sum_\nu g_\nu^2 \gamma_{f,\nu} = (1 + C_2^2/\rho_{m+1}) J_f(g).
\end{aligned}
$$

Therefore, $J_f(g) \geq (1 + C_2^2/\rho_{m+1})^{-1} C_2^{-2} J(g)$. Proof is completed. $\qquad \square$

The equivalence of $\|\cdot\|$ and $\|\cdot\|_f$ stated in Lemma A.7 leads to that $S^m(\mathbb{I})$ is a RKHS under $\langle \cdot, \cdot \rangle_f$ for any $f \in \mathcal{F}(C)$. Let $K^f(x, x')$ be the corresponding reproducing kernel function. In particular, define $K = K^{f_0}$ for simplicity. By [44, Proposition 2.1] (see an online supplement document therein for its proof) we have the following series representation.

**Proposition A.3.** *For any $f \in \mathcal{F}(C)$, $g \in S^m(\mathbb{I})$ and $x \in \mathbb{I}$, we have $\|g\|_f^2 = \sum_\nu |V_f(g, \varphi_{f,\nu})|^2 (1 + \lambda \gamma_{f,\nu})$, $K_x^f(\cdot) \equiv K^f(x, \cdot) = \sum_\nu \frac{\varphi_{f,\nu}(x)}{1 + \lambda \gamma_{f,\nu}} \varphi_{f,\nu}(\cdot)$, and $\mathcal{P}_\lambda^f \varphi_{f,\nu}(\cdot) = \frac{\lambda \gamma_{f,\nu}}{1 + \lambda \gamma_{f,\nu}} \varphi_{f,\nu}(\cdot)$.*

The following lemma demonstrates a uniform bound for the kernel $K^f$.

**Lemma A.8.** *It holds that*

$$
c_K(C) \equiv \sup_{f \in \mathcal{F}(C)} \sup_{0 < h \leq 1} \sup_{x \in \mathbb{I}} h^{1/2} \|K_x^f\|_f \leq c_m \sqrt{\frac{C_2}{\pi} + 1},
$$

*where $c_m > 0$ is a universal constant depending on $m$ only.*

*Proof of Lemma A.8.* For any $f \in \mathcal{F}(C)$, $g \in S^m(\mathbb{I})$ and $x \in \mathbb{I}$, it follows by [13, Lemma (2.11), pp. 54] that

$$
|\langle K_x^f, g \rangle_f| = |g(x)| \leq c_m h^{-1/2} \sqrt{\|g\|_{L^2}^2 + \lambda \|g^{(m)}\|_{L^2}^2},
$$

where $c_m > 0$ is a universal constant depending on $m$ only, and $\|\cdot\|_{L^2}$ denotes the usual $L^2$-norm. Since $\|g\|_{L^2}^2 \leq \frac{C_2}{\pi} V_f(g)$ and $\|g^{(m)}\|_{L^2}^2 = U(g) \leq J_f(g)$ (see proof of Lemma A.7 for the last inequality). Then

$$
|\langle K_x^f, g \rangle_f| \leq c_m \sqrt{\frac{C_2}{\pi} + 1} h^{-1/2} \|g\|_f,
$$

implying that $\|K_x^f\|_f \leq c_m \sqrt{\frac{C_2}{\pi} + 1} h^{-1/2}$. So $c_K(C) \leq c_m \sqrt{\frac{C_2}{\pi} + 1}$. $\qquad \square$

The lemma below directly comes from Lemma A.8, which relates the norms $\|\cdot\|_f$ and $\|\cdot\|_\infty$.

**Lemma A.9.** *For any $f \in \mathcal{F}(C)$ and $g \in S^m(\mathbb{I})$, $\|g\|_\infty \leq c_K(C) h^{-1/2} \|g\|_f$.*

Suppose that $(Y, X)$ follows model (2.1) based on $f$. The following conditional expectation can be found based on [28]:

$$
E_f\{Y|X\} = \dot{A}(f(X)). \tag{A.49}
$$

Let $g, g_k \in S^m(\mathbb{I})$ for $k = 1, 2, 3$. The Fréchet derivative of $\ell_{n,\lambda}$ can be identified as

$$
\begin{aligned}
D\ell_{n,\lambda}(g)g_1 &= \frac{1}{n}\sum_{i=1}^{n}(Y_i - \dot{A}(g(X_i)))\langle K_{X_i}^f, g_1\rangle_f - \langle \mathcal{P}_\lambda^f g, g_1\rangle_f \\
&\equiv \langle S_{n,\lambda}(g), g_1\rangle_f.
\end{aligned}
$$

Define $S_\lambda(g) = E_f\{S_{n,\lambda}(g)\}$. We also use $DS_\lambda$ and $D^2 S_\lambda$ to represent the second- and third-order Fréchet derivatives of $S_\lambda$. Note that $S_{n,\lambda}(\widehat{f}_{n,\lambda}) = 0$, and $S_{n,\lambda}(f)$ can be expressed as

$$
S_{n,\lambda}(f) = \frac{1}{n}\sum_{i=1}^{n}(Y_i - \dot{A}(f(X_i)))K_{X_i}^f - \mathcal{P}_\lambda^f f. \tag{A.50}
$$

The Fréchet derivatives of $S_{n,\lambda}$ and $DS_{n,\lambda}$ are denoted $DS_{n,\lambda}(g)g_1 g_2$ and $D^2 S_{n,\lambda}(g)g_1 g_2 g_3$. These derivatives can be explicitly written as

$$
\begin{aligned}
D^2\ell_{n,\lambda}(g)g_1 g_2 &\equiv DS_{n,\lambda}(g)g_1 g_2 \\
&= -\frac{1}{n}\sum_{i=1}^{n}\ddot{A}(g(X_i))g_1(X_i)g_2(X_i) - \langle \mathcal{P}_\lambda^f g_1, g_2\rangle_f, \\
D^3\ell_{n,\lambda}(g)g_1 g_2 g_3 &\equiv D^2 S_{n,\lambda}(g)g_1 g_2 g_3 \\
&= -\frac{1}{n}\sum_{i=1}^{n}\dddot{A}(g(X_i))g_1(X_i)g_2(X_i)g_3(X_i), \\
DS_\lambda(g)g_1 &= -E\{\ddot{A}(g(X))g_1(X)K_X^f\} - \mathcal{P}_\lambda^f g_1, \\
D^2 S_\lambda(g)g_1 g_2 &= -E\{\dddot{A}(g(X))g_1(X)g_2(X)K_X^f\}.
\end{aligned}
$$

Consider a function class

$$
\mathcal{G}(C) = \{g \in S^m(\mathbb{I}) : \|g\|_\infty \le 1, J(g,g) \le 2C_2^2 c_K(C)^{-2}h^{-2m+1}\}.
$$

Let $N(\varepsilon, \mathcal{G}(C), \|\cdot\|_\infty)$ be $\varepsilon$ packing number in terms of supremum norm. The following result can be found in [55].

**Lemma A.10.** *There exists a universal constant $c_0 > 0$ s.t. for any $\varepsilon > 0$,*

$$
\log N(\varepsilon, \mathcal{G}(C), \|\cdot\|_\infty) \le c_0(\sqrt{2}C_2 c_K(C)^{-1})^{1/m}h^{-\frac{2m-1}{2m}}\varepsilon^{-1/m}.
$$

In the future, for notational simplicity, we will simply drop $C$ from $c_K(C)$ and $\mathcal{G}(C)$ if there is no confusion.

For $r \ge 0$, define $\Psi(r) = \int_0^r \sqrt{\log(1 + \exp(x^{-1/m}))}dx$. For arbitrary $\varepsilon > 0$, define

$$
\begin{aligned}
A(h, \varepsilon) &= \frac{32\sqrt{6}}{\tau}\sqrt{2}C_2 c_K^{-1}c_0^m h^{-(2m-1)/2}\Psi\left(\frac{1}{2\sqrt{2}C_2}c_K c_0^{-m}h^{(2m-1)/2}\varepsilon\right) \\
&\quad + \frac{10\sqrt{24}\varepsilon}{\tau}\sqrt{\log\left(1 + \exp\left(2c_0((\sqrt{2}C_2)^{-1}c_K h^{(2m-1)/2}\varepsilon)^{-1/m}\right)\right)},
\end{aligned}
$$

where $\tau = \sqrt{\log 1.5} \approx 0.6368$. We have the following useful lemma.

**Lemma A.11.** *For any $f \in S^m(\mathbb{I})$, suppose that $\psi_{n,f}(z; g)$ is a measurable function defined upon $z = (y, x) \in \mathcal{Y} \times \mathbb{I}$ and $g \in \mathcal{G}$ satisfying $\psi_{n,f}(z; 0) = 0$ and the following Lipschitz continuity condition: for any $1 \le i \le n$ and $g_1, g_2 \in \mathcal{G}$,*

$$|\psi_{n,f}(Z_i; g_1) - \psi_{n,f}(Z_i; g_2)| \le c_K^{-1} h^{1/2} \|g_1 - g_2\|_\infty. \tag{A.51}$$

*Then for any constant $t \ge 0$ and $n \ge 1$,*

$$\sup_{f \in S^m(\mathbb{I})} P_f \left( \sup_{g \in \mathcal{G}} \|Z_{n,f}(g)\|_f > t \right) \le 2 \exp \left( -\frac{t^2}{B(h)^2} \right),$$

*where $B(h) = A(h, 2)$ and*

$$Z_{n,f}(g) = \frac{1}{\sqrt{n}} \sum_{i=1}^n [\psi_{n,f}(Z_i; g) K_{X_i}^f - E_f \{\psi_{n,f}(Z_i; g) K_{X_i}^f\}].$$

*Proof of Lemma A.11.* For any $f \in S^m(\mathbb{I})$ and $n \ge 1$, and any $g_1, g_2 \in \mathcal{G}$, we get that

$$\|(\psi_{n,f}(Z_i; g_1) - \psi_{n,f}(Z_i; g_2)) K_{X_i}^f\|_f$$
$$\le c_K^{-1} h^{1/2} \|g_1 - g_2\|_\infty c_K h^{-1/2} = \|g_1 - g_2\|_\infty.$$

By Theorem 3.5 of [38], for any $t > 0$, $P_f (\|Z_{n,f}(g_1) - Z_{n,f}(g_2)\|_f \ge t) \le 2 \exp \left( -\frac{t^2}{8\|g_1 - g_2\|_\infty^2} \right)$. Then by Lemma 8.1 in [24], we have

$$\|\|Z_{n,f}(g_1) - Z_{n,f}(g_2)\|_f\|_{\psi_2} \le \sqrt{24} \|g_1 - g_2\|_\infty,$$

where $\|\cdot\|_{\psi_2}$ denotes the Orlicz norm associated with $\psi_2(s) := \exp(s^2) - 1$. Recall $\tau = \sqrt{\log 1.5} \approx 0.6368$. Define $\phi(x) = \psi_2(\tau x)$. Then it can be shown by elementary calculus that $\phi(1) \le 1/2$, and for any $x, y \ge 1$, $\phi(x)\phi(y) \le \phi(xy)$. By a careful examination of the proof of Lemma 8.2, it can be shown that for any random variables $\xi_1, \ldots, \xi_l$,

$$\| \max_{1 \le i \le l} \xi_i \|_{\psi_2} \le \frac{2}{\tau} \psi_2^{-1}(l) \max_{1 \le i \le l} \|\xi_i\|_{\psi_2}. \tag{A.52}$$

Next we use a "chaining" argument. Let $T_0 \subset T_1 \subset T_2 \subset \cdots \subset T_\infty := \mathcal{G}$ be a sequence of finite nested sets satisfying the following properties:

- for any $T_q$ and any $s, t \in T_q$, $\|s - t\|_\infty \ge \varepsilon 2^{-q}$; each $T_q$ is "maximal" in the sense that if one adds any point in $T_q$, then the inequality will fail;

- the cardinality of $T_q$ is upper bounded by

$$\log |T_q| \le \log N(\varepsilon 2^{-q}, \mathcal{G}, \|\cdot\|_\infty)$$
$$\le c_0 (\sqrt{2} C_2 c_K^{-1})^{1/m} h^{-(2m-1)/(2m)} (\varepsilon 2^{-q})^{-1/m},$$

where $c_0 > 0$ is absolute constant;

- each element $t_{q+1} \in T_{q+1}$ is uniquely linked to an element $t_q \in T_q$ which satisfies $\|t_q - t_{q+1}\|_\infty \leq \varepsilon 2^{-q}$.

For arbitrary $s_{k+1}, t_{k+1} \in T_{k+1}$ with $\|s_{k+1} - t_{k+1}\|_\infty \leq \varepsilon$, choose two chains (both being of length $k+2$) $t_q$ and $s_q$ with $t_q, s_q \in T_q$ for $0 \leq q \leq k+1$. The ending points $s_0$ and $t_0$ satisfy

$$
\begin{aligned}
\|s_0 - t_0\|_\infty &\leq \sum_{q=0}^{k} [\|s_q - s_{q+1}\|_\infty + \|t_q - t_{q+1}\|_\infty] + \|s_{k+1} - t_{k+1}\|_\infty \\
&\leq 2\sum_{q=0}^{k} \varepsilon 2^{-q} + \varepsilon \leq 5\varepsilon,
\end{aligned}
$$

and hence, $\|\|Z_{n,f}(s_0) - Z_{n,f}(t_0)\|_f\|_{\psi_2} \leq 5\sqrt{24}\varepsilon$. It follows by the proof of Theorem 8.4 of [24] and (A.52) that

$$
\begin{aligned}
&\left\|\max_{s_{k+1}, t_{k+1} \in T_{k+1}} \|Z_{n,f}(s_{k+1}) - Z_{n,f}(t_{k+1}) - (Z_{n,f}(s_0) - Z_{n,f}(t_0))\|_f\right\|_{\psi_2} \\
&\leq 2\sum_{q=0}^{k} \left\|\max_{\substack{u \in T_{q+1}, v \in T_q \\ u,v \text{ link each other}}} \|Z_{n,f}(u) - Z_{n,f}(v)\|_f\right\|_{\psi_2} \\
&\leq \frac{4}{\tau}\sum_{q=0}^{k} \psi_2^{-1}(N(2^{-q-1}\varepsilon, \mathcal{G}, \|\cdot\|_\infty)) \\
&\quad \times \max_{\substack{u \in T_{q+1}, v \in T_q \\ u,v \text{ link each other}}} \|\|Z_{n,f}(u) - Z_{n,f}(v)\|_f\|_{\psi_2} \\
&\leq \frac{4\sqrt{24}}{\tau}\sum_{q=0}^{k} \sqrt{\log\left(1 + N(\varepsilon 2^{-q-1}, \mathcal{G}, \|\cdot\|_\infty)\right)}\varepsilon 2^{-q} \\
&\leq \frac{8\sqrt{24}}{\tau}\sum_{q=1}^{k+1} \sqrt{\log\left(1 + \exp\left(c_0(\sqrt{2}C_2 c_K^{-1})^{1/m} h^{-(2m-1)/(2m)}(\varepsilon 2^{-q})^{-1/m}\right)\right)}\varepsilon 2^{-q} \\
&\leq \frac{32\sqrt{6}}{\tau}\int_0^{\varepsilon/2} \sqrt{\log\left(1 + \exp\left(c_0(\sqrt{2}C_2 c_K^{-1})^{1/m} h^{-(2m-1)/(2m)} x^{-1/m}\right)\right)}dx \\
&= \frac{32\sqrt{6}}{\tau}\sqrt{2}C_2 c_K^{-1} c_0^m h^{-(2m-1)/2}\Psi\left(\frac{1}{2\sqrt{2}C_2} c_K c_0^{-m} h^{(2m-1)/2}\varepsilon\right).
\end{aligned}
$$

On the other hand,

$$
\begin{aligned}
\left\|\max_{\substack{u,v \in T_0 \\ \|u-v\|_\infty \leq 5\varepsilon}} \|Z_{n,f}(u) - Z_{n,f}(v)\|_f\right\|_{\psi_2} &\leq \frac{2}{\tau}\psi_2(|T_0|^2) \max_{\substack{u,v \in T_0 \\ \|u-v\|_\infty \leq 5\varepsilon}} \|\|Z_{n,f}(u) - Z_{n,f}(v)\|_f\|_{\psi_2} \\
&\leq \frac{2}{\tau}\psi_2^{-1}(N(\varepsilon, \mathcal{G}, \|\cdot\|_\infty)^2)(5\sqrt{24}\varepsilon).
\end{aligned}
$$

Therefore,

$$\left\| \max_{\substack{s,t \in T_{k+1} \\ \|s-t\|_\infty \leq \varepsilon}} \|Z_{n,f}(s) - Z_{n,f}(t)\|_f \right\|_{\psi_2}$$

$$\leq \frac{32\sqrt{6}}{\tau} \sqrt{2} C_2 c_K^{-1} c_0^m h^{-(2m-1)/2} \Psi \left( \frac{1}{2\sqrt{2}C_2} c_K c_0^{-m} h^{(2m-1)/2} \varepsilon \right)$$

$$+ \frac{2}{\tau} \psi_2^{-1} (N(\varepsilon, \mathcal{G}, \|\cdot\|_\infty)^2)(5\sqrt{24}\varepsilon)$$

$$\leq \frac{32\sqrt{6}}{\tau} \sqrt{2} c_2 c_K^{-1} c_0^m h^{-(2m-1)/2} \Psi \left( \frac{1}{2\sqrt{2}C_2} c_K c_0^{-m} h^{(2m-1)/2} \varepsilon \right)$$

$$+ \frac{10\sqrt{24}\varepsilon}{\tau} \sqrt{\log \left( 1 + \exp \left( 2c_0((\sqrt{2}C_2)^{-1} c_K h^{(2m-1)/2} \varepsilon)^{-1/m} \right) \right)}$$

$$= A(h, \varepsilon).$$

Now for any $g_1, g_2 \in \mathcal{G}$ with $\|g_1 - g_2\|_\infty \leq \varepsilon/2$. Let $k \geq 2$, hence, $2^{1-k} \leq 1 - \|g_1 - g_2\|_\infty/\varepsilon$. Since $T_k$ is "maximal", there exist $s_k, t_k \in T_k$ s.t. $\max\{\|g_1 - s_k\|_\infty, \|g_2 - t_k\|_\infty\} \leq \varepsilon 2^{-k}$. It is easy to see that $\|s_k - t_k\|_\infty \leq \varepsilon$. So

$$\|Z_{n,f}(g_1) - Z_{n,f}(g_2)\|_f$$

$$\leq \|Z_{n,f}(g_1) - Z_{n,f}(s_k)\|_f + \|Z_{n,f}(g_2) - Z_{n,f}(t_k)\|_f$$

$$+ \|Z_{n,f}(s_k) - Z_{n,f}(t_k)\|_f$$

$$\leq 4\sqrt{n}\varepsilon 2^{-k} + \max_{\substack{u,v \in T_k \\ \|u-v\|_\infty \leq \varepsilon}} \|Z_{n,f}(u) - Z_{n,f}(v)\|_f.$$

Therefore, letting $k \to \infty$ we get that

$$\left\| \sup_{\substack{g_1, g_2 \in \mathcal{G} \\ \|g_1-g_2\|_\infty \leq \varepsilon/2}} \|Z_{n,f}(g_1) - Z_{n,f}(g_2)\|_f \right\|_{\psi_2}$$

$$\leq 4\sqrt{n}\varepsilon 2^{-k}/\sqrt{\log 2} + \left\| \max_{\substack{u,v \in T_k \\ \|u-v\|_\infty \leq \varepsilon}} \|Z_{n,f}(u) - Z_{n,f}(v)\|_f \right\|_{\psi_2}$$

$$\leq 4\sqrt{n}\varepsilon 2^{-k}/\sqrt{\log 2} + A(h, \varepsilon) \to A(h, \varepsilon).$$

Taking $\varepsilon = 2$ in the above inequality, we get that

$$\left\| \sup_{\substack{g_1, g_2 \in \mathcal{G} \\ \|g_1-g_2\|_\infty \leq 1}} \|Z_{n,f}(g_1) - Z_{n,f}(g_2)\|_f \right\|_{\psi_2} \leq A(h, 2) = B(h).$$

By Lemma 8.1 in [24], we have

$$P_f \left( \sup_{g \in \mathcal{G}} \|Z_{n,f}(g)\|_f \geq t \right) \leq 2 \exp \left( -\frac{t^2}{B(h)^2} \right).$$

Note that the right hand side in the above does not depend on $f$. This completes the proof. □

Define
$$H^m(C) = \{f \in S^m(\mathbb{I}) : J(f) \le C^2/C_3^2\}.$$

It follows from (A.46) that for any $g \in H^m(C)$, $\|g\|_\infty \le C_3\sqrt{J(g)} \le C$, implying that $g \in \mathcal{F}(C)$. Thus, we have proved the following inclusion:

$$H^m(C) \subseteq \mathcal{F}(C). \tag{A.53}$$

It is easy to see that when $C > C_3\sqrt{J(f_0)}$, then $f_0 \in H^m(C)$, and hence, $f_0 \in \mathcal{F}(C)$.

**Lemma A.12.** *Suppose that Assumption A1 holds. For any constant $C$ satisfying $C > C_3\sqrt{J(f_0)}$, let $C_0, C_1, C_2$ be positive constants satisfying Assumption A1, and define*

$$b = \frac{C_2 C}{C_3}\sqrt{1 + \frac{1}{\rho_{m+1}}}. \tag{A.54}$$

*If $r, h, M$ are positives satisfying the following Rate Condition (H):*

*(i)* $(4C_2c_K^2 + 5)bh^{m-1/2} \le 2(\log 2)C_0 c_K$, $C_2^2 c_K bh^{m-1/2} \le 1/4$,
   $2bC_2 h^{m+1/2} \le c_K$,

*(ii)* $h^{1/2}r \le 1$,

*(iii)* $C_2 c_K^2 M^{1/2} rh^{-1/2} B(h) \le 1/6$,

*(iv)* $12C_0 C_2 c_K^4 (4C_1 + M)h^{-1}r(M^{1/2}rB(h) + C_2^{1/2} c_K^{-1}) \le 1/6$,

*then, for any $1 \le j \le s$, the following two results hold:*

*(a)*
$$\sup_{f \in H^m(C)} P_f\left(\|\widehat{f}_{n,\lambda} - f\|_f \ge \delta_n\right) \le 6\exp(-Mnhr^2),$$

   *where $\delta_n = 2bh^m + 24C_0 c_K(4C_1 + M)r$;*

*(b) if in addition, $c_K h^{-1/2}\delta_n < C$, then*

$$\sup_{f \in H^m(C)} P_f\left(\|\widehat{f}_{n,\lambda} - f - S_{n,\lambda}(f)\|_f > a_n + b_n\right) \le 8\exp(-Mnhr^2),$$

   *where*
$$a_n = C_2 c_K^2 M^{1/2} h^{-1/2} rB(h)\delta_n, \ \text{and} \ b_n = C_2^2 c_K h^{-1/2}\delta_n^2.$$

We remark that Part (b) of Lemma A.12 can be viewed as a uniform extension of the functional Bahadur representation established by [43, 44].

*Proof of Lemma A.12.* Let $f \in H^m(C)$ be the parameter based on which the data are drawn.

It is easy to see that

$$DS_\lambda(f)g = -E\{A(f(X))g(X)K_X^f\} - \mathcal{P}_\lambda^f g,$$

for any $g \in S^m(\mathbb{I})$. Therefore, for any $g, \tilde{g} \in S^m(\mathbb{I})$, $\langle DS_\lambda(f)g, \tilde{g}\rangle_f = -\langle g, \tilde{g}\rangle_f$, leading to $DS_\lambda(f) = -id$.

The proof of (a) is finished in two parts.

**Part I**: Define an operator mapping $S^m(\mathbb{I})$ to $S^m(\mathbb{I})$:

$$T_{1f}(g) = g + S_\lambda(f + g), \ g \in S^m(\mathbb{I}).$$

First observe that

$$\|S_\lambda(f)\|_f = \|\mathcal{P}_\lambda^f f\|_f = \sup_{\|g\|_f = 1} |\langle \mathcal{P}_\lambda^f f, g\rangle_f| \le \sqrt{\lambda J_f(f)} \le h^m b,$$

where the last inequality follows by Lemma A.7 and $f \in H^m(C)$. Let $r_{1n} = 2bh^m$. Let $\mathbb{B}(r_{1n}) = \{g \in S^m(\mathbb{I}) : \|g\|_f \le r_{1n}\}$ be the $r_{1n}$-ball. For any $g \in \mathbb{B}(r_{1n})$, using $DS_\lambda(f) = -id$ and $\|g\|_\infty \le c_K h^{-1/2} r_{1n} = 2c_K bh^{m-1/2} \le C$, it is easy to see that

$$\begin{aligned}
&\|T_{1f}(g)\|_f \\
\le\ & \|g + S_\lambda(f + g) - S_\lambda(f)\|_f + \|S_\lambda(f)\|_f \\
=\ & \|g + DS_\lambda(f)g + \int_0^1 \int_0^1 sD^2 S_\lambda(f + ss'g)ggds ds'\|_f + \|S_\lambda(f)\|_f \\
=\ & \|\int_0^1 \int_0^1 sD^2 S_\lambda(f + ss'g)ggds ds'\|_f + \|S_\lambda(f)\|_f \\
=\ & \|\int_0^1 \int_0^1 sE\{\ddot{A}(f(X) + ss'g(X))g(X)^2 K_X^f\}ds ds'\|_f + r_{1n}/2 \\
\le\ & C_2 c_K h^{-1/2} \int_0^1 \int_0^1 sE\{g(X)^2\}ds ds' + r_{1n}/2 \\
\le\ & C_2^2 c_K h^{-1/2}\|g\|_f^2/2 + r_{1n}/2 \\
\le\ & C_2^2 c_K h^{-1/2} r_{1n}^2/2 + r_{1n}/2 = C_2^2 c_K bh^{m-1/2} r_{1n} + r_{1n}/2 \le 3r_{1n}/4,
\end{aligned}$$

where the last step follows from the assumption $C_2^2 c_K bh^{m-1/2} \le 1/4$. Therefore, $T_{1f}$ maps $\mathbb{B}(r_{1n})$ to itself.

For any $g_1, g_2 \in \mathbb{B}(r_{1n})$, denote $g = g_1 - g_2$. Note that for any $0 \le s \le 1$, $\|g_2 + sg\|_f \le s\|g_1\|_f + (1-s)\|g_2\|_f \le r_{1n}$. By rate assumption we get that $\|g_2 + sg\|_\infty \le c_K h^{-1/2} r_{1n} = 2bc_K h^{m-1/2} < C$, and hence $|f(X) + s'(g_2(X) + sg(X))| \le 2C$ for any $s, s' \in [0, 1]$. By Taylor's expansion and

Lemma A.9 we have

$$\|T_{1f}(g_1) - T_{1f}(g_2)\|_f$$
$$= \|g_1 - g_2 + S_\lambda(f + g_1) - S_\lambda(f + g_2)\|_f$$
$$= \|g_1 - g_2 + \int_0^1 DS_\lambda(f + g_2 + sg)g ds\|_f$$
$$= \|\int_0^1 [DS_\lambda(f + g_2 + sg) - DS_\lambda(f)]g ds\|_f$$
$$= \|\int_0^1 \int_0^1 D^2 S_\lambda(f + s'(g_2 + sg))(g_2 + sg)g ds ds'\|_f$$
$$\leq \int_0^1 \int_0^1 \|E\{\ddot{A}(f(X) + s'(g_2(X) + sg(X)))(g_2(X) + sg(X))g(X)K_X^f\}\|_f ds ds'$$
$$\leq C_2 c_K h^{-1/2} \int_0^1 E\{|g_2(X) + sg(X)| \times |g(X)|\} ds$$
$$\leq C_2^2 c_K h^{-1/2} \int_0^1 \|g_2 + sg\|_f ds \times \|g\|_f$$
$$\leq 2C_2^2 c_K b h^{m-1/2}\|g_1 - g_2\|_f \leq \|g_1 - g_2\|_f/2.$$

This shows that $T_{1f}$ is a contraction mapping which maps $\mathbb{B}(r_{1n})$ into $\mathbb{B}(r_{1n})$. By contraction mapping theorem (see [39]), $T_{1f}$ has a unique fixed point $g' \in \mathbb{B}(r_{1n})$ satisfying $T_{1f}(g') = g'$. Let $f_\lambda = f + g'$. Then $S_\lambda(f_\lambda) = 0$ and $\|f_\lambda - f\|_f \leq r_{1n}$.

**Part II**: For any $f \in H^m(C)$, under (2.1) with $f$ being the truth, let $f_\lambda$ be the function obtained in **Part I** s.t. $\|f_\lambda - f\|_f \leq r_{1n}$, and hence, $\|f_\lambda - f\|_\infty \leq c_K h^{-1/2}\|f_\lambda - f\|_f \leq c_K h^{-1/2} r_{1n} \leq C/4$ so that $|f(X) + s(f_\lambda(X) - f(X))| \leq 2C$ a.s. for any $s \in [0, 1]$. It can be shown that for all $g_1, g_2 \in S^m(\mathbb{I})$,

$$|[DS_\lambda(f_\lambda) - DS_\lambda(f)]g_1 g_2|$$
$$= |\int_0^1 D^2 S_\lambda(f + s(f_\lambda - f))(f_\lambda - f)g_1 g_2 ds|$$
$$\leq \int_0^1 |E\{\ddot{A}(f(X) + s(f_\lambda - f)(X))(f_\lambda - f)(X)g_1(X)g_2(X)\}|ds$$
$$\leq C_2 E\{|f_\lambda(X) - f(X)| \cdot |g_1(X)g_2(X)|\}$$
$$\leq 2C_2^2 c_K b h^{m-1/2}\|g_1\|_f\|g_2\|_f \leq \|g_1\|_f\|g_2\|_f/2.$$

where the last inequality follows by $C_2^2 c_K b h^{m-1/2} \leq 1/4$. Together with the fact $DS_\lambda(f) = -id$, we get that the operator norm $\|DS_\lambda(f_\lambda) + id\|_{\text{operator}} \leq 1/2$. This implies that $DS_\lambda(f_\lambda)$ is invertible with operator norm within $[1/2, 3/2]$, and hence, $\|DS_\lambda(f_\lambda)^{-1}\|_{\text{operator}} \leq 2$.

Define an operator

$$T_{2f}(g) = g - [DS_\lambda(f_\lambda)]^{-1} S_{n,\lambda}(f_\lambda + g), \ g \in S^m(\mathbb{I}).$$

Rewrite $T_{2f}$ as

$$
\begin{aligned}
T_{2f}(g) \;=\; & -DS_\lambda(f_\lambda)^{-1}[DS_{n,\lambda}(f_\lambda)g - DS_\lambda(f_\lambda)g] \\
& -DS_\lambda(f_\lambda)^{-1}[S_{n,\lambda}(f_\lambda + g) - S_{n,\lambda}(f_\lambda) - DS_{n,\lambda}(f_\lambda)g] \\
& -DS_\lambda(f_\lambda)^{-1}S_{n,\lambda}(f_\lambda).
\end{aligned}
$$

Denote the above three terms by $I_{1f}, I_{2f}, I_{3f}$, respectively.

For any $1 \le i \le n$, let

$$
R_i = (Y_i - \dot{A}(f_\lambda(X_i)))K^f_{X_i} - E_f\{(Y - \dot{A}(f_\lambda(X)))K^f_X\}.
$$

Since $E_f\{Y - \dot{A}(f(X))|X\} = 0$, it can be shown that for some (random) $s \in [0, 1]$,

$$
\begin{aligned}
& \|E_f\{(Y - \dot{A}(f_\lambda(X)))K^f_X\}\|_f \\
=\; & \sup_{\|g\|_f = 1} |\langle E_f\{(Y - \dot{A}(f_\lambda(X)))K^f_X\}, g\rangle_f| \\
=\; & \sup_{\|g\|_f = 1} |E_f\{(Y - \dot{A}(f_\lambda(X)))g(X)\}| \\
=\; & \sup_{\|g\|_f = 1} |E_f\{(\dot{A}(f_\lambda(X)) - \dot{A}(f(X)))g(X)\}| \\
=\; & \sup_{\|g\|_f = 1} \left| E_f\left\{\ddot{A}(f(X))(f_\lambda(X) - f(X))g(X)\right\} \right. \\
& \left. +\frac{1}{2}E_f\left\{\dddot{A}(f(X) + s(f_\lambda(X) - f(X)))(f_\lambda(X) - f(X))^2 g(X)\right\} \right| \\
=\; & \sup_{\|g\|_f = 1} \left| \langle f_\lambda - f, g\rangle_f \right. \\
& \left. +\frac{1}{2}E_f\left\{\dddot{A}(f(X) + s(f_\lambda(X) - f(X)))(f_\lambda(X) - f(X))^2 g(X)\right\} \right| \\
\le\; & \|f_\lambda - f\|_f + \frac{C_2}{2}E\{(f_\lambda(X) - f(X))^2|g(X)|\} \\
\le\; & \|f_\lambda - f\|_f + \frac{1}{2}C_2^2 c_K h^{-1/2}\|f_\lambda - f\|_f^2 \\
\le\; & r_{1n} + C_2^2 c_K b h^{m-1/2} r_{1n} \le 5r_{1n}/4.
\end{aligned}
$$

Therefore,

$$
\begin{aligned}
\|R_i\|_f \;\le\; & c_K h^{-1/2}|Y_i - \dot{A}(f_\lambda(X_i))| + 5r_{1n}/4 \\
\le\; & c_K h^{-1/2}\left(|Y_i - \dot{A}(f(X_i))| + 2C_2 c_K b h^{m-1/2}\right) + 5r_{1n}/4,
\end{aligned}
$$

which leads to that

$$
E\left\{\exp\left(\frac{\|R_i\|_f}{C_0 c_K h^{-1/2}}\right)\right\} \le C_1 \exp\left(\frac{(4C_2 c_K^2 + 5)b h^{m-1/2}}{2C_0 c_K}\right) \le 2C_1,
$$

where the last inequality follows by condition

$$(4C_2c_K^2 + 5)bh^{m-1/2} \le 2(\log 2)C_0c_K.$$

Let $\delta = hr/(2C_0c_K)$. Recall the condition $h^{1/2}r \le 1$ which implies $\delta \le (2C_0c_Kh^{-1/2})^{-1}$. Therefore,

$$E\{\exp(2\delta\|R_i\|_f)\} \le E\{\exp(\|R_i\|_f/(C_0c_Kh^{-1/2}))\} \le 2C_1.$$

Moreover, $\|R_i\|_f^2 \le 8C_0^2c_K^2h^{-1}\exp(\|R_i\|_f/(2C_0c_Kh^{-1/2}))$, which leads to that

$$
\begin{aligned}
&E\{\exp(\delta\|R_i\|_f) - 1 - \delta\|R_i\|_f\} \\
\le{} & E\{(\delta\|R_i\|_f)^2\exp(\delta\|R_i\|_f)\} \\
\le{} & 8C_0^2c_K^2h^{-1}\delta^2E\left\{\exp\left(\left(\delta + \frac{1}{2C_0c_Kh^{-1/2}}\right)\|R_i\|_f\right)\right\} \\
\le{} & 16C_0^2C_1c_K^2h^{-1}\delta^2.
\end{aligned}
$$

It follows by Theorem 3.2 of [38] that, for $L(M) := 2C_0c_K(4C_1 + M)$,

$$
\begin{aligned}
P_f\left(\|\sum_{i=1}^n R_i\|_f \ge L(M)nr\right) &\le 2\exp\left(-L(M)\delta nr + 16C_0^2C_1c_K^2nh^{-1}\delta^2\right) \\
&= 2\exp(-Mnhr^2),
\end{aligned}
\tag{A.55}
$$

We note that the right hand side in the above inequality does not depend on $f$. It is easy to see that $S_{n,\lambda}(f_\lambda) = S_{n,\lambda}(f_\lambda) - S_\lambda(f_\lambda) = \frac{1}{n}\sum_{i=1}^n R_i$. Let

$$\mathcal{E}_{n,1} = \{\|S_{n,\lambda}(f_\lambda)\|_f \le L(M)r\},$$

then $\sup_{f \in H^m(C)} P_f(\mathcal{E}_{n,1}^c) \le 2\exp(-Mnhr^2)$.

Define

$$\psi_{n,f}^{(1)}(X_i; g) = [C_2c_K]^{-1}h^{1/2}\ddot{A}(f_\lambda(X_i))g(X_i), \ i = 1, \dots, n,$$

and $Z_{n,f}^{(1)}(g) = \frac{1}{\sqrt{n}}\sum_{i=1}^n[\psi_{n,f}^{(1)}(X_i; g)K_{X_i}^f - E_f\{\psi_{n,f}^{(1)}(X_i; g)K_{X_i}^f\}]$. It follows by Lemma A.11 that $\sup_{f \in H^m(C)} P_f(\mathcal{E}_{n,2}^c) \le 2\exp(-Mnhr^2)$, where $\mathcal{E}_{n,2} = \{\sup_{g \in \mathcal{G}}\|Z_{n,f}^{(1)}(g)\|_f \le \sqrt{Mnhr^2}B(h)\}$.

For any $g \in S^m(\mathbb{I})\backslash\{0\}$, let $\bar{g} = g/d_n'$, where $d_n' = c_Kh^{-1/2}\|g\|_f$. It follows by Lemma A.7 that

$$\|\bar{g}\|_\infty \le c_Kh^{-1/2}\|\bar{g}\|_f = c_Kh^{-1/2}\|g\|_f/d_n' = 1, \text{ and}$$

$$
\begin{aligned}
J(\bar{g}, \bar{g}) &= d_n'^{-2}J(g, g) \\
&= h^{-2m}\frac{\lambda J(g, g)}{c_K^2h^{-1}\|g\|_f^2} \le h^{-2m}\frac{\|g\|^2}{c_K^2h^{-1}\|g\|_f^2} \le 2C_2^2c_K^{-2}h^{-2m+1}.
\end{aligned}
$$

Therefore, $\bar{g} \in \mathcal{G}$. Consequently, on $\mathcal{E}_{n,2}$, for any $g \in S^m(\mathbb{I}) \backslash \{0\}$, we get $\|Z_{n,f}^{(1)}(\bar{g})\|_f \leq \sqrt{Mnhr^2}B(h)$, which leads to that

$$
\begin{aligned}
& \|DS_{n,\lambda}(f_\lambda)g - DS_\lambda(f_\lambda)g\|_f \\
= \ & \frac{1}{n}\|\sum_{i=1}^{n}[\ddot{A}(f_\lambda(X_i))g(X_i)K_{X_i}^f - E_f\{\ddot{A}(f_\lambda(X_i))g(X_i)K_{X_i}^f\}]\|_f \\
\leq \ & C_2 c_K^2 M^{1/2} r h^{-1/2} B(h) \|g\|_f \leq \|g\|_f/6,
\end{aligned} \tag{A.56}
$$

where the last inequality follows by condition $C_2 c_K^2 M^{1/2} r h^{-1/2} B(h) \leq 1/6$. Note that the above inequality also holds for $g = 0$.

Next we define $T_{3f}(g) = S_{n,\lambda}(f_\lambda + g) - S_{n,\lambda}(f_\lambda) - DS_{n,\lambda}(f_\lambda)g$. Let $r_{2n} = 6L(M)r$. For any $g \in \mathcal{G}$ and $1 \leq i \leq n$, define $\tilde{\psi}_{n,i}(g) = |g(X_i)|$, and let $\tilde{Z}_n(g) = \frac{1}{\sqrt{n}}\sum_{i=1}^{n}[\tilde{\psi}_{n,i}(g) - E\{\tilde{\psi}_{n,i}(g)\}]$. It is easy to see that for any $g_1, g_2 \in \mathcal{G}$, $|\tilde{\psi}_{n,i}(g_1) - \tilde{\psi}_{n,i}(g_2)| \leq \|g_1 - g_2\|_\infty$. Following the proof of Lemma A.11 it can be shown that for any $t \geq 0$,

$$
P\left(\sup_{g \in \mathcal{G}} |\tilde{Z}_n(g)| \geq t\right) \leq 2\exp\left(-\frac{t^2}{B(h)^2}\right),
$$

and hence, we get that $P(\mathcal{E}_{n,3}^c) \leq 2\exp(-Mnhr^2)$, where

$$
\mathcal{E}_{n,3} = \{\sup_{g \in \mathcal{G}} |\tilde{Z}_n(g)| \leq \sqrt{Mnhr^2}B(h)\}.
$$

On $\mathcal{E}_{n,2} \cap \mathcal{E}_{n,3}$, for any $g_1, g_2 \in \mathbb{B}(r_{2n})$ (with $g_1 \neq g_2$) and letting $g = g_1 - g_2$ (and hence $\|g_2 + sg\|_\infty \leq c_K h^{-1/2} r_{2n} \leq C/4$ for any $s \in [0,1]$), together with $\|f_\lambda - f\|_\infty \leq C/4$, we have

$$
\begin{aligned}
& \|T_{3f}(g_1) - T_{3f}(g_2)\|_f \\
= \ & \|S_{n,\lambda}(f_\lambda + g_1) - S_{n,\lambda}(f_\lambda + g_2) - DS_{n,\lambda}(f_\lambda)g\|_f \\
= \ & \|\int_0^1 \int_0^1 D^2 S_{n,\lambda}(f_\lambda + s'(g_2 + sg))(g_2 + sg)g\,ds\,ds'\|_f \\
\leq \ & \int_0^1 \int_0^1 \|D^2 S_{n,\lambda}(f_\lambda + s'(g_2 + sg))(g_2 + sg)g\|_f\,ds\,ds' \\
\leq \ & \int_0^1 \int_0^1 \|\frac{1}{n}\sum_{i=1}^{n}\dddot{A}(f_\lambda(X_i) + s'(g_2(X_i) + sg(X_i))) \\
& (g_2(X_i) + sg(X_i))g(X_i)K_{X_i}^f\|_f\,ds\,ds' \\
\leq \ & \int_0^1 \int_0^1 \frac{C_2}{n}\sum_{i=1}^{n} \|g_2 + sg\|_\infty \times |g(X_i)| \times \|K_{X_i}^f\|_f\,ds\,ds' \\
\leq \ & \frac{C_2(c_K h^{-1/2})^2 r_{2n}}{n}\sum_{i=1}^{n} |g(X_i)| \\
= \ & \frac{C_2(c_K h^{-1/2})^3 r_{2n}}{n}\left(\sum_{i=1}^{n} \tilde{\psi}_{n,i}(\bar{g})\right)\|g\|_f,
\end{aligned} \tag{A.57}
$$

where $\bar{g} = g/(c_K h^{-1/2} \|g\|_f)$. Recalling the previous arguments we get $\bar{g} \in \mathcal{G}$. It can be shown by Cauchy-Schwartz inequality that

$$E\{\tilde{\psi}_{n,i}(\bar{g})\} = \frac{E\{|g(X_i)|\}}{c_K h^{-1/2} \|g\|_f} \leq \frac{C_2^{1/2} V_f(g,g)^{1/2}}{c_K h^{-1/2} \|g\|_f} \leq C_2^{1/2} c_K^{-1} h^{1/2}.$$

Since $\mathcal{E}_{n,3}$ implies $|\tilde{Z}_n(\bar{g})| \leq \sqrt{Mnhr^2} B(h)$, we get that

$$\frac{1}{n} \sum_{i=1}^{n} \tilde{\psi}_{n,i}(\bar{g}) \leq \sqrt{Mhr^2} B(h) + C_2^{1/2} c_K^{-1} h^{1/2}.$$

Therefore, (A.57) has upper bound

$$
\begin{aligned}
(A.57) &\leq C_2(c_K h^{-1/2})^3 r_{2n} (\sqrt{Mhr^2} B(h) + C_2^{1/2} c_K^{-1} h^{1/2}) \|g\|_f \\
&= 12 C_0 C_2 c_K^4 (4C_1 + M) h^{-1} r (M^{1/2} r B(h) + C_2^{1/2} c_K^{-1}) \|g\|_f \\
&\leq \|g_1 - g_2\|_f / 6,
\end{aligned}
\tag{A.58}
$$

where the last inequality follows by condition

$$12 C_0 C_2 c_K^4 (4C_1 + M) h^{-1} r (M^{1/2} r B(h) + C_2^{1/2} c_K^{-1}) \leq 1/6.$$

Taking $g_2 = 0$ in (A.58) we get that $\|T_{3f}(g_1)\|_f \leq \|g_1\|_f / 6$ for any $g_1 \in \mathbb{B}(r_{2n})$. Therefore, it follows by (A.56) that, for any $f \in H^m(C)$, on $\mathcal{E}_n := \mathcal{E}_{n,1} \cap \mathcal{E}_{n,2} \cap \mathcal{E}_{n,3}$ and for any $g \in \mathbb{B}(r_{2n})$,

$$\|T_{2f}(g)\|_f \leq 2(\|g\|_f / 6 + \|g\|_f / 6 + r_{2n}/6) \leq 2(r_{2n}/6 + r_{2n}/6 + r_{2n}/6) = r_{2n}.$$

Meanwhile, for any $g_1, g_2 \in \mathbb{B}(r_{2n})$, replacing $g$ by $g_1 - g_2$ in (A.56), together with (A.57) and (A.58), we get that

$$\|T_{2f}(g_1) - T_{2f}(g_2)\|_f \leq 2(\|g_1 - g_2\|_f / 6 + \|g_1 - g_2\|_f / 6) = 2\|g_1 - g_2\|_f / 3.$$

Therefore, for any $f \in H^m(C)$, on $\mathcal{E}_n$, $T_{2f}$ is a contraction mapping from $\mathbb{B}(r_{2n})$ to itself. By contraction mapping theorem, there exists uniquely an element $g'' \in \mathbb{B}(r_{2n})$ s.t. $T_{2f}(g'') = g''$. Let $\hat{f}_{n,\lambda} = f_\lambda + g''$. Clearly, $S_{n,\lambda}(\hat{f}_{n,\lambda}) = 0$, and hence, $\hat{f}_{n,\lambda}$ is the maximizer of $\ell_{n,\lambda}$; see (3.6). So we get that, on $\mathcal{E}_n$, $\|\hat{f}_{n,\lambda} - f\|_f \leq \|f_\lambda - f\|_f + \|\hat{f}_{n,\lambda} - f_\lambda\|_f \leq r_{1n} + r_{2n} = 2bh^m + 6L(M)r$. The desired conclusion follows by the trivial fact: $\sup_{f \in H^m(C)} P_f(\mathcal{E}_n^c) \leq 6 \exp(-Mnhr^2)$. Proof of (a) is completed.

Next we show (b).

For any $f \in H^m(C)$, let $\hat{f}_{n,\lambda}$ be the penalized MLE of $f$ obtained by (3.6). Let $g_n = \hat{f}_{n,\lambda} - f$, $\delta_n = 2bh^m + 6L(M)r$, $d'_n = c_K h^{-1/2} \delta_n$, and for $g \in \mathcal{G}$ define

$$\psi_{n,f}^{(2)}(X_i; g) = c_K^{-1} h^{1/2} [C_2 d'_n]^{-1} (\dot{A}(f(X_i) + d'_n g(X_i)) - \dot{A}(f(X_i))).$$

It can be seen that for any $g_1, g_2 \in \mathcal{G}$, by $\delta'_n = c_K h^{-1/2} \delta_n < C$, we have

$$|\psi^{(2)}_{n,f}(X_i; g_1) - \psi^{(2)}_{n,f}(X_i; g_2)|$$
$$\leq c_K^{-1} h^{1/2} [C_2 d'_n]^{-1} C_2 d'_n \|g_1 - g_2\|_\infty = c_K^{-1} h^{1/2} \|g_1 - g_2\|_\infty.$$

Let $\mathcal{E}_{n,4} = \{\sup_{g \in \mathcal{G}} \|Z^{(2)}_{n,f}(g)\|_f \leq \sqrt{Mnhr^2} B(h)\}$, where

$$Z^{(2)}_{n,f}(g) = \frac{1}{\sqrt{n}} \sum_{i \in I_j} [\psi^{(2)}_{n,f}(X_i; g) K^f_{X_i} - E^X_f \{\psi^{(2)}_{n,f}(X; g) K^f_X\}],$$

$E^X_f$ denotes the expectation with respect to $X$ (under $P_f$). It follows by Lemma A.11 that $\sup_{f \in H^m(C)} P_f(\mathcal{E}^c_{n,4}) \leq 2 \exp(-Mnhr^2)$.

On $\widetilde{\mathcal{E}}_n := \mathcal{E}_n \cap \mathcal{E}_{n,4}$, we have $\|g_n\|_f \leq \delta_n$. Let $\bar{g} = g_n/d'_n$. Clearly, $\bar{g} \in \mathcal{G}$. Then we get that

$$\|S_{n,\lambda}(f + g_n) - S_{n,\lambda}(f) - (S_\lambda(f + g_n) - S_\lambda(f))\|_f$$
$$= \frac{1}{n} \| \sum_{i=1}^n [(\dot{A}(f(X_i) + g_n(X_i)) - \dot{A}(f(X_i))) K^f_{X_i}$$
$$- E^X_f \{(\dot{A}(f(X) + g_n(X)) - \dot{A}(f(X))) K^f_X\}] \|_f$$
$$= \frac{1}{n} \| \sum_{i \in I_j} [(\dot{A}(f(X_i) + d'_n \bar{g}(X_i)) - \dot{A}(f(X_i))) K^f_{X_i}$$
$$- E^X_f \{(\dot{A}(f(X) + d'_n \bar{g}(X)) - \dot{A}(f(X))) K^f_X\}] \|_f$$
$$= \frac{C_2 c_K h^{-1/2} d'_n}{n} \| \sum_{i \in I_j} [\psi^{(2)}_{n,f}(X_i; \bar{g}) K^f_{X_i} - E^X_f \{\psi^{(2)}_{n,f}(X; \bar{g}) K^f_X\}] \|_f$$
$$= \frac{C_2 c_K h^{-1/2} d'_n}{\sqrt{n}} \|Z^{(2)}_{n,f}(\bar{g})\|_f \leq C_2 c_K^2 M^{1/2} h^{-1/2} r B(h) \delta_n = a_n. \tag{A.59}$$

It is easy to show that

$$\| \int_0^1 \int_0^1 s D^2 S_\lambda(f + ss' g_n) g_n g_n ds ds' \|_f$$
$$= \| \int_0^1 \int_0^1 s E^X_f \{\ddot{A}(f(X) + ss' g_n(X)) g_n(X)^2 K_X\} ds ds' \|_f$$
$$\leq C_2 c_K h^{-1/2} \int_0^1 \int_0^1 s E^X_f \{g_n(X)^2\} ds ds'$$
$$\leq C_2^2 c_K h^{-1/2} \|g_n\|_f^2 \leq C_2^2 c_K h^{-1/2} \delta_n^2 = b_n. \tag{A.60}$$

Since $S_{n,\lambda}(f + g_n) = 0$ and $DS_\lambda(f) = -id$, from (A.59) and (A.60) we have on $\widetilde{\mathcal{E}}_n$,

$$a_n \geq \|S_{n,\lambda}(f) + DS_\lambda(f) g_n + \int_0^1 \int_0^1 s D^2 S_\lambda(f + ss' g_n) g_n g_n ds ds' \|_f$$
$$= \|S_{n,\lambda}(f) - g_n + \int_0^1 \int_0^1 s D^2 S_\lambda(f + ss' g_n) g_n g_n ds ds' \|_f$$
$$\geq \|S_{n,\lambda}(f) - g_n\|_f - \| \int_0^1 \int_0^1 s D^2 S_\lambda(f + ss' g_n) g_n g_n ds ds' \|_f,$$

which implies that

$$\|\widehat{f}_{n,\lambda} - f - S_{n,\lambda}(f)\|_f \le a_n + b_n.$$

Since $\sup_{f \in H^m(C)} P_f(\widetilde{\mathcal{E}}_n^c) \le 8 \exp(-Mnhr^2)$, proof of (b) is completed. $\qquad\square$

### *A.6. Proof of Proposition A.1*

The proof requires the following result.

**Proposition A.4.** *(An initial contraction rate) Under Assumption A1, if $r_n = o(h^{3/2})$, $h^{1/2} \log N = o(1)$, $nh^{2m+1} \geq 1$, and $f_0 = \sum_{\nu=1}^{\infty} f_\nu^0 \varphi_\nu$ satisfies Condition (S), then there exists a universal constant $M > 0$ s.t. $P(\|f - f_0\| \geq Mr_n | \boldsymbol{D}_j) = o_{P_{f_0}^n}(1)$ as $n \to \infty$.*

*Proof of Proposition A.4.* Note that there exists a universal constant $c' > 0$ such that $\Psi(x) \leq c' x^{1-1/(2m)}$ for any $0 < x < 1$. Therefore, there exists a universal constant $c'' > 0$ s.t. $B(h) \leq c'' h^{-(2m-1)/(4m)}$.

To prove the theorem, we first show the following posterior consistency: for any $\varepsilon > 0$, as $n \to \infty$,

$$P(\|f - f_0\|_\infty \geq \varepsilon | \mathbf{D}_n) \to 0, \text{ in } P_{f_0}^n\text{-probability.} \tag{A.61}$$

We can rewrite the posterior density of $f$ by

$$p(f | \mathbf{D}_n) = \frac{\prod_{i=1}^n (p_f/p_{f_0})(Z_i) \exp(-n\lambda J(f)/2) d\Pi(f)}{\int_{S^m(\mathbb{I})} \prod_{i=1}^n (p_f/p_{f_0})(Z_i) \exp(-n\lambda J(f)/2) d\Pi(f)},$$

where recall that $p_f(z)$ is the probability density of $Z = (Y, X)$ under $f$.

First of all, we give a lower bound for

$$I_1 = \int_{S^m(\mathbb{I})} \prod_{i=1}^n (p_f/p_{f_0})(Z_i) \exp(-n\lambda J(f)/2) d\Pi(f).$$

Define $B_n = \{f \in S^m(\mathbb{I}) : \|f - f_0\| \leq r_n\}$. Then

$$\begin{aligned}
I_1 &\geq \int_{B_n} \prod_{i=1}^n (p_f/p_{f_0})(Z_i) \exp(-n\lambda J(f)/2) d\Pi(f) \\
&= \int_{B_n} \exp(\sum_{i=1}^n R_i(f, f_0)) \exp(-n\lambda J(f)/2) d\Pi(f),
\end{aligned}$$

where $R_i(f, f_0) = \log(p_f(Z_i)/p_{f_0}(Z_i)) = Y_i(f(X_i) - f_0(X_i)) - A(f(X_i)) + A(f_0(X_i))$ for any $1 \leq i \leq n$. Define $d\Pi^*(f) = d\Pi(f)/\Pi(B_n)$, a reduced probability measure on $B_n$. By Jensen's inequality,

$$\begin{aligned}
&\log \int_{B_n} \exp(\sum_{i=1}^n R_i(f, f_0)) \exp(-n\lambda J(f)/2) d\Pi^*(f) \\
&\geq \int_{B_n} \left( \sum_{i=1}^n R_i(f, f_0) - n\lambda J(f)/2 \right) d\Pi^*(f) \\
&= \int_{B_n} \sum_{i=1}^n [R_i(f, f_0) - E_{f_0}\{R_i(f, f_0)\}] d\Pi^*(f) \\
&\quad + n \int_{B_n} E_{f_0}\{R_i(f, f_0)\} d\Pi^*(f) - \int_{B_n} \frac{n\lambda J(f)}{2} d\Pi^*(f) \\
&\equiv J_1 + J_2 + J_3.
\end{aligned}$$

For any $f \in B_n$, $\|f - f_0\| \le r_n$. By Lemma A.9 and the condition $h^{-3/2} r_n = o(1)$, we can choose $n$ to be sufficiently large so that $\|f - f_0\|_\infty \le ch^{-1/2}\|f - f_0\| \le ch^{-1/2} r_n \le 1$.

It follows from Assumption A1 that for $C = 1 + C_3\sqrt{J(f_0)}$, there exist positives $C_0', C_1', C_2'$ s.t. (2.2) and (2.3) hold with $C_0, C_1, C_2$ therein replaced by $C_0', C_1', C_2'$, respectively.

It follows by Taylor's expansion, $E_{f_0}\{Y_i - \dot{A}(f_0(X_i))|X_i\} = 0$, $1/C_2' \le \ddot{A}(z) \le C_2'$ for $|z| \le 2C$ and Assumption A1 that for any $f \in B_n$,

$$|E_{f_0}\{R_i(f, f_0)\}| \le C_2' E_{f_0}\{(f(X) - f_0(X))^2\} \le (C_2')^2 V(f - f_0) \le (C_2')^2 r_n^2.$$

Therefore, $J_2 \ge -(C_2')^2 n r_n^2$.

Since $r_n^2 = o(1)$, we can choose $n$ to be large so that $|E_{f_0}\{R_i(f, f_0)\}| \le 1$. Meanwhile, for any $f \in B_n$, for some $s \in [0, 1]$, we have

$$
\begin{aligned}
& |R_i(f, f_0)| \\
= & |Y_i(f(X_i) - f_0(X_i)) - A(f(X_i)) + A(f_0(X_i))| \\
= & |Y_i - \dot{A}(f_0(X_i)) \\
& - \frac{1}{2}\ddot{A}(f_0(X_i) + s(f(X_i) - f_0(X_i)))(f - f_0)(X_i)| \times |(f - f_0)(X_i)| \\
\le & |Y_i - \dot{A}(f_0(X_i))| + C_2'/2.
\end{aligned}
$$

We have used $\|f - f_0\|_\infty \le 1$ in the above inequalities.

For any $1 \le i \le n$, define $A_i = \{|Y_i - \dot{A}(f_0(X_i))| \le 2C_0' \log n\}$. It follows by Assumption A1 that $P_{f_0}^n(\cup_{i=1}^n A_i^c) \le C_1'/n \to 0$, as $n \to \infty$. Define $\xi_i = \int_{B_n} R_i(f, f_0)d\Pi^*(f) \times I_{A_i}$, we get that $|\xi_i| \le 2C_0' \log n + C_2'/2$, a.s. It can also be shown by $r_n^2 \ge 1/n$ that

$$
\begin{aligned}
& |E_{f_0}\{\int_{B_n} R_i(f, f_0)d\Pi^*(f) \times I_{A_i^c}\}| \\
\le & E_{f_0}\{(|Y_i - f_0(X_i)| + C_2'/2) \times I_{A_i^c}\} \\
= & E_{f_0}\{|Y_i - f_0(X_i)| \times I_{A_i^c}\} + \frac{C_2'}{2}P_{f_0}^n(A_i^c) \\
\le & C_0'\sqrt{2C_1' P_{f_0}^n(A_i^c)^{1/2}} + \frac{C_2'}{2}P_{f_0}^n(A_i^c) \\
\le & \frac{\sqrt{2}C_0'C_1'}{n} + \frac{C_1'C_2'}{2n^2} \le (\sqrt{2}C_0'C_1' + C_1'C_2')r_n^2.
\end{aligned}
$$

Let $\delta = 1/(\sqrt{n}r_n)$. Note that by the condition $h^{1/2}\log n = o(1)$ we have $\delta \log n = (\log n)/(\sqrt{n}r_n) \le h^{1/2}\log n = o(1)$, we can let $n$ be large so that $\delta(4C_0'\log n + C_2') \le 1$. Let $d_i = \xi_i - E_{f_0}\{\xi_i\}$ for $1 \le i \le n$, then it is easy to see that

$$|d_i| \le |\xi_i| + |E_{f_0}\{\xi_i\}| \le 4C_0'\log n + C_2', \ a.s.$$

Let $e_i = E_{f_0}\{\exp(\delta|d_i|) - 1 - \delta|d_i|\}$. It can be shown using inequality $\exp(x) - 1 - x \le x^2\exp(x)$

for $x \geq 0$ and Cauchy-Schwartz inequality that

$$
\begin{aligned}
|e_i| &\leq E_{f_0}\{\delta^2 d_i^2 \exp(\delta |d_i|)\} \\
&\leq e\delta^2 E_{f_0}\{d_i^2\} \\
&\leq e\delta^2 E_{f_0}\{\xi_i^2\} \\
&\leq e\delta^2 \int_{B_n} E_{f_0}\{R_i(f, f_0)^2\} d\Pi^*(f) \\
&\leq e\delta^2 \int_{B_n} E_{f_0}\{(|Y_i - \dot{A}(f_0(X_i))| + C_2'/2)^2 (f - f_0)(X_i)^2\} d\Pi^*(f) \\
&\leq e(4(C_0')^2 C_1' C_2' + (C_2')^3)\delta^2 r_n^2,
\end{aligned}
$$

where the last step follows from $V(f - f_0) \leq r_n^2$ for any $f \in B_n$. Therefore, it follows by [38, Theorem 3.2] that

$$
\begin{aligned}
& P_{f_0}^n \left( |\sum_{i=1}^n [\xi_i - E_{f_0}\{\xi_i\}]| \geq (e(4(C_0')^2 C_1' C_2' + (C_2')^3) + 2)\sqrt{n} r_n \log n \right) \\
&\leq 2 \exp(-(e(4(C_0')^2 C_1' C_2' + (C_2')^3) + 2)\sqrt{n} r_n (\log n)\delta \\
& \quad + e(4(C_0')^2 C_1' C_2' + (C_2')^3)\delta^2 n r_n^2) \\
&\leq 2/n^2 \to 0, \text{ as } n \to \infty.
\end{aligned}
\tag{A.62}
$$

Since $\sqrt{n} r_n \gg \log n$, we can let $n$ be large so that $(e(4(C_0')^2 C_1' C_2' + (C_2')^3) + 2)\sqrt{n} r_n \log n \leq n r_n^2$. Since on $\cap_{i=1}^n A_i$,

$$
J_1 = \sum_{i=1}^n [\xi_i - E_{f_0}\{\xi_i\}] - n E_{f_0}\{\int_{B_n} R_i(f, f_0) d\Pi^*(f) \times I_{A_i^c}\},
$$

we get from (A.62) that with $P_{f_0}^n$-probability approaching one,

$$
J_1 \geq -(e(4(C_0')^2 C_1' C_2' + (C_2')^3) + 2)\sqrt{n} r_n \log n - n r_n^2 \geq -2n r_n^2.
$$

Meanwhile, for any $f \in B_n$, $\lambda J(f - f_0) \leq r_n^2$. Therefore, $J_3 \geq -\frac{(1 + J(f_0)^{1/2})^2}{2} n r_n^2$. So, with probability approaching one,

$$
I_1 \geq \exp\left(-(2 + (C_2')^2) n r_n^2 - \frac{(1 + J(f_0)^{1/2})^2}{2} n r_n^2\right) \Pi(B_n).
$$

By Assumption A2,

$$
\Pi(B_n) \geq c_1 \exp(-c_0 r_n^{-2/(2m+\psi)}).
\tag{A.63}
$$

Since $\psi > 0$ and $r_n^2 = (nh)^{-1} + \lambda \geq n^{-2m/(2m+1)}$, we get $r_n^2 \geq \lambda$ and $n r_n^{\frac{2(2m+\beta)}{2m+\psi}} \geq n^{1 - \frac{2m(2m+\psi+1)}{(2m+1)(2m+\psi)}} > 1$, so $n r_n^2 > r_n^{-\frac{2}{2m+\psi}}$. Consequently, with $P_{f_0}^n$-probability approaching one,

$$
I_1 \geq c_1 \exp(-c_2 n r_n^2),
\tag{A.64}
$$

where $c_2 = 2 + (C_2')^2 + (1 + J(f_0)^{1/2})^2/2 + c_0$.

Now we choose a different constant $C$:

$$C = \max\{2C_3\sqrt{c_2+1}, c_2+1, 2(1+C_3\sqrt{J(f_0)})\}. \tag{A.65}$$

It follows by Assumption A1 that there exist positives $C_0, C_1, C_2$ s.t. (2.2) and (2.3) hold. Next we examine

$$I_2 := \int_{A_n} \prod_{i=1}^{n} (p_f/p_{f_0})(Z_i) \exp(-\frac{n\lambda}{2} J(f)) d\Pi(f),$$

where $A_n = \{f \in S^m(\mathbb{I}) : \|f - f_0\| \geq 3C_2\delta_n\}$, $\delta_n = 2bh^m + 24C_0c_K(C)(4C_1+C)r$, $r = r_n h^{-1/2}$, and $b = \frac{C_2C}{C_3}\sqrt{1 + \frac{1}{\rho_{m+1}}}$. By the condition $h^{-3/2}r_n = o(1)$ and $B(h) \lesssim h^{-(2m-1)/(4m)}$ it can be easily checked that the Rate Condition (**H**): (i)–(iv) are satisfied (when $n$ becomes large) with $M$ therein replaced by $C$. Define test $\phi_n = I(\|\widehat{f}_{n,\lambda} - f_0\| \geq C_2\delta_n)$. Since $C_2 \geq 1$, it follows by part (a) of Lemma A.12 that

$$E_{f_0}\{\phi_n\} = P_{f_0}^n(\|\widehat{f}_{n,\lambda} - f_0\| \geq C_2\delta_n) \leq P_{f_0}^n(\|\widehat{f}_{n,\lambda} - f_0\| \geq \delta_n) \leq 6\exp(-Cnr_n^2),$$

and by (A.48),

$$
\begin{aligned}
\sup_{\substack{f \in H^m(C) \\ \|f-f_0\| \geq 3C_2\delta_n}} E_f\{1 - \phi_n\} &= \sup_{\substack{f \in H^m(C) \\ \|f-f_0\| \geq 3C_2\delta_n}} P_f^n(\|\widehat{f}_{n,\lambda} - f_0\| < C_2\delta_n) \\
&\leq \sup_{\substack{f \in H^m(C) \\ \|f-f_0\| \geq 3C_2\delta_n}} P_f^n(\|\widehat{f}_{n,\lambda} - f\| \geq 2C_2\delta_n) \\
&\leq \sup_{\substack{f \in H^m(C) \\ \|f-f_0\| \geq 3C_2\delta_n}} P_f^n(\|\widehat{f}_{n,\lambda} - f\|_f \geq \delta_n) \\
&\leq 6\exp(-Cnr_n^2),
\end{aligned}
$$

where the second last inequality follows by Lemma A.7.

Note that for any $f \in A_n \backslash H^m(C)$,

$$J(f) > (1 + 1/\rho_{m+1})^{-1} C_2^{-2} b^2 = C^2/C_3^2 \geq 4(c_2+1).$$

Since $nh^{2m+1} \geq 1$ leads to $r_n^2 = (nh)^{-1} + \lambda \leq 2\lambda$, it then holds that,

$$
\begin{aligned}
&E_{f_0}\{I_2(1 - \phi_n)\} \\
&= \int_{A_n} E_f\{1 - \phi_n\} \exp(-n\lambda J(f)/2) d\Pi(f) \\
&= \int_{A_n \backslash H^m(C)} E_f\{1 - \phi_n\} \exp(-n\lambda J(f)/2) d\Pi(f) \\
&\quad + \int_{A_n \cap H^m(C)} E_f\{1 - \phi_n\} \exp(-n\lambda J(f)/2) d\Pi(f) \\
&\leq \exp(-2n\lambda(c_2+1)) + 6\exp(-(c_2+1)nr_n^2) \\
&\leq \exp(-(c_2+1)nr_n^2) + 6\exp(-(c_2+1)nr_n^2) = 7\exp(-(c_2+1)nr_n^2),
\end{aligned}
$$

so

$$E_{f_0}\{I_2(1-\phi_n)\} \le 7\exp(-(c_2+1)nr_n^2),$$

which implies $I_2(1-\phi_n) = O_{P_{f_0}^n}(\exp(-(c_2+1)nr_n^2))$. On the other hand,

$$E_{f_0}\{P(A_n|\mathbf{D}_n)\phi_n\} \le P_{f_0}^n(\|\widehat{f}_{n,\lambda}-f_0\| \ge C_2\delta_n) \le 6\exp(-(c_2+1)nr_n^2),$$

so as $n \to \infty$,

$$E_{f_0}\{P(A_n|\mathbf{D}_n)\phi_n\} \le 6\exp(-(c_2+1)nr_n^2) \to 0,$$

which implies that $P(A_n|\mathbf{D}_n)\phi_n = o_{P_{f_0}^n}(1)$. By the above arguments and (A.64), we have

$$
\begin{aligned}
& P(A_n|\mathbf{D}_n) \\
={} & P(A_n|\mathbf{D}_n)\phi_n + P(A_n|\mathbf{D}_n)(1-\phi_n) \\
\le{} & P(A_n|\mathbf{D}_n)\phi_n + \frac{I_2(1-\phi_n)}{I_1} \\
={} & o_{P_{f_0}^n}(1) + O_{P_{f_0}^n}(\exp(-(c_2+1)nr_n^2)\exp(c_2nr_n^2)) = o_{P_{f_0}^n}(1),
\end{aligned}
$$

where the last step follows from $\exp(-nr_n^2) \le \exp(-h^{-1}) = o(1)$. By condition $r_n h^{-3/2} = o(1)$ and the trivial fact $\delta_n \asymp r_n h^{-1/2}$, we have that $h^{-1/2}\delta_n = o(1)$, together with Lemma A.9 we have that (A.61) holds.

To prove the theorem, we let

$$I_2' := \int_{A_n'} \prod_{i=1}^n (p_f/p_{f_0})(Z_i)\exp(-\frac{n\lambda}{2}J(f))d\Pi(f),$$

where $A_n' = \{f \in S^m(\mathbb{I}) : \|f-f_0\| \ge \sqrt{2}Mr_n\}$ for a fixed number

$$M > \max\{2, J(f_0)^{1/2} + \sqrt{2(c_2+1)}, 1 + \|f_0\|_\infty\}$$

to be further described later. Let

$$A_{n1}' = \{f \in S^m(\mathbb{I}) : V(f-f_0) \ge M^2 r_n^2, \lambda J(f-f_0) \le M^2 r_n^2\}$$

and

$$A_{n2}' = \{f \in S^m(\mathbb{I}) : \lambda J(f-f_0) \ge M^2 r_n^2\}.$$

For any $f \in A_{n2}'$, it can be shown that

$$Mr_n \le \sqrt{\lambda J(f-f_0)} \le \sqrt{\lambda}(J(f)^{1/2} + J(f_0)^{1/2}) \le (\lambda J(f))^{1/2} + J(f_0)^{1/2}r_n,$$

which leads to $\lambda J(f) \ge (M - J(f_0)^{1/2})^2 r_n^2$. So we have

$$
\begin{aligned}
& E_{f_0}\{\int_{A_{n2}'} \prod_{i=1}^n (p_f/p_{f_0})(Z_i)\exp(-\frac{n\lambda}{2}J(f))d\Pi(f)\} \\
={} & \int_{A_{n2}'} \exp(-\frac{n\lambda}{2}J(f))d\Pi(f) \le \exp(-(M-J(f_0)^{1/2})^2 nr_n^2/2),
\end{aligned}
$$

which leads to that

$$
\begin{aligned}
&\int_{A'_{n2}} \prod_{i=1}^n (p_f/p_{f_0})(Z_i) \exp(-\frac{n\lambda}{2}J(f)) d\Pi(f) \\
&= O_{P^n_{f_0}}(\exp(-(M - J(f_0)^{1/2})^2 n r_n^2/2)).
\end{aligned} \tag{A.66}
$$

To continue, we need to build uniformly consistent test. Let $d_H^2(P_f, P_g) = \frac{1}{2}\int(\sqrt{dP_f} - \sqrt{dP_g})^2$ be the squared Hellinger distance between the two probability measures $P_f(z)$ and $P_g(z)$. Recall that their corresponding probability density functions are $p_f$ and $p_g$, respectively. Next we present a lemma showing the local equivalence of $V$ and $d_H^2$.

**Lemma A.13.** *Let $C$ be chosen as (A.65) and $C_0, C_1, C_2$ be positives satisfying Assumption A1. Let $\varepsilon > 0$ satisfy $\varepsilon < \min\{1, 1/C_0, C\}$ and*

$$
\frac{1}{12}C_2^2\varepsilon + \frac{1}{32}C_2^3\varepsilon^2 + C_0^3 C_1 C_2 \varepsilon \exp\left(\frac{\varepsilon}{4}C_2 + \frac{C_2}{4C_0}\right) < \frac{1}{16}.
$$

*Then for any $f, g \in \mathcal{F}(C)$ satisfying $\|f - g\|_\infty \le \varepsilon$,*

$$
V(f - g)/16 \le d_H^2(P_f, P_g) \le 3V(f - g)/16.
$$

*Proof of Lemma A.13.* For any $f, g \in \mathcal{F}(C)$ with $\|f - g\|_\infty \le \varepsilon$, define $\Delta_Z(f, g) = \frac{1}{2}[Y(f(X) - g(X)) - A(f(X)) + A(g(X))]$, where recall and $Z = (Y, X)$. It is easy to see by direct calculations that

$$
d_H^2(P_f, P_g) = 1 - E_g\{\exp(\Delta_Z(f, g))\}.
$$

By Taylor's expansion, for some random $t \in [0, 1]$,

$$
\begin{aligned}
&1 - E_g\{\exp(\Delta_Z(f, g))\} \\
&= -E_g\{\Delta_Z(f, g)\} - \frac{1}{2}E_g\{\Delta_Z(f, g)^2\} - \frac{1}{6}E_g\{\exp(t\Delta_Z(f, g))\Delta_Z(f, g)^3\}.
\end{aligned}
$$

We will analyze the terms on the right side of the equation.

Define $\xi = Y - \dot{A}(g(X))$. By [28] we get $E_g\{\xi|X\} = 0$ and $E_g\{\xi^2|X\} = \ddot{A}(g(X))$. By Taylor's expansion,

$$
\begin{aligned}
\Delta_Z(f, g) &= \frac{1}{2}[\xi(f(X) - g(X)) - \frac{1}{2}\ddot{A}(g(X))(f(X) - g(X))^2 \\
&\quad - \frac{1}{6}\dddot{A}(f_{1*}(X))(f(X) - g(X))^3],
\end{aligned}
$$

$$
\Delta_Z(f, g) = \frac{1}{2}[\xi(f(X) - g(X)) - \frac{1}{2}\ddot{A}(f_{2*}(X))(f(X) - g(X))^2],
$$

where $f_{k*}(X)$ is between $g(X)$ and $f(X)$ for $k = 1, 2$. It clearly holds that $\|f_{k*}\|_\infty \le \|f\|_\infty + \|g - f\|_\infty < 2C$. Then we get that

$$
-E_g\{\Delta_Z(f, g)\} = \frac{1}{4}V(f - g) + \frac{1}{12}E_g\{\dddot{A}(f_{1*}(X))(f(X) - g(X))^3\},
$$

and

$$E_g\{\Delta_Z(f,g)^2\}$$
$$= E_g\{(\frac{1}{2}\xi(f(X)-g(X))-\frac{1}{4}\ddot{A}(f_{2*}(X))(f(X)-g(X))^2)\}$$
$$= \frac{1}{4}E_g\{\xi^2(f(X)-g(X))^2\}-\frac{1}{4}E_g\{\xi(f(X)-g(X))^3\ddot{A}(f_{2*}(X))\}$$
$$+\frac{1}{16}E_g\{\ddot{A}(f_{2*}(X))^2(f(X)-g(X))^4\}$$
$$= \frac{1}{4}V(f-g)+\frac{1}{16}E_g\{\ddot{A}(f_{2*}(X))^2(f(X)-g(X))^4\}.$$

Since $\|f-g\|_\infty \leq \varepsilon < \min\{1, 1/C_0, C\}$ and $0 < \ddot{A}(z) \leq C_2$ for any $z \in [-2C, 2C]$, implying $|\Delta_Z(f,g)| \leq \frac{1}{2}(|\xi|+C_2/2)|f(X)-g(X)|$, we get

$$|E_g\{\exp(t\Delta_Z(f,g))\Delta_Z(f,g)^3\}|$$
$$\leq E_g\{\exp(|\Delta_Z(f,g)|)|\Delta_Z(f,g)|^3\}$$
$$\leq E_g\{\exp(\varepsilon|\xi|/2+C_2\varepsilon/4)(|\xi|/2+C_2/4)^3|f(X)-g(X)|^3\}$$
$$= 6C_0^3E_g\left\{\exp(\varepsilon|\xi|/2+C_2\varepsilon/4)\times\frac{1}{3!}\left(\frac{|\xi|/2+C_2/4}{C_0}\right)^3|f(X)-g(X)|^3\right\}$$
$$\leq 6C_0^3E_g\{\exp(\varepsilon|\xi|/2+C_2\varepsilon/4)\exp(|\xi|/(2C_0)+C_2/(4C_0))|f(X)-g(X)|^3\}$$
$$\leq 6C_0^3\exp(C_2\varepsilon/4+C_2/(4C_0))E_g\{\exp(|\xi|/C_0)|f(X)-g(X)|^3\}$$
$$\leq 6C_0^3C_1C_2\exp(C_2\varepsilon/4+C_2/(4C_0))\varepsilon V(f-g).$$

It also holds that

$$|E_g\{\ddot{A}(f_{1*}(X))(f(X)-g(X))^3\}| \leq C_2^2\varepsilon V(f-g),$$
$$|E_g\{\ddot{A}(f_{2*}(X))^2(f(X)-g(X))^4\}| \leq C_2^3\varepsilon^2 V(f-g).$$

Therefore, by the above argument it holds that, for any $f, g \in \mathcal{F}(C)$ with $\|f-g\|_\infty \leq \varepsilon$,

$$|d_H^2(P_f, P_g)-V(f-g)/8|$$
$$= |\frac{1}{12}E_g\{\ddot{A}(f_{1*}(X))(f(X)-g(X))^3\}$$
$$-\frac{1}{32}E_g\{\ddot{A}(f_{2*}(X))^2(f(X)-g(X))^4\}$$
$$-\frac{1}{6}E_g\{\exp(t\Delta_Z(f,g))\Delta_Z(f,g)^3\}|$$
$$\leq \left(\frac{1}{12}C_2^2\varepsilon+\frac{1}{32}C_2^3\varepsilon^2+C_0^3C_1C_2\exp(C_2\varepsilon/4+C_2/(4C_0))\varepsilon\right)V(f-g)$$
$$< V(f-g)/16,$$

which implies $V(f-g)/16 \leq d_H^2(P_f, P_g) \leq 3V(f-g)/16$. This proves Lemma A.13. $\qquad\square$

Let $\varepsilon$ satisfy the conditions in Lemma A.13. Define $\mathcal{F}_n = \{f \in S^m(\mathbb{I}) : \|f - f_0\|_\infty \leq \varepsilon/2, J(f) \leq (M + J(f_0)^{1/2})^2 r_n^2 \lambda^{-1}\}$. Note that for any $f \in \mathcal{F}_n$, we have $\|f\|_\infty \leq \|f_0\|_\infty + \varepsilon/2 < C$. Therefore, $\mathcal{F}_n \subseteq \mathcal{F}(C)$. Let $\mathcal{P}_n = \{P_f^n : f \in \mathcal{F}_n\}$ and $D(\delta, \mathcal{P}_n, d_H)$ be the $\delta$-packing number in terms of $d_H$. Since $r_n^2 \geq \lambda$ which leads to $(M + J(f_0)^{1/2}) r_n h^{-m} > M + J(f_0)^{1/2} > \varepsilon + \|f_0\|_\infty$, it can be easily checked that $\mathcal{F}_n \subset (M + J(f_0)^{1/2}) r_n h^{-m} \mathcal{T}$, where $\mathcal{T} = \{f \in S^m(\mathbb{I}) : \|f\|_\infty \leq 1, J(f) \leq 1\}$.

For any $f, g \in \mathcal{F}_n$ (implying $f, g \in \mathcal{F}(C)$) with $\|f - g\|_\infty \leq \varepsilon$, it follows by Lemma A.13 that $D(\delta, \mathcal{P}_n, d_H) \leq D(4\delta/\sqrt{3}, \mathcal{F}_n, d_V)$, where $d_V$ is the distance induced by $V$, i.e., $d(f, g) = V^{1/2}(f - g)$. And hence, it follows by [24, Theorem 9.21] that

$$
\begin{aligned}
\log D(\delta, \mathcal{P}_n, d_H) &\leq \log D(4\delta/\sqrt{3}, \mathcal{F}_n, d_V) \\
&\leq \log D(4\delta/\sqrt{3}, (M + J(f_0)^{1/2}) r_n h^{-m} \mathcal{T}, d_V) \\
&\leq c_V \left( \frac{\delta}{(M + J(f_0)^{1/2}) r_n h^{-m}} \right)^{-1/m},
\end{aligned}
$$

where $c_V$ is a universal constant only depending on the regularity level $m$. This implies that for any $\delta > 2r_n$,

$$
\begin{aligned}
\log D(\delta/2, \mathcal{P}_n, d_H) &\leq \log D(r_n, \mathcal{P}_n, d_H) \\
&\leq c_V (M + J(f_0)^{1/2})^{1/m} h^{-1} \\
&\leq c_V (M + J(f_0)^{1/2})^{1/m} n r_n^2,
\end{aligned}
$$

where the last inequality follows by the fact $r_n^2 \geq (nh)^{-1}$. Thus, the right side of the above inequality is constant in $\delta$. By [18, Theorem 7.1], with $\delta = M r_n/4$, there exists test $\tilde{\phi}_n$ and a universal constant $k_0 > 0$ satisfying

$$
\begin{aligned}
E_{f_0}\{\tilde{\phi}_n\} &= P_{f_0}^n \tilde{\phi}_n \\
&\leq \frac{\exp(c_V (M + J(f_0)^{1/2})^{1/m} n r_n^2) \exp(-k_0 n \delta^2)}{1 - \exp(-k_0 n \delta^2)} \\
&= \frac{\exp(c_V (M + J(f_0)^{1/2}) n r_n^2 - k_0 M^2 n r_n^2/16)}{1 - \exp(-k_0 M^2 n r_n^2/16)},
\end{aligned}
$$

and, combined with Lemma A.13,

$$
\begin{aligned}
\sup_{\substack{f \in \mathcal{F}_n \\ d_V(f, f_0) \geq 4\delta}} E_f\{1 - \tilde{\phi}_n\} &= \sup_{\substack{f \in \mathcal{F}_n \\ d_V(f, f_0) \geq 4\delta}} P_f^n \{1 - \tilde{\phi}_n\} \\
&\leq \sup_{\substack{f \in \mathcal{F}_n \\ d_H(P_f^n, P_{f_0}^n) \geq \delta}} P_f^n \{1 - \tilde{\phi}_n\} \\
&\leq \exp(-k_0 n \delta^2) = \exp(-k_0 M^2 n r_n^2/16).
\end{aligned}
$$

This implies that

$$E_{f_0}\Big\{\int_{\substack{f\in\mathcal{F}_n \\ d_V(f,f_0)\geq 4\delta}} \prod_{i=1}^n (p_f/p_{f_0})(Z_i)\exp(-n\lambda J(f)/2)d\Pi(f)(1-\tilde{\phi}_n)\Big\}$$

$$\leq \int_{\substack{f\in\mathcal{F}_n \\ d_V(f,f_0)\geq 4\delta}} E_{f_0}\Big\{\prod_{i=1}^n (p_f/p_{f_0})(Z_i)(1-\tilde{\phi}_n)\Big\}d\Pi(f)$$

$$= \int_{\substack{f\in\mathcal{F}_n \\ d_V(f,f_0)\geq 4\delta}} E_f\{1-\tilde{\phi}_n\}d\Pi(f)$$

$$\leq \exp(-k_0 M^2 n r_n^2/16).$$

Therefore,

$$\int_{\substack{f\in\mathcal{F}_n \\ d_V(f,f_0)\geq 4\delta}} \prod_{i=1}^n (p_f/p_{f_0})(Z_i)\exp(-n\lambda J(f)/2)d\Pi(f)(1-\tilde{\phi}_n)$$

$$= O_{P_{f_0}^n}\left(\exp(-k_0 M^2 n r_n^2/16)\right). \tag{A.67}$$

It follows from (A.64) and (A.66) that

$$P(A'_{n2}|\mathbf{D}_n) = O_{P_{f_0}^n}\left(\exp(-(M-J(f_0)^{1/2})^2 n r_n^2/2 + c_2 n r_n^2)\right) = o_{P_{f_0}^n}(1),$$

where the last inequality follows by $(M-J(f_0)^{1/2})^2 > 2(c_2+1)$ and $\exp(-n r_n^2) = o(1)$. Together with (A.61), we get that

$$P(A'_n|\mathbf{D}_n)$$

$$\leq P(A'_{n1}|\mathbf{D}_n) + P(A'_{n2}|\mathbf{D}_n)$$

$$\leq P(A'_{n1}, \|f-f_0\|_\infty \leq \varepsilon/2|\mathbf{D}_n) + P(\|f-f_0\|_\infty > \varepsilon/2|\mathbf{D}_n) + P(A'_{n2}|\mathbf{D}_n)$$

$$\leq P(A'_{n1}, \|f-f_0\|_\infty \leq \varepsilon/2|\mathbf{D}_n) + o_{P_{f_0}^n}(1)$$

$$\leq P(A'_{n1}, \|f-f_0\|_\infty \leq \varepsilon/2|\mathbf{D}_n)\tilde{\phi}_n$$

$$\quad + P(A'_{n1}, \|f-f_0\|_\infty \leq \varepsilon/2|\mathbf{D}_n)(1-\tilde{\phi}_n) + o_{P_{f_0}^n}(1).$$

Choose the constant $M$ to be even bigger so that

$$c_V(M+J(f_0)^{1/2})+1 < k_0 M^2/16, \quad 1+c_2 \leq k_0 M^2/16.$$

Then we get that

$$E_{f_0}\{P(A'_{n1}, \|f-f_0\|_\infty \leq \varepsilon/2|\mathbf{D}_n)\tilde{\phi}_n\}$$

$$\lesssim \exp(c_V(M+J(f_0)^{1/2})n r_n^2 - k_0 M^2 n r_n^2/16)$$

$$\leq \exp(-n r_n^2) = o(1),$$

leading to $P(A'_{n1}, \|f - f_0\|_\infty \le \varepsilon/2 | \mathbf{D}_n) \tilde{\phi}_n = o_{P^n_{f_0}}(1)$. Meanwhile, it follows by (A.64) and (A.67) that

$$
\begin{aligned}
& P(A'_{n1}, \|f - f_0\|_\infty \le \varepsilon/2 | \mathbf{D}_n)(1 - \tilde{\phi}_n) \\
\le\ & P(f \in \mathcal{F}_n, d_V(f, f_0) \ge 4\delta | \mathbf{D}_n)(1 - \tilde{\phi}_n) \\
\le\ & \frac{\int_{\substack{f \in \mathcal{F}_n \\ d_V(f, f_0) \ge 4\delta}} \prod_{i=1}^n (p_f/p_{f_0})(Z_i) \exp(-n\lambda J(f)/2) d\Pi(f)(1 - \tilde{\phi}_n)}{I_1} \\
=\ & O_{P^n_{f_0}}\left(\exp(-k_0 M^2 n r_n^2/16 + c_2 n r_n^2)\right) \\
=\ & O_{P^n_{f_0}}\left(\exp(-n r_n^2)\right) = o_{P^n_{f_0}}(1).
\end{aligned}
$$

Thus, we have shown that $P(\|f - f_0\| \ge \sqrt{2} M r_n | \mathbf{D}_n) = o_{P^n_{f_0}}(1)$. This completes the proof. □

*Proof of Proposition A.1.* Fix any $\varepsilon_1, \varepsilon_2 \in (0, 1)$. Let $C = C_3 \sqrt{J(f_0)} + 1$, and $C_0, C_1, C_2$ be positive constants satisfying (2.2) and (2.3) in Assumptions A1. It follows by Lemma A.12 that for any fixed constant $M > 1$, if we set

$$
b = \frac{C_2 C}{C_3} \sqrt{1 + \frac{1}{\rho_{m+1}}}, r = (nh/\log 2s)^{-1/2}, \delta_n = 2bh^m + 24 C_0 c_K (4C_1 + M) r, \tag{A.68}
$$

$$
a_n = C_2 c_K^2 M^{1/2} h^{-1/2} r B(h) \delta_n, \text{ and } b_n = C_2^2 c_K h^{-1/2} \delta_n^2, \tag{A.69}
$$

then as $n \to \infty$,

$$
P^n_{f_0}\left(\|\hat{f}_{n,\lambda} - f_0\| \ge \delta_n\right) \le 6n^{-M} \to 0,
$$

and

$$
P^n_{f_0}\left(\|\hat{f}_{n,\lambda} - f_0 - S_{n,\lambda}(f_0)\| > a_n + b_n\right) \le 8n^{-M} \to 0.
$$

By $B(h) \lesssim h^{-\frac{2m-1}{4m}}$ and the simple fact $a_n + b_n \lesssim D_n$, we get that

$$
\|\hat{f}_{n,\lambda} - f_0 - S_{n,\lambda}(f_0)\| = O_{P^n_{f_0}}(a_n + b_n) = O_{P^n_{f_0}}(D_n). \tag{A.70}
$$

Recall that

$$
S_{n,\lambda}(f_0) = \frac{1}{n} \sum_{i=1}^n (Y_i - \dot{A}(f_0(X_i))) K_{X_i} - \mathcal{P}_\lambda f_0.
$$

It was shown by [44] that $\mathcal{P}_\lambda \varphi_\nu = \frac{\lambda \rho_\nu}{1 + \lambda \rho_\nu} \varphi_\nu$. Since $f_0$ satisfies Condition (**S**),

$$
\begin{aligned}
\|\mathcal{P}_\lambda f_0\|^2 &= \langle \sum_{\nu=1}^\infty f_\nu^0 \frac{\lambda \rho_\nu}{1 + \lambda \rho_\nu} \varphi_\nu, \sum_{\nu=1}^\infty f_\nu^0 \frac{\lambda \rho_\nu}{1 + \lambda \rho_\nu} \varphi_\nu \rangle \\
&= \sum_{\nu=1}^\infty |f_\nu^0|^2 \frac{\lambda^2 \rho_\nu^2}{1 + \lambda \rho_\nu} \\
&= \lambda^{1 + \frac{\beta-1}{2m}} \sum_{\nu=1}^\infty |f_\nu^0|^2 \rho_\nu^{1 + \frac{\beta-1}{2m}} \frac{(\lambda \rho_\nu)^{1 - \frac{\beta-1}{2m}}}{1 + \lambda \rho_\nu} = O(h^{2m+\beta-1}),
\end{aligned}
$$

where the last equation follows by $\lambda = h^{2m}$, $\sup_{x \geq 0} \frac{x^{1-\frac{\beta-1}{2m}}}{1+x} < \infty$, and Condition (**S**). On the other side, it follows by the proof of (A.55) that

$$P_{f_0}^n \left( \| \sum_{i=1}^n (Y_i - \dot{A}(f_0(X_i)))K_{X_i} \| \geq L(M)n(nh/\log 2)^{-1/2} \right)$$
$$\leq 2\exp\left(-Mnh(nh/\log 2)^{-1}\right) = 2^{1-M} \to 0, \text{ as } M \to \infty,$$

implying that

$$\| \sum_{i=1}^n (Y_i - \dot{A}(f_0(X_i)))K_{X_i} \| = O_{P_{f_0}^n}\left(n(nh/\log 2)^{-1/2}\right),$$

and hence,

$$\|S_{n,\lambda}(f_0)\| = O_{P_{f_0}^n}\left((nh)^{-1/2} + h^{m+\frac{\beta-1}{2}}\right) = O_{P_{f_0}^n}(\widetilde{r}_n).$$

Together with (A.70) and the rate condition $D_n \lesssim \widetilde{r}_n$, we get that

$$\|\widehat{f}_{n,\lambda} - f_0\| = O_{P_{f_0}^n}(\widetilde{r}_n).$$

Let $M_1$ be large constant so that event

$$\mathcal{E}'_n = \{\|\widehat{f}_{n,\lambda} - f_0\| \leq M_1\widetilde{r}_n\} \tag{A.71}$$

has probability approaching one. Meanwhile, form some positive constant $M_0$, it follows by Theorem A.4 that $P(\|f - f_0\| \geq M_0 r_n | \mathbf{D}_n)$ converges to zero in $P_{f_0}^n$-probability. Let $C' > M_1$ be a constant to be further determined later, then we have that

$$P(\|f - f_0\| \geq 2C'\widetilde{r}_n | \mathbf{D}_n)$$
$$\leq P(\|f - f_0\| \geq M_0 r_n | \mathbf{D}_n) + P(2C'\widetilde{r}_n \leq \|f - f_0\| \leq M_0 r_n | \mathbf{D}_n).$$

Thanks to Theorem A.4, the first term converges to zero in $P_{f_0}^n$-probability. Thus, when $n$ is sufficiently large,

$$P_{f_0}^n\left(P(\|f - f_0\| \geq M_0 r_n | \mathbf{D}_n) \geq \varepsilon_2/2\right) \leq \varepsilon_1/2.$$

We only need to handle the second term.

Define

$$\mathcal{E}''_n = \left\{\sup_{g \in \mathcal{G}} \|Z_{n,f_0}^{(l)}(g)\| \leq B(h)\sqrt{M\log n}, \; l = 1, 2\right\}, \tag{A.72}$$

where

$$Z_{n,f_0}^{(l)}(g) = \frac{1}{\sqrt{n}}\sum_{i=1}^n [\psi_{n,f_0}^{(l)}(Z_i; g)K_{X_i} - E_f\{\psi_{n,f_0}^{(l)}(Z_i; g)K_{X_i}\}] \text{ for } l = 1, 2,$$

and

$$\psi_{n,f_0}^{(1)}(Z_i; g) = c_K^{-1}h^{1/2}g(X_i),$$
$$\psi_{n,f_0}^{(2)}(Z_i; g) = C_2^{-1}c_K^{-1}h^{1/2}\ddot{A}(f_0(X_i))g(X_i).$$

It is easy to see that $\psi_{n,f_0}^{(l)}(Z_i; g)$ satisfies (A.51). By Lemma A.11 we have that $\mathcal{E}_n''$ has $P_{f_0}^n$-probability approaching one. Thus, it holds that, when $n$ becomes large, $P_{f_0}^n(\mathcal{E}_n) \geq 1 - \varepsilon_1/2$. In the rest of the proof we simply assume that $\mathcal{E}_n \equiv \mathcal{E}_n' \cap \mathcal{E}_n''$ holds.

Let $I_n(f) = \int_0^1 \int_0^1 s D S_{n,\lambda}(\widehat{f}_{n,\lambda} + ss'(f - \widehat{f}_{n,\lambda}))(f - \widehat{f}_{n,\lambda})(f - \widehat{f}_{n,\lambda}) ds ds'$ for any $f \in S^m(\mathbb{I})$. Let $\Delta f = f - \widehat{f}_{n,\lambda}$. Therefore,

$$
\begin{aligned}
& I_n(f) \\
=\ & -\frac{1}{n} \int_0^1 \int_0^1 s \sum_{i=1}^n \ddot{A}(\widehat{f}_{n,\lambda}(X_i) + ss'(\Delta f)(X_i))(\Delta f)(X_i)^2 ds ds' \\
& -\lambda J(\Delta f, \Delta f)/2 \\
=\ & -\frac{1}{n} \int_0^1 \int_0^1 s \sum_{i=1}^n [\ddot{A}(\widehat{f}_{n,\lambda}(X_i) + ss'(\Delta f)(X_i))(\Delta f)(X_i)^2 \\
& -\ddot{A}(f_0(X_i))(\Delta f)(X_i)^2] ds ds' \\
& -\frac{1}{2n} \sum_{i=1}^n [\ddot{A}(f_0(X_i))(\Delta f)(X_i)^2 - E_{f_0}^X\{\ddot{A}(f_0(X))(\Delta f)(X)^2\}] - \frac{1}{2}\|\Delta f\|^2 \\
\equiv\ & T_1(f) + T_2(f) - \frac{1}{2}\|\Delta f\|^2,
\end{aligned}
$$

where recall that

$$
\begin{aligned}
T_1(f) &= -\frac{1}{n} \int_0^1 \int_0^1 s \sum_{i=1}^n [\ddot{A}(\widehat{f}_{n,\lambda}(X_i) + ss'(\Delta f)(X_i))(\Delta f)(X_i)^2 \\
& \qquad -\ddot{A}(f_0(X_i))(\Delta f)(X_i)^2] ds ds', \\
T_2(f) &= -\frac{1}{2n} \sum_{i=1}^n [\ddot{A}(f_0(X_i))(\Delta f)(X_i)^2 - E_{f_0}^X\{\ddot{A}(f_0(X))(\Delta f)(X)^2\}].
\end{aligned}
$$

By Taylor's expansion in terms of Fréchet derivatives,

$$
\ell_{n,\lambda}(f) - \ell_{n,\lambda}(\widehat{f}_{n,\lambda}) = S_{n,\lambda}(\widehat{f}_{n,\lambda})(f - \widehat{f}_{n,\lambda}) + I_n(f) = I_n(f).
$$

Therefore,

$$
\begin{aligned}
P(A_n | \mathbf{D}_n) &= \frac{\int_{A_n} \exp(n(\ell_{n,\lambda}(f) - \ell_{n,\lambda}(\widehat{f}_{n,\lambda}))) d\Pi(f)}{\int_{S^m(\mathbb{I})} \exp(n(\ell_{n,\lambda}(f) - \ell_{n,\lambda}(\widehat{f}_{n,\lambda}))) d\Pi(f)} \\
&= \frac{\int_{A_n} \exp(n I_n(f)) d\Pi(f)}{\int_{S^m(\mathbb{I})} \exp(n I_n(f)) d\Pi(f)},
\end{aligned}
$$

where $A_n = \{f \in S^m(\mathbb{I}) : 2C'\widetilde{r}_n \leq \|f - f_0\| \leq M_0 r_n\}$.

Let

$$
J_1 = \int_{S^m(\mathbb{I})} \exp(n I_n(f)) d\Pi(f), \ J_2 = \int_{A_n} \exp(n I_n(f)) d\Pi(f).
$$

Then on $\mathcal{E}_n$ and for $\|f - f_0\| \leq \widetilde{r}_n$, we have $\|f - \widehat{f}_{n,\lambda}\| \leq \|f - f_0\| + \|\widehat{f}_{n,\lambda} - f_0\| \leq (M_1 + 1)\widetilde{r}_n$.

Let $d_n = c_K(M_1 + 1)h^{-1/2}\widetilde{r}_n$. It follows by similar arguments above (A.56) that $d_n^{-1}\Delta f \in \mathcal{G}$. It follows by Lemma A.9 that $\|\Delta f\|_\infty \leq c_K h^{-1/2}\|\Delta f\| \leq c_K(M_1 + 1)h^{-1/2}\widetilde{r}_n$. By rate assumption $r_n = o(h^{3/2})$ and $h^{1/2}\log n = o(1)$ and the simple fact $\widetilde{r}_n \leq r_n\sqrt{\log 2n}$, we get that

$$h^{-1/2}\widetilde{r}_n \leq h^{-1/2}r_n\sqrt{\log 2N} = o(h\sqrt{\log n}) = o(1).$$

Therefore, we can let $n$ be large so that, on $\mathcal{E}_n$ and $\|f_0\|_\infty + \|\widehat{f}_{n,\lambda} - f_0\|_\infty + \|\Delta f\|_\infty < C$. Then on $\mathcal{E}_n$, we have

$$
\begin{aligned}
&|T_1(f)| \\
\leq\ & C_2 \frac{\|\widehat{f}_{n,\lambda} - f_0\|_\infty + \|\Delta f\|_\infty}{2n} \sum_{i=1}^n (\Delta f)(X_i)^2 \\
=\ & C_2 \frac{\|\widehat{f}_{n,\lambda} - f_0\|_\infty + \|\Delta f\|_\infty}{2n} \sum_{i=1}^n [(\Delta f)(X_i)^2 - E^X\{(\Delta f)(X)^2\}] \\
& + C_2 \frac{\|\widehat{f}_{n,\lambda} - f_0\|_\infty + \|\Delta f\|_\infty}{2} E^X\{(\Delta f)(X)^2\} \\
\leq\ & C_2 \frac{\|\widehat{f}_{n,\lambda} - f_0\|_\infty + \|\Delta f\|_\infty}{2n} \|\Delta f\| \\
& \times \|\sum_{i=1}^n [(\Delta f)(X_i)K_{X_i} - E^X\{(\Delta f)(X)K_X\}]\| \\
& + C_2 \frac{\|\widehat{f}_{n,\lambda} - f_0\|_\infty + \|\Delta f\|_\infty}{2} E^X\{(\Delta f)(X)^2\} \\
\leq\ & C_2 d_n \frac{\|\widehat{f}_{n,\lambda} - f_0\|_\infty + \|\Delta f\|_\infty}{2n} \|\Delta f\| \\
& \times \|\sum_{i=1}^n [d_n^{-1}(\Delta f)(X_i)K_{X_i} - E^X\{d_n^{-1}(\Delta f)(X)K_X\}]\| \\
& + C_2^2 \frac{\|\widehat{f}_{n,\lambda} - f_0\|_\infty + \|\Delta f\|_\infty}{2} \|\Delta f\|^2 \\
\leq\ & C_2 d_n \frac{\|\widehat{f}_{n,\lambda} - f_0\|_\infty + \|\Delta f\|_\infty}{2n} \|\Delta f\| \cdot c_K \sqrt{n} h^{-1/2} B(h)\sqrt{M \log N} \\
& + C_2^2 \frac{\|\widehat{f}_{n,\lambda} - f_0\|_\infty + \|\Delta f\|_\infty}{2} \|\Delta f\|^2 \\
\leq\ & \frac{1}{2} C_2 M^{1/2} c_K^3 (2M_1 + 1)^3 h^{-3/2} \widetilde{r}_n^3 n^{-1/2} B(h)\sqrt{\log N} \\
& + \frac{1}{2} C_2^2 c_K (2M_1 + 1)^3 h^{-1/2} \widetilde{r}_n^3 \\
\leq\ & D_1(C_2, c_K, M, M_1) \times \widetilde{r}_n^3 (n^{-1/2} h^{-\frac{8m-1}{4m}}\sqrt{\log N} + h^{-1/2}) \\
\leq\ & D_1(C_2, c_K, M, M_1) \times \widetilde{r}_n^3 b_{n1},
\end{aligned}
$$

$$\tag{A.73}$$

where $D_1(C_2, c_K, M, M_1)$ is constant depending only on $C_2, c_K, M, M_1$.

We can use similar empirical processes techniques to handle $T_2$. Note that on $\mathcal{E}_n$ and for $\|f - f_0\| \leq \widetilde{r}_n$, using Assumption A1,

$$
\begin{aligned}
&|T_2(f)| \\
&= \frac{1}{2n} \left| \sum_{i=1}^n [\ddot{A}(f_0(X_i))(\Delta f)(X_i)^2 - E_{f_0}^X \{\ddot{A}(f_0(X))(\Delta f)(X)^2\}] \right| \\
&= \frac{1}{2n} \left| \left\langle \sum_{i=1}^n [\ddot{A}(f_0(X_i))(\Delta f)(X_i) K_{X_i} - E_{f_0}^X \{\ddot{A}(f_0(X))(\Delta f)(X) K_X\}], \Delta f \right\rangle \right| \\
&\leq \frac{1}{2n} \|\Delta f\| \\
&\qquad \times \| \sum_{i=1}^n [\ddot{A}(f_0(X_i))(\Delta f)(X_i) K_{X_i} - E_{f_0}^X \{\ddot{A}(f_0(X))(\Delta f)(X) K_X\}] \| \\
&= \frac{C_2 c_K h^{-1/2} d_n \|\Delta f\|}{2\sqrt{n}} \times \|Z_{n,f_0}^{(2)}(d_n^{-1} \Delta f)\| \\
&\leq \frac{C_2 c_K h^{-1/2} d_n \|\Delta f\|}{2\sqrt{n}} B(h) \sqrt{M \log N} \\
&\leq D_2(C_2, c_K, M, M_1) \times n^{-1/2} h^{-\frac{6m-1}{4m}} \widetilde{r}_n^2 \sqrt{\log N} \\
&\leq D_2(C_2, c_K, M, M_1) \times \widetilde{r}_n^2 b_{n2},
\end{aligned}
\tag{A.74}
$$

where $D_2(C_2, c_K, M, M_1)$ is constant depending only on $C_2, c_K, M_1, M$.

It follows that on $\mathcal{E}_n$,

$$
\begin{aligned}
J_1 &\geq \int_{\|f - f_0\| \leq \widetilde{r}_n} \exp(n I_n(f)) d\Pi(f) \\
&= \int_{\|f - f_0\| \leq \widetilde{r}_n} \exp\left( n T_1(f) + n T_2(f) - \frac{n}{2} \|f - \widehat{f}_{n,\lambda}\|^2 \right) d\Pi(f) \\
&\geq \exp\left( -[D_1(C_2, c_K, M, M_1) \widetilde{r}_n b_{n1} + D_2(C_2, c_K, M, M_1) b_{n2} \right. \\
&\qquad \left. + (M_1 + 1)^2/2] n \widetilde{r}_n^2 \right) \Pi(\|f - f_0\| \leq \widetilde{r}_n).
\end{aligned}
$$

To continue, we provide a lower bound for $\Pi(\|f - f_0\| \leq \widetilde{r}_n)$ using the same arguments as in (A.63). Note that $\lambda \leq \widetilde{r}_n^{\frac{4m}{2m+\beta-1}}$. Then it follows by Assumption A2, with $\varepsilon$ replaced by $\widetilde{r}_n$, that

$$
\Pi(\|f - f_0\| \leq \widetilde{r}_n) \geq c_1 \exp(-c_0 \widetilde{r}_n^{-\frac{2}{2m+\psi}}).
$$

Note that

$$
\widetilde{r}_n \geq (nh)^{-1/2} + h^{m+\frac{\psi}{2}} \geq 2n^{-\frac{2m+\psi}{2(2m+\psi+1)}},
$$

we get that

$$
n \widetilde{r}_n^{2 + \frac{2}{2m+\psi}} \geq n \left( 4n^{-\frac{2m+\psi}{2m+\psi+1}} \right)^{1 + \frac{1}{2m+\psi}} = 4.
$$

Therefore, $\widetilde{r}_n^{-\frac{2}{2m+\psi}} \leq n\widetilde{r}_n^2/4$, leading to

$$\Pi(\|f - f_0\| \leq \widetilde{r}_n) \geq c_1 \exp\left(-\frac{c_0}{4}n\widetilde{r}_n^2\right). \tag{A.75}$$

This implies by rate conditions $\widetilde{r}_n b_{n1} \leq 1$ and $b_{n2} \leq 1$ that, on $\mathcal{E}_n$,

$$
\begin{aligned}
J_1 &\geq c_1 \exp\left(-[D_1(C_2, c_K, M, M_1)\widetilde{r}_n b_{n1} + D_2(C_2, c_K, M, M_1)b_{n2}\right. \\
&\qquad \left. + (M_1 + 1)^2/2 + c_0/4]n\widetilde{r}_n^2\right) \\
&\geq c_1 \exp\left(-[D_1(C_2, c_K, M, M_1) + D_2(C_2, c_K, M, M_1)\right. \\
&\qquad \left. + (M_1 + 1)^2/2 + c_0/4]n\widetilde{r}_n^2\right).
\end{aligned}
$$

Next we handle $J_2$. The idea is similar to how we handle $J_1$ but with technical difference. Let $\Delta f = f - \widehat{f}_{n,\lambda}$. Note that $\widetilde{r}_n \leq r_n\sqrt{\log n}$, and hence, on $\mathcal{E}_n$, for any $f \in A_n$, i.e., $\|f - f_0\| \leq M_0 r_n$, we get that $\|\Delta f\| = \|\widehat{f}_{n,\lambda} - f\| \leq \|\widehat{f}_{n,\lambda} - f_0\| + \|f - f_0\| \leq M_1\widetilde{r}_n + M_0 r_n \leq (M_0 + M_1)r_n\sqrt{\log n}$. This implies that on $\mathcal{E}_n$, $\|\Delta f\|_\infty \leq c_K(M_0 + M_1)h^{-1/2}r_n\sqrt{\log n}$, where the last term by our rate assumption is $o(1)$, and hence, we can choose $n$ to be large enough so that $\|f_0\|_\infty + \|\widehat{f}_{n,\lambda} - f_0\|_\infty + \|\Delta f\|_\infty < C$. Let $d_{*n} = c_K(M_0 + M_1)h^{-1/2}r_n\sqrt{\log n}$. Then $d_{*n}^{-1}\Delta f \in \mathcal{G}$. Using previous similar arguments handling $T_1(f)$, i.e., (A.73) , we have that on $\mathcal{E}_n$, for any $f \in A_n$,

$$
\begin{aligned}
&|T_1(f)| \\
&\leq \frac{C_2 c_K(2M_1 + M_0)}{2n}h^{-1/2}r_n\sqrt{\log n} \\
&\quad \times \left(d_{*n}\|\sum_{i=1}^n[d_{*n}^{-1}(\Delta f)(X_i)K_{X_i} - E^X\{d_{*n}^{-1}(\Delta f)(X)K_X\}]\| \cdot \|\Delta f\| \right. \\
&\quad \left. + nE^X\{(\Delta f)(X)^2\}\right) \\
&\leq \frac{C_2 c_K(2M_1 + M_0)}{2n}h^{-1/2}r_n\sqrt{\log n} \\
&\quad \times (\sqrt{n}c_K h^{-1/2}d_{*n} \cdot (M_0 + M_1)r_n\sqrt{\log n} \cdot B(h)\sqrt{M \log N} \\
&\quad + nC_2[(M_0 + M_1)r_n\sqrt{\log n}]^2) \\
&= \frac{1}{2}C_2 c_K^3(2M_1 + M_0)^3 M^{1/2}h^{-3/2}r_n^3 n^{-1/2}B(h)(\log n)^2 \\
&\quad + \frac{1}{2}C_2^2 c_K(2M_1 + M_0)^3 h^{-1/2}r_n^3(\log n)^{3/2} \\
&\leq D_3(C_2, c_K, M, M_0, M_1) \times r_n^3\left(n^{-1/2}h^{-\frac{8m-1}{4m}}(\log n)^2 + h^{-1/2}(\log n)^{3/2}\right) \\
&= D_3(C_2, c_K, M, M_0, M_1) \times r_n^3 b_{n1} \leq D_3(C_2, c_K, M, M_0, M_1) \times \widetilde{r}_n^2,
\end{aligned}
$$

where $D_3(C_2, c_K, M, M_0, M_1)$ is constant depending only on $C_2, c_K, M, M_0, M_1$ and the last in-equality follows by rate condition $r_n^3 b_{n1} \leq \tilde{r}_n^2$. Likewise, on $\mathcal{E}_n$ and for any $f \in A_n$,

$$
\begin{aligned}
|T_2(f)| &\leq \frac{\|\Delta f\|}{2\sqrt{n}} C_2 c_K h^{-1/2} d_{*n} \cdot B(h) \sqrt{M \log n} \\
&\leq \frac{1}{2} C_2 c_K^2 (M_0 + M_1)^2 M^{1/2} n^{-1/2} h^{-1} r_n^2 B(h) (\log n)^{3/2} \\
&\leq D_4(C_2, c_K, M, M_0, M_1) \times n^{-1/2} r_n^2 h^{-\frac{6m-1}{4m}} (\log n)^{3/2} \\
&= D_4(C_2, c_K, M, M_0, M_1) \times r_n^2 b_{n2} \leq D_4(C_2, c_K, M, M_0, M_1) \times \tilde{r}_n^2,
\end{aligned}
$$

where $D_4(C_2, c_K, M, M_0, M_1)$ is constant only depending on $C_2, c_K, M, M_0, M_1$ and the last inequality follows by rate condition $r_n^2 b_{n2} \leq \tilde{r}_n^2$. It is easy to see that on $\mathcal{E}_n$ and for any $f \in A_n$,

$$
\|\widehat{f}_{n,\lambda} - f\| \geq \|f - f_0\| - \|\widehat{f}_{n,\lambda} - f_0\| \geq (2C' - M_1)\tilde{r}_n,
$$

leading to that

$$
\begin{aligned}
J_2 &\leq \\
&\exp\left(-\left(\frac{(2C' - M_1)^2}{2} - D_3(C_2, c_K, M, M_0, M_1) - D_4(C_2, c_K, M, M_0, M_1)\right) n\tilde{r}_n^2\right).
\end{aligned}
$$

Choose $C' > M_1$ to be large s.t.

$$
\begin{aligned}
\frac{(2C' - M_1)^2}{2} &\geq \\
&1 + D_1(C_2, c_K, M, M_1) + D_2(C_2, c_K, M, M_1) + D_3(C_2, c_K, M, M_0, M_1) \\
&+ D_4(C_2, c_K, M, M_0, M_1) + (M_1 + 1)^2/2 + c_3/4.
\end{aligned}
$$

Therefore, on $\mathcal{E}_n$,

$$
P(A_n | \mathbf{D}_n) \leq \frac{J_2}{J_1} \leq \exp(-n\tilde{r}_n^2).
$$

When $n$ becomes large s.t. $\exp(-n\tilde{r}_n^2) \leq \varepsilon_2/2$, we get that

$$
P_{f_0}^n(P(A_n | \mathbf{D}_n) \geq \varepsilon_2/2) \leq P_{f_0}^n(\mathcal{E}_n^c) \leq \varepsilon_1/2.
$$

This shows that

$$
P_{f_0}^n\left(P(\|f - f_0\| \geq 2C'\tilde{r}_n | \mathbf{D}_n) \geq \varepsilon_2\right) \leq \varepsilon_1.
$$

Proof is completed. □

*Verification of (3.7) $\Longrightarrow$ Rate Condition (**R**).* Consider two cases.

$$
\begin{aligned}
Case\,1: \quad &\max\left\{\frac{2}{6m + 3\psi - 1}, \frac{2m}{2m(4m + 2\psi - 3) + 1}, \frac{1}{4m}\right\} < a \leq \frac{1}{2m + \psi + 1} \\
Case\,2: \quad &\frac{1}{2m + \psi + 1} < a \leq \frac{1}{2m + 1}
\end{aligned}
$$

We only verify that Case 1 satisfies Rate Condition (**R**). The verification of Case 2 is similar. By Case 1, we have

$$r_n \asymp h^m \asymp n^{-ma}, \quad \widetilde{r}_n \asymp h^{m+\psi/2} \asymp n^{-(m+\psi/2)a}, \quad D_n \asymp (n^{-\frac{1}{2}+\frac{6m-1}{4m}a - ma} + n^{\frac{a}{2}-2ma}) \log n$$

$$b_{n1} \asymp n^{-\frac{1}{2}+\frac{8m-1}{4m}a}(\log n)^2 + n^{\frac{a}{2}}(\log n)^{3/2}, \quad b_{n2} \asymp n^{-\frac{1}{2}+\frac{6m-1}{4m}a}(\log n)^{3/2}$$

Then the following hold:

- $m > 3/2 \implies r_n = o(h^{3/2})$
- $h^{1/2}\log n = n^{-a/2}\log n = o(1)$
- $a < 1/(2m+1) \implies nh^{2m+1} \geq 1$
- $2m/(6m-1+2m\psi) > 1/(2m+\psi+1) \geq a$ and $\psi < m-1/2 \implies D_n = O(\widetilde{r}_n)$
- $4m^2 - 8m + 1 > 0, \psi > 0, 2m+\psi > 1 \implies \widetilde{r}_n b_{n1} \leq 1$
- $4m^2 + 2m\psi - 4m + 1 > 0 \implies b_{n2} \leq 1$
- $\psi < m-1/2, -1/2 + (8m-1)a/(4m) < a/2 < (m-\psi)a \implies r_n^3 b_{n1} \leq \widetilde{r}_n^2$
- $\psi < m-1/2, a < 1/(2m+\psi+1) \implies r_n^2 b_{n2} \leq \widetilde{r}_n^2$
- $a > 2/(6m+3\psi-1), a > 2m/[2m(4m+2\psi-3)+1] \implies n\widetilde{r}_n^3 b_{n1} = o(1)$
- $a > 2m/[2m(4m+2\psi-3)+1] \implies n\widetilde{r}_n^2 b_{n2} = o(1)$

Hence, Rate Condition (**R**) holds. □